\documentstyle[12pt,amsmath,amsfonts,amssymb]{article}
\newtheorem{theorem}{Theorem}[section]
\newtheorem{definition}[theorem]{Definition}

\newtheorem{proposition}[theorem]{Proposition}
\newtheorem{remark}[theorem]{Remark}
\newtheorem{lemma}[theorem]{Lemma}

\newtheorem{example}[theorem]{Example}
\newcommand {\Kc}      {{\mathcal K}}
\newcommand {\Mc}      {{\mathcal M}}

\newcommand {\Hc}      {{\mathcal H}}
\newcommand {\Ac}      {{\mathcal A}}

\newcommand {\Bc}      {{\mathcal B}}
\newcommand {\Cc}      {{\mathcal C}}
\newcommand {\Pc}      {{\mathcal P}}
\newcommand {\Ec}      {{\mathcal E}}
\newcommand {\Qc}      {{\mathcal Q}}
\newcommand {\Jc}      {{\mathcal J}}
\newcommand {\WP}      {W^1_p(\RN)}
\newcommand {\WPO}     {W^1_p(\Omega)}
\newcommand {\LOPO}    {L^1_p(\Omega)}
\newcommand {\LOP}     {L^1_p(\RN)}

\newcommand {\LKP}     {L^k_p(\RN)}
\newcommand {\BS}      {B^s_{p,q}(\RN)}
\newcommand {\Bspq}    {B^s_{p,q}}

\newcommand {\PK}      {\Pc_k}
\newcommand {\TW}      {{\mathbb W}(S)}
\newcommand {\TWM}     {{\mathbb W}_\ve(S)}
\newcommand {\TWU}     {{\mathbb W}(A)}
\newcommand {\TWE}     {{\mathbb W}(S:\ve)}
\newcommand {\TWB}[1]  {{\mathbb W}(S:#1)}
\newcommand {\tK}      {\widetilde{K}}

\newcommand {\tQ}      {\widetilde{Q}}

\newcommand {\tf}      {\tilde{f}}

\newcommand {\tmu}     {\tilde{\mu}}
\newcommand {\R}       {{\bf R}}

\newcommand {\RN}      {\R^n}

\newcommand {\ve}      {\varepsilon}
\newcommand {\tve}     {\tilde{\varepsilon}}
\newcommand {\dlt}     {\delta}
\newcommand {\Se}      {S_\varepsilon}

\newcommand {\emp}     {\emptyset}

\newcommand {\fc}      {f^{\curlywedge}}

\newcommand {\KRN}     {\Kc(\RN)}
\newcommand {\bc}      {\bar{c}}
\newcommand {\tpi}     {\widetilde{\pi}}
\newcommand {\Es}      {\Ex_{\delta,\bc,S}}
\newcommand {\EFQ}     {\Ec(f;Q)}
\newcommand {\AP}      {A_p}
\newcommand {\PRM}     {\Ac_{p,\alpha}}
\newcommand {\PRMA}[1] {\Ac_{p,#1}}
\newcommand {\SPR}     {{\mathbb A}_{p,\alpha}}

\newcommand {\TRS}     {\Jc_{p,\alpha}}
\newcommand {\omp}     {\omega_{1,p}}
\newcommand {\fq}      {(f\circ T)_Q}
\newcommand {\pcw}     {\pi^\curlywedge}
\newcommand {\WS}      {W^\sharp}
\newcommand {\Ect}     {\widetilde{\Ec}}
\newcommand {\intl}    {\int\limits}
\newcommand {\bQ}      {\overline{Q}}
\newcommand {\ot}      {\tilde{t}}

\newcommand {\rxy}     {\rho_{\alpha,S}(x,y)}
\newcommand {\rxx}     {\rho_{\alpha,S}(x,x')}
\newcommand {\rrxy}[1] {\rho_{\alpha,#1}(x,y)}
\newcommand {\card}    {\operatorname{card}}
\newcommand {\supp}    {\operatorname{supp}}
\newcommand {\Ex}      {\operatorname{Ext}}
\newcommand {\diam}    {\operatorname{diam}}
\newcommand {\dist}    {\operatorname{dist}}

\newcommand {\bx}      {\hspace{10mm}$\Box$}
\newcommand {\BX}      {\hspace{10mm}\Box}
\newcommand {\nn}      {\nonumber}
\newcommand {\rf}[1]   {(\ref{#1})}   
\newcommand {\reff}[1] {\ref{#1}}     
\newcommand {\SECT}[2] {\section*{\centerline{\normalsize
{\bf #1}}} \setcounter{section}{#2}
\setcounter{theorem}{0}\setcounter{equation}{0}}
\newcommand {\SECTLONG}[3]
{\section*{\centerline{\normalsize {\bf #1}}
\centerline{\normalsize {\bf #2}}} \setcounter{section}{#3}
\setcounter{theorem}{0}\setcounter{equation}{0}}
\newcommand{\lbl}[1]        {\label{#1}}            
\newcommand{\be}            {\begin{eqnarray}}
\newcommand{\bel}[1]        {\begin{eqnarray} \label{#1}}
\newcommand{\ee}            {\end{eqnarray}}
\begin{document}
\medskip
\centerline{\large{\bf Sobolev $W^1_p$-spaces on closed
subsets of $\RN$}}
\vspace*{10mm} \centerline{By  {\it Pavel Shvartsman}}
\vspace*{12 mm}
\renewcommand{\thefootnote}{ }
\footnotetext[1]{{\it\hspace{-6mm}Math Subject
Classification} 46E35}
\begin{abstract} For each $p>n$ we use local oscillations
and doubling measures to give intrinsic characterizations
of the restriction of the Sobolev space
$W_{p}^{1}({\R}^{n})$ to an arbitrary closed subset of
${\R}^{n}$.
\end{abstract}
\renewcommand{\thefootnote}{\arabic{footnote}}
\setcounter{footnote}{0}
\SECT{1. Introduction.}{1}
\indent
{\bf 1.1 Main definitions and results.} Let $S$ be an
arbitrary closed subset of ${\R}^{n}$. In this paper we
describe the restrictions to $S$ of the functions in the
Sobolev space $W_{p}^{1}({\bf R}^{n})$ whenever $p>n$.
\par We recall that, for each choice of the open set
$\Omega \subset {\R}^{n}$ and of the exponent $p\in
[1,\infty]$, the Sobolev space $W_{p}^{1}(\Omega ) $
consists of all (equivalence classes of) real valued
functions $f\in L_{p}(\Omega )$ whose first order
distributional partial derivatives on $\Omega $ belong to
$L_{p}(\Omega )$. $W_{p}^{1}(\Omega )$ is normed by
$$ \Vert f\Vert _{W_{p}^{1}(\Omega )}:=\Vert f\Vert
_{L_{p}(\Omega )}+\Vert \nabla f\Vert _{L_{p}(\Omega )}. $$
\par We also recall that $L_{p}^{1}(\Omega )$ is the
corresponding homogeneous Sobolev space, defined by the
finiteness of the seminorm
$\Vert f\Vert _{L_{p}^{1}(\Omega )} :=\Vert \nabla f\Vert
_{L_{p}(\Omega )}.$
\par There is an extensive literature devoted to describing
the restrictions of Sobolev functions to different classes
of subsets of ${\bf R}^{n}$. (We refer the reader to
\cite{A,BIN,BK,Bur,FJ,GV1,GV2,Jn,JW,K,M,MP,St,S4,T3} and
references therein for numerous results in this direction
and techniques for obtaining them.)
\par In this paper we consider restrictions of Sobolev
functions to an {\it arbitrary} closed subset $S$ of
${\R}^{n}$ and ways of extending such functions back to
Sobolev functions on all of ${\R}^{n}$. In the case where
$p=\infty $, the Sobolev space $W_{\infty }^{1}({\R}^{n})$
can be identified with the space $\operatorname{Lip}(\RN)$
of Lipschitz functions on ${\R}^{n}$ and it is known that
the restriction $\operatorname{Lip}(\RN)|_{S}$ coincides
with the space $\operatorname{Lip}(S)$ of Lipschitz
functions on $S$ and that, furthermore, the classical
Whitney extension operator linearly and continuously maps
the space $\operatorname{Lip}(S)$ into the space
$\operatorname{Lip}(\RN)$ (see e.g., \cite{St}, Chapter 6).
\par Here we will show that, perhaps surprisingly, an
analogous result also holds for all $p$ in the range
$n<p<\infty $, namely, that {\it very same classical linear
Whitney extension operator also provides an extension to
$W_{p}^{1}({\R}^{n})$, even an almost optimal extension, of
each function on $S$ which is the restriction to $S$ of a
function in $W_{p}^{1}({\R}^{n})$.} We also provide various
intrinsic characterizations of the restriction of
$W_{p}^{1}({\R}^{n})$ to the given set $S$. To the best of
our knowledge, these are the first results of this type to
have been obtained for the range $n<p<\infty $ without the
imposition of some extra conditions on the set $S$.
\par Recall that, when $p>n$, it follows from the Sobolev
embedding theorem that every function $f\in L_{p}^{1}({\bf
R}^{n})$ coincides almost everywhere with a function
satisfying the H\"{o}lder condition of order $\alpha
:=1-\frac{n}{p}$. I.e., after possibly modifying $f$ on a
set of Lebesgue measure zero, we have
\bel{S-Lip} |f(x)-f(y)|\le C(n,p) \Vert
f\Vert_{L_{p}^{1}(\RN)} \Vert x-y\Vert
^{1-\frac{n}{p}}~~\text{ for all }~~x,y\in \RN. \ee
\par This fact enables us {\it to identify each element }
$f\in W_{p}^{1}(\RN)$, {\it $p>n$, with its unique
continuous representative}. This identification, in
particular, implies that $f$ has a well defined restriction
to any given subset of ${\bf R}^{n}$. Furthermore it means
that there is no loss of generality in only considering
closed subsets of ${\bf R}^{n}$. Thus, given an arbitrary
closed subset $S$ of ${\bf R}^{n}$, we define the trace
space $W_{p}^{1}({\bf R}^{n})|_{S},$ $p>n,$ of all
restrictions of $W_{p}^{1}({\bf R}^{n})$-functions to $S$
by
$$ W_{p}^{1}({\bf R}^{n})|_{S}:=\{f:S\to {\bf R}:{\rm
there~is~a~continuous}\ \ F\in W_{p}^{1}({\bf R}^{n})\ \
{\rm such~that}\ \ F|_{S}=f\}. $$
$W_{p}^{1}({\bf R}^{n})|_{S}$ is equipped with the standard
quotient space norm
$$ \Vert f\Vert _{W_{p}^{1}({\bf R}^{n})|_{S}}
:=\inf \{\Vert F\Vert _{W_{p}^{1}(%
{\bf R}^{n})}:F\in W_{p}^{1}({\bf R}^{n})~
{\rm and~continuous,~~}F|_{S}=f\}.
$$
The trace space $L_{p}^{1}({\bf R}^{n})|_{S}$ is defined in
an analogous way.
\par Our first main result, Theorem \reff{EXT-SIMP}, gives
an intrinsic characterization of the space
$L_{p}^{1}(\RN)|_{S}$. Before formulating this result we
need to fix some notation: Throughout this paper, the
terminology ``cube'' will mean a closed cube in ${\bf
R}^{n}$ whose sides are parallel to the coordinate axes. We
let $Q(x,r)$ denote the cube in $\RN$ centered at $x$ with
side length $2r$. Given $\lambda >0$ and a cube $Q$ we let
$\lambda Q$ denote the dilation of $Q$ with respect to its
center by a factor of $\lambda $. (Thus $\lambda
Q(x,r)=Q(x,\lambda r)$.) The Lebesgue measure of a
measurable set $A\subset \RN$ will be denoted by
$\left|A\right|$.
\begin{theorem} \lbl{EXT-SIMP} A function $f$ defined on a
closed set $S\subset {\bf R}^{n}$ can be extended to a
(continuous) function $F\in L_{p}^{1}({\bf
R}^{n}),n<p<\infty ,$ if and only if there exists a
constant $\lambda >0$ such that for every finite family
$\{Q_{i}:i=1,...,m\}$ of pairwise disjoint cubes in ${\bf
R}^{n}$ and every choice of points
\bel{11Q} x_{i},y_{i}\in (11Q_{i})\cap S \ee
the inequality
\bel{CONSTR} \sum_{i=1}^{m}\frac{|f(x_{i})-f(y_{i})|^{p}}
{(\operatorname{diam}Q_{i})^{p-n}}\le \lambda \ee
holds. Moreover,
\bel{EQVF} \Vert f\Vert _{L_{p}^{1}({\bf R}^{n})|_{S}} \sim
\inf \lambda ^{\frac{1}{p}} \ee
with constants of equivalence depending only on $n$ and
$p$.
\end{theorem}
\par The special case of this result where $S={\bf R}^{n}$
has been treated by Yu. Brudnyi \cite{Br1}. He proved that
a function $f$ is in $L_{p}^{1}({\bf R }^{n})$, $p>n$, if
and only if the inequality \rf{CONSTR} holds for every
finite family $\{Q_{i}\}$ of disjoint cubes and all
$x_{i},y_{i}\in Q_{i}$.
\par  In particular, Brudnyi's result shows that, for
$S={\bf R}^{n}$, the condition \rf{11Q} can be replaced by
$x_{i},y_{i}\in Q_{i}$. We do not know whether the same is
true in the case of an arbitrary subset $S$. We do know
that a more elaborate variant of our proof of Theorem
\reff{EXT-SIMP} enables the factor $11$ in \rf{11Q} to be
decreased to $3+\varepsilon $ where $\varepsilon >0$ is
arbitrary. (Of course, in this case the constants of
equivalence in \rf{EQVF} will also depend on $\varepsilon
$).
\par We now formulate the second main result of the paper,
Theorem \reff{EXT-NORM}, which describes the restrictions
of $W_{p}^{1}({\bf R}^{n})$-functions to $S$.
\begin{theorem} \lbl{EXT-NORM} Let $S$ be an arbitrary
closed subset of ${\bf R}^{n}$. Let $\varepsilon $ be an
arbitrary positive number and let $S_{\varepsilon }$ denote
the $\varepsilon $-neighborhood of $S$. Fix a number
$\theta \ge 1$ and let $T:S_{\varepsilon }\to S$ be a
measurable mapping such that
\bel{TR} \Vert T(x)-x\Vert \le \theta
\operatorname{dist}(x,S)~~ \text{ for every }~~x\in
S_{\varepsilon }. \ee
The function $f:S\rightarrow {\bf R}$ is an element of
$W_{p}^{1}({\bf R}^{n})|_{S}$ if and only if $f\circ T\in
L_{p}(S_{\varepsilon })$ and there exists a constant
$\lambda >0$ such that the inequality \rf{CONSTR} holds for
every finite family $\{Q_{i}:i=1,...,m\}$ of disjoint cubes
contained in $S_{\varepsilon }$, and every choice of points
$x_{i},y_{i}\in (\eta \,Q_{i})\cap S$, where
$\eta:=10\,\theta +1$.
Moreover,
\bel{TRN} \Vert f\Vert _{W_{p}^{1}({\bf R}^{n})|_{S}}\sim
\Vert f\circ T\Vert _{L_{p}(S_{\varepsilon })}+\inf \lambda
^{\frac{1}{p}} \ee
with constants of equivalence depending only on
$n,p,\varepsilon $ and $\theta $.
\end{theorem}
\begin{remark} \lbl{T-NORM}{\em The mapping $T$ appearing
in this theorem must of course satisfy $T(x)=x$ for all
$x\in S$.
\par In fact, a simple example of a mapping $T$ satisfying
\rf{TR} can be easily constructed with the help of a
Whitney decomposition $W(S)=\{Q_{k}:~k\in I\}$ of the open
set ${\bf R}^{n}\setminus S$. (Recall that $W(S)$ is a
covering of ${\bf R}^{n}\setminus S$ by non-overlapping
cubes such that $\operatorname{diam}Q_{k}\sim
\operatorname{dist}(Q_k,S),k\in I.$) Let $a_{Q}$ be a point
in $S$ which is nearest to $Q$ and let
\bel{TWQ} T(x):=\left\{ \begin{array}{ll}
x, & x\in S, \\
&  \\
a_{Q}, & x\in int(Q),~~Q\in W(S), \end{array} \right. \ee
where, for each set $A\subset {\bf R}^{n}$,  $int(A)$
denotes the interior of $A$. Clearly, $T$ is defined a.e.
on ${\bf R}^{n}$ and satisfies the inequality \rf{TR} (with
a constant $\theta $ depending only on parameters of the
Whitney decomposition $W(S)$). If we use this particular
mapping $T$ then \rf{TRN} yields the following formula for
the trace norm of a function:
$$ \Vert f\Vert _{W_{p}^{1}({\bf R}^{n})|_{S}} \sim \Vert
f\Vert _{L_{p}(S)}+\left\{ \sum_{Q\in
W(S),\,\operatorname{diam}Q \le \varepsilon
}|f(a_{Q})|^{p}\,|Q|\right\} ^{\frac{1}{p}}+ \inf \lambda
^{\frac{1}{p}} $$
(with constants of equivalence depending only on
$n,p,\varepsilon,$ and parameters of $W(S))$.} \end{remark}
\par Our next result provides a different kind of
description of trace spaces. It is expressed in terms of a
certain kind of maximal function which is a variant of
$f_{p}^{\sharp }$, the familiar {\it sharp maximal function
in} $L_{p}$. Let us first recall that, for each $p\in
[1,\infty )$ and each $f\in L_{p,loc}({\bf R}^{n})$, one
defines $f_{p}^{\sharp}$ by
$$
f_{p}^{\sharp }(x):=\sup_{r>0}\frac{1}{r}
\left( \frac{1}{|Q(x,r)|}%
\int\limits_{Q(x,r)}|f(y)-f_{Q(x,r)}|^{p}\,dy\right)
^{\frac{1}{p}}.
$$
Here $f_{Q}:=|Q|^{-1}\int_{Q}f\,dx$ denotes the average of
$f$ on the cube $Q$. For $p=\infty $ we put
$$ f_{\infty }^{\sharp }(x):=\sup
\{|f(y)-f(z)|/r:~r>0,\,y,z\in Q(x,r)\}. $$
In \cite{C1} Calder\'{o}n proved that, for $1<p\le \infty$,
a function $f$ is in $W_{p}^{1}({\bf R}^{n})$ if and only
if $f$ and $f_{1}^{\sharp }$ are both in $L_{p}({\bf
R}^{n})$, and that, furthermore,
\bel{W-MF1}
\Vert f\Vert _{W_{p}^{1}({\bf R}^{n})}\sim
\Vert f\Vert _{L_{p}({\bf R}^{n})}+\Vert f_{1}^{\sharp
}\Vert _{L_{p}({\bf R}^{n})}
\ee
with constants of equivalence depending only on $p$ and
$n$. See also \cite{CS}. We observe that the Sobolev
imbedding theorem (for $p>n$) allows us to replace
$f_{1}^{\sharp }$ in \rf{W-MF1}) with the (bigger) sharp
maximal function $f_{\infty }^{\sharp }$. Thus, for $p>n$
and for every $f\in W_{p}^{1}({\bf R}^{n})$, we have
$$ \Vert f\Vert _{W_{p}^{1}({\bf R}^{n})}
\sim \Vert f\Vert _{L_{p}({\bf R}%
^{n})}+\Vert f_{\infty }^{\sharp }\Vert _{L_{p}
({\bf R}^{n})}.
$$
This equivalence motivates us to introduce the new variant
of the sharp maximal function $f_{\infty }^{\sharp }$ that
we need here. It is defined with respect to an {\it
arbitrary} closed set $S\subset {\bf R}^{n}$ for functions
$f:S\rightarrow {\bf R}$. It is given by
\[ f_{\infty ,S}^{\sharp }(x):=\sup
\{|f(y)-f(z)|/r:~r>0,\,y,z\in Q(x,r)\cap S\},~~~~x\in {\bf
R}^{n}. \]
\begin{theorem} \lbl{MF-C} Let $S\subset {\bf R}^{n}$ be a
closed set and let $p>n$.
\par (i). A function $f$ defined on $S$ can be extended to
a (continuous) function $F\in L_{p}^{1}({\bf R}^{n})$ if
and only if $f_{\infty ,S}^{\sharp }\in L_{p}({\bf
R}^{n}).$ Moreover,
\bel{FirstNewNumber} \Vert f\Vert _{L_{p}^{1}({\bf
R}^{n})|_{S}}\sim \Vert f_{\infty ,S}^{\sharp }\Vert
_{L_{p}({\bf R}^{n})}. \ee
\par (ii). Given $\varepsilon >0$ and $\theta \ge 1$, let
$S_{\varepsilon }$ be the $\varepsilon $-neighborhood of
$S$, and let $T:S_{\varepsilon }\to S$ be a measurable
mapping satisfying the condition \rf{TR}. A function $f$
defined on $S$ is an element of $W_{p}^{1}({\bf
R}^{n})|_{S}$ if and only if $f\circ T$ and $f_{\infty
,S}^{\sharp }$ are both in $L_{p}(S_{\varepsilon }) $.
Furthermore,
\bel{SecondNewNumber} \Vert f\Vert _{W_{p}^{1}({\bf
R}^{n})|_{S}}\sim \Vert f\circ T\Vert
_{L_{p}(S_{\ve})}+\Vert f_{\infty ,S}^{\sharp }\Vert
_{L_{p}(S_{\ve})}. \ee
\par The constants of equivalence in \rf{FirstNewNumber}
depend only on $n$ and $p$, and in \rf{SecondNewNumber}
they only depend on $n$, $p$, $\varepsilon $ and $\theta $.
\end{theorem}
\par As already mentioned above, the classical Whitney
extension operator provides an almost optimal extension of
a function $f$ defined on $S$ to a function $F$ in
$W_{p}^{1}({\bf R}^{n})$ or in $L_{p}^{1}({\bf R}^{n})$,
whenever $p>n$. Since this extension operator is {\it
linear}, we obtain the following
\begin{theorem} \lbl{LINEXT} For every closed subset
$S\subset {\bf R}^{n}$ and every $p>n$ there exists a
linear extension operator which maps the trace space $
W_{p}^{1}({\bf R}^{n})|_{S}$ continuously into
$W_{p}^{1}({\bf R}^{n})$ and also maps the trace space
$L_{p}^{1}({\bf R}^{n})|_{S}$ continuously into
$L_{p}^{1}({\bf R}^{n})$. Its operator norms are bounded by
constants depending only on $n$ and $p$.
\end{theorem}
\medskip
{\bf 1.2 Extensions of Sobolev functions and doubling
measures.} In this subsection we present a series of
results which describe the restrictions of $W_{p}^{1}({\bf
R}^{n})$-functions to a closed set $S$ via certain doubling
measures supported on $S$. As before, we assume that $p>n$.
Our approach here is inspired by a paper of Jonsson
\cite{J} devoted to characterization of traces of Besov
spaces to arbitrary closed subsets of ${\bf R}^{n}$, see
Section 6.
\par Let $\mu $ be a positive regular Borel measure on
${\bf R}^{n}$ such that
$$ \operatorname{supp}\mu =S\,. $$
Following \cite{VK}, we say that $\mu $ satisfies the
condition $(D_{n})$, if there exists a constant
$C=C_{\mu }>0$ such that, for all $x\in S$, $%
0<r\le 1$ and $1\le k\le 1/r$, the following inequality
\bel{DN}
 \mu (Q(x,kr))\le C_{\mu }k^{n}\mu (Q(x,r))
\ee
holds. The condition $(D_{n})$ is a generalization of the
familiar {\it doubling condition}:
\bel{DOUBC} \mu (Q(x,2r))\le c_{\mu }\mu (Q(x,r)),~~~x\in
S,~0<r\le 1/2. \ee
We refer to the constant $c_{\mu }$ in \rf{DOUBC} as a {\it
doubling constant}.
\par We also assume that the measure $\mu $ satisfies the
following condition: there exist constants $C_{\mu}^{\prime
},C_{\mu }^{\prime \prime }>0$ such that
\bel{B1} C_{\mu }^{\prime }\le \mu (Q(x,1))\le C_{\mu
}^{\prime \prime },~~~\text{for all }x\in S. \ee
\par Dynkin \cite{D1,D2} conjectured that {\it every
compact set } $S\subset {\bf R}^{n}$ {\it carries a
non-trivial doubling measure $\mu $}. He constructed a
doubling measure on every compact set $S\subset {\bf
R}^{1}$ which satisfies a certain ``porosity'' condition
(see the so-called ``the ball condition'' mentioned below).
Dynkin's conjecture was subsequently proved by Volberg and
Konyagin \cite{VK}. (See also Wu \cite{Wu}.) Moreover,
Volberg and Konyagin showed that every compact set in ${\bf
R}^{n}$ carries a non-trivial measure $\mu $ satisfying the
$(D_{n})$-condition \rf{DN}. Using their argument,
Luukkainen and Saksman \cite{LS} extended this result to
all closed subsets of ${\bf R}^{n}$. Jonsson \cite{J}
showed that it may be also assumed that $\mu $ satisfies
\rf{B1}.
\par  Throughout this paper $\sigma$ will denote a measure
{\it supported on the boundary $\partial S$ of the set
$S$.} It will usually be assumed that $\sigma$ satisfies
the generalized doubling condition $(D_{n})$ and also the
condition \rf{B1}. Thus $\operatorname{supp}\sigma
=\partial S$, and there exist positive constants $C_{\sigma
},C_{\sigma }^{\prime },C_{\sigma }^{\prime \prime },$ such
that
\bel{DNSG} \sigma (Q(x,kr))\le C_{\sigma }k^{n}\sigma
(Q(x,r)),~~~~~x\in \partial S,~1\le k\le 1/r, \ee
and
\bel{B1SIGMA} C_{\sigma }^{\prime }\le \sigma (Q(x,1))\le
C_{\sigma }^{\prime \prime },~~~~x\in \partial S. \ee
\par As preparation for formulating our the next main
result, Theorem \reff {EXT-DIFF-15}, we need to define a
special kind of ``distance'' between points of some subset
of ${\bf R}^{n}$. Given a set $A\subset {\bf R}^{n}$ and
points $x,y\in A$, let
\bel{DX15} \rho _{A}(x,y):=\inf \{\operatorname{diam}Q:
~Q\text{ is a cube, }~Q\ni x,y,~%
\tfrac{1}{15}\,Q\subset {\bf R}^{n}\setminus A\}. \ee
Thus we have $\rho _{A}(x,y):=+\infty $ whenever there does
not exist any cube $Q$ for which $x,y\in Q$ and
$\tfrac{1}{15}\,Q\subset {\bf R} ^{n}\setminus A$. Clearly,
$\Vert x-y\Vert \le \rho _{A}(x,y)$ for every $%
x,y\in A$.
\begin{theorem} \lbl{EXT-DIFF-15} Let $S$ be an arbitrary
closed subset of ${\bf R}^{n}$. Let $\mu$ be a measure with
$\supp \mu=S$ satisfying \rf{DN} and \rf{B1}, and let
$\sigma$ be a measure supported on $\partial S$ which
satisfies \rf{DNSG} and \rf{B1SIGMA}. Then, for every $p\in
(n,\infty )$ and every function $f\in C(S)$, we have
\begin{eqnarray} \Vert f\Vert _{W_{p}^{1}({\bf
R}^{n})|_{S}} &\sim &\Vert f\Vert _{L_{p}(\mu
)}+\sup_{0<t\le 1} \left( ~\iint\limits_{\Vert x-y\Vert
<t}\frac{ |f(x)-f(y)|^{p}}{t^{p-n}\mu (Q(x,t))^{2}} \,d\mu
(x)\,d\mu (y)\right) ^{\frac{
1}{p}}  \nonumber \\
&+&\Vert f\Vert _{L_{p}(\sigma )}+
\left( ~\iint\limits_{\rho _{\partial
S}(x,y)<1}\frac{|f(x)-f(y)|^{p}}
{\rho _{\partial S}(x,y)^{p-n}\sigma
(Q(x,\rho _{\partial S}(x,y)))^{2}}
\,d\sigma (x)\,d\sigma (y)\right) ^{\frac{
1}{p}}  \nonumber
\end{eqnarray}
with constants of equivalence depending only on $n,p$, and
the constants $C_{\mu },C_{\mu }^{\prime },C_{\mu }^{\prime
\prime }$ and $C_{\sigma },C_{\sigma }^{\prime },C_{\sigma
}^{\prime \prime }$.
\end{theorem}
\begin{remark} \lbl{SB-TRN}{\em Let us make some comments
about the expressions
$$ I_{1}(f;S):=\sup_{0<t\le 1} \left( ~\iint\limits_{\Vert
x-y\Vert <t}\frac{ |f(x)-f(y)|^{p}}{t^{p-n}\mu
(Q(x,t))^{2}}\,d\mu (x) \,d\mu (y)\right) ^{\frac{1}{p}} $$
and
$$ I_{2}(f;S):=\left( ~\iint\limits_{\rho _{\partial
S}(x,y)<1} \frac{|f(x)-f(y)|^{p}}{\rho _{\partial
S}(x,y)^{p-n} \sigma (Q(x,\rho _{\partial
S}(x,y)))^{2}}\,d\sigma (x)\,d\sigma (y)
\right)^{\frac{1}{p}} $$
which appear in the statement of the preceding theorem.
In particular $I_{1}(f;S)$ can be rewritten in the form
\bel{I1S} I_{1}(f;S)=\sup_{0<t\le 1}{\bf
w}_{1}(f;t)_{L_{p}(S;\mu )}/t, \ee
where
$$ {\bf w}_{1}(f;t)_{L_{p}(S;\mu )} :=\left(
\frac{1}{t^{n}}\,\iint\limits_ {\Vert x-y\Vert
<t}|f(x)-f(y)|^{p} \left( \frac{t^{n}}{\mu (Q(x,t))}\right)
^{2}\,d\mu (x)\,d\mu (y)\right) ^{\frac{1}{p}},~~~~t>0. $$
We will refer to the function ${\bf
w}_{1}(f;\cdot)_{L_{p}(S;\mu )}$  as {\it the averaged
modulus of smoothness}  of $f$ in $L_{p}(S;\mu)$. This name
is motivated by the following example: If we choose $S={\bf
R}^{n}$ and let $\mu $ be Lebesgue measure, then the
function ${\bf w}_{1}(f;\cdot )_{L_{p}({\bf R}^{n};\mu )}$
coincides with the function
$$ {\bf w}_{1}(f;t)_{L_{p}}:=\left(
\frac{1}{t^{n}}\,\int\limits_{\Vert h\Vert
<t}\int\limits_{{\bf
R}^{n}}|f(x+h)-f(x)|^{p}\,dx\,dh\right) ^{\frac{1}{p}}, $$
which is known in the literature as the averaged modulus of
smoothness in $L_{p}$. It is also known (see, e.g.
\cite{DL}) that ${\bf w}_{1}(f;\cdot )_{L_{p}}\sim \omega
(f;\cdot )_{L_{p}}$ where $\omega (f;\cdot )_{L_{p}}$ is
the classical modulus of smoothness in $L_{p}$:
\bel{DOM} \omega (f;t)_{L_{p}}:=\sup_{\Vert h\Vert <t}
\left( ~\int\limits_{{\bf R}
^{n}}|f(x+h)-f(x)|^{p}\,dx\right) ^{\frac{1}{p}}. \ee
It is known that
\bel{L1RN}
\sup_{t>0}\omega (f;t)_{L_{p}}/t\sim
\Vert f\Vert _{L_{p}^{1}({\bf R}%
^{n})},
\ee
see, e.g. \cite{Z}, p. 45, so that
$$
\|f\|_{\WP}
\sim \|f\|_{L_p(\RN)}+
\sup_{0<t\le 1}\omega (f;t)_{L_{p}}/t\,.
$$
Consequently, we have
\bel{FourthNewNumber}
\Vert f\Vert _{W_{p}^{1}({\bf
R}^{n})}\sim \Vert f\Vert _{L_{p}({\bf R}
^{n})}+\sup_{0<t\le 1}{\bf w}_{1}(f;t)_{L_{p}}/t\,.
\ee
In view of \rf{FourthNewNumber} and \rf{I1S}, the quantity
$\Vert f\Vert _{L_{p}(\mu )}+I_{1}(f;S)$ can be interpreted
as a certain Sobolev-type norm on $S$. This motivates us to
call $I_{1}(f;S)$ the {\it Sobolev part} of the trace norm
$\Vert f\Vert_{W_{p}^{1}({\bf R}^{n})|_{S}}$.
\par We now discuss the second expression $I_{2}(f;S)$.
This quantity depends only on the values of $f$ on
$\partial S$, the measure $\sigma $ supported on $\partial
S$ and the function $\rho _{\partial S}$ defined on
$\partial S$. So we can consider $I_{2}(f;S)$ as a ``
characteristic'' of the trace $f|_{S}$ which controls its
behavior on the boundary of the set $S$. Let us consider
the example when $S={\bf R}^{d}$ for some $d$, $\,0<d<n$.
Then, clearly, $S\sim \partial S$ and $\rho _{\partial
S}(x,y)\sim\Vert x-y\Vert .$ Suppose that, for this choice
of $S$, we take $\sigma $ to be Lebesgue measure on ${\bf
R}^{d}$ so that $\sigma (Q(x,t))\sim t^{d}$ for every $x\in
S$ and every $t>0$. It can be readily shown that in this
case
\begin{eqnarray} \|f\|_{L_p(\sigma)}+I_{2}(f;S) &\sim
&\|f\|_{L_p(\R^d)}+
\left(~\iint\limits_{\Vert x-y\Vert <1}\frac{
|f(x)-f(y)|^{p}}{\Vert x-y\Vert ^{p-n} \Vert x-y\Vert
^{2d}}\,d\sigma
(x)\,d\sigma (y)\right) ^{\frac{1}{p}}  \nonumber \\
&=&\|f\|_{L_p(\R^d)}+
\left( ~\iint\limits_{\Vert x-y\Vert <1}
\frac{|f(x)-f(y)|^{p}}{\Vert
x-y\Vert ^{d+\left( 1-\frac{n-d}{p}\right) p}}
\,d\sigma (x)\,d\sigma
(y)\right) ^{\frac{1}{p}}.  \nonumber
\end{eqnarray}
In fact this last quantity is the norm in a certain Besov
space of functions defined on ${\bf R}^{d}$. We recall that
the Besov space $B_{pp}^{s}({\bf R}^{d})$, for $0<s<1$,\ is
defined by the finiteness of the norm
$$ \Vert f\Vert _{B_{p,p}^{s}({\bf R}^{d})}
:=\Vert f\Vert _{L_{p}({\bf R}%
^{d})}+\left( ~\iint\limits_
{\Vert x-y\Vert <1}\frac{|f(x)-f(y)|^{p}}{\Vert
x-y\Vert ^{d+sp}}\,d\sigma (x)\,d\sigma (y)\right)
 ^{\frac{1}{p}}.
$$
(For the general theory of Besov spaces we refer the reader
to the monographs \cite{BIN,N,T2,T3} and references
therein.)
\par Thus for $S={\bf R}^{d}$ the quantity $\Vert f\Vert
_{L_{p}(\sigma )}+I_{2}(f;S)$ is an equivalent norm for the
Besov space $B_{p,p}^{1- \frac{n-d}{p}}({\bf R}^{d})$. This
suggests that we should interpret $ \Vert f\Vert
_{L_{p}(\sigma )}+I_{2}(f;S)$ for an {\it arbitrary} set $
S\subset {\bf R}^{n}$ as a certain ``Besov-type'' norm on $
\partial S$. This example also motivates us to call $\Vert
f\Vert _{L_{p}(\sigma )}+I_{2}(f;S)$ the {\it Besov part}
of the trace norm $\Vert f\Vert _{W_{p}^{1}({\bf
R}^{n})|_{S}}$.
\par To summarize, these observations and examples give us
some idea of the structure of the trace space
$W_{p}^{1}({\bf R}^{n})|_{S}$ for $p>n$. In particular,
they show that, in general, $W_{p}^{1}({\bf R}^{n})|_{S}$
is a certain ``Sobolev-type''  space of functions on $S$
whose values on the boundary belong to a certain
``Besov-type'' space of functions defined on $\partial S$.}
\end{remark}
\par The next theorem provides additional tools for the
calculation of the trace norm $\Vert \cdot \Vert
_{W_{p}^{1}({\bf R}^{n})|_{S}}$. In particular it offers us
the options of using either a suitable measure $\sigma $
supported on $\partial S$ or a suitable measure $\mu $
supported on $S$, for this calculation. It is convenient to
have this flexibility since, on the one hand, there are
many choices of $S$ for which the measure $\sigma $ can be
easily determined, whereas it may be a very non-trivial
task to specify a suitable measure $\mu $, and, on the
other hand, there are sets $S$ for which it is obvious how
to choose $\mu $ but very difficult to specify $\sigma $.
\begin{theorem} \lbl{EXT-DIFF-S15} Let $p\in (n,\infty )$.
\par (i). Let $\sigma $ be a measure supported on $\partial
S$ which satisfies conditions \rf{DNSG} and \rf{B1SIGMA},
and let $\Omega :=int(S)$ be the interior of $S$. Then for
every $f\in C(S)$ the following equivalence
\begin{eqnarray} \Vert f\Vert _{W_{p}^{1}({\bf
R}^{n})|_{S}} &\sim &\Vert f|\,_{\Omega }\Vert
_{W_{p}^{1}(\Omega )}+\Vert f\Vert _{L_{p}(\sigma
)}+\sup_{0<t\le 1}\left( \,\,\,\iint\limits_{\Vert x-y\Vert
<t}\frac{|f(x)-f(y)|^{p}}{t^{p-n}\sigma
(Q(x,t))^{2}}\,d\sigma (x)d\sigma (y)\right)
^{\frac{1}{p}}  \nonumber \\
&+&\left( ~\iint\limits_{\rho _{\partial S}(x,y)<1}
\frac{|f(x)-f(y)|^{p}}{%
\rho _{\partial S}(x,y)^{p-n}
\sigma (Q(x,\rho _{\partial S}(x,y)))^{2}}%
\,d\sigma (x)\,d\sigma (y)\right) ^{\frac{1}{p}}  \nonumber
\end{eqnarray}
holds.
\par (ii). For every function $f\in C(S)$ and every measure
$\mu $ supported on $S$ and satisfying \rf{DN} and \rf{B1},
we have
\begin{eqnarray} \Vert f\Vert _{W_{p}^{1}({\bf
R}^{n})|_{S}} &\sim &\Vert f\Vert _{L_{p}(\mu
)}+\sup_{0<t\le 1} \left( ~\iint\limits_{\Vert x-y\Vert
<t}\frac{ |f(x)-f(y)|^{p}}{t^{p-n}\mu (Q(x,t))^{2}}\, d\mu
(x)\,d\mu (y)\right) ^{\frac{
1}{p}}  \nonumber \\
&+&\left( ~\iint\limits_{\rho _{S}(x,y)<1}
\frac{|f(x)-f(y)|^{p}}{\rho
_{S}(x,y)^{p-n}\mu (Q(x,\rho _{S}(x,y)))^{2}}
\,d\mu (x)\,d\mu (y)\right) ^{
\frac{1}{p}}.  \nonumber
\end{eqnarray}
The constants of these equivalences depend only on $n,p,$
and $C_{\sigma },C_{\sigma }^{\prime },C_{\sigma }^{\prime
\prime }$ (for the case (i)) and $C_{\mu },C_{\mu }^{\prime
},C_{\mu }^{\prime \prime }$ (for (ii)).
\end{theorem}
\par The results of Theorems \reff{EXT-DIFF-15} and
\reff{EXT-DIFF-S15} can be simplified essentially whenever
$S$ or ${\bf R}^{n}\setminus S$ possesses certain
``plumpness'' properties.
\par Following Triebel \cite{T3}, Definition 9.16, we
introduce the next
\begin{definition} \lbl{BALLCND} {\em We say that a closed
set $A$ satisfies the {\it ball condition} if there exists
a constant $\beta _{A}>0$ such that the following statement
is true: For each Euclidean ball $B$ of diameter at most
$1$ whose center lies in $A,$ there exists a ball
$B^{\prime }\subset B\setminus A$ such that
$\operatorname{diam}B\le \beta
_{A}\operatorname{diam}B^{\prime }.$}
\end{definition}
\par It can be readily seen that $A$ satisfies the ball
condition if and only if $\rho _{A}(x,y)\sim \Vert x-y\Vert
$ for all $x,y$ such that $\Vert x-y\Vert \le 1$.
\par As the next theorem shows, if $\partial S$ satisfies
the ball condition then, in the equivalence of Theorem
\reff{EXT-DIFF-S15},\thinspace (i), one can omit the third
term on the right-hand side.
\begin{theorem}\lbl{EXT-CB-DM} Let $p\in (n,\infty )$,
$\Omega :=int(S)$, and let $\sigma $ be a measure with
$\operatorname{supp}\sigma =\partial S$ satisfying
conditions \rf{DNSG} and \rf{B1SIGMA}. Suppose that
$\partial S$ satisfies the ball condition. Then, for every
$f\in C(S)$, we have
\begin{eqnarray} \Vert f\Vert _{W_{p}^{1}({\bf
R}^{n})|_{S}} &\sim &\Vert f|\,_{\Omega }\Vert
_{W_{p}^{1}(\Omega )}+\Vert f\Vert _{L_{p}(\sigma )}
\label{Felafel} \\
&+&\left( ~\iint\limits_{\Vert x-y\Vert <1}
\frac{|f(x)-f(y)|^{p}}{\Vert
x-y\Vert ^{p-n}\sigma (Q(x,\Vert x-y\Vert ))^{2}}
\,d\sigma (x)\,d\sigma
(y)\right) ^{\frac{1}{p}}
\nonumber
\end{eqnarray}
with constants of equivalence depending only on
$n,p,C_{\sigma },C_{\sigma }^{\prime },C_{\sigma }^{\prime
\prime },$ and the ball condition constant $\beta_{\partial
S}$.
\end{theorem}
\par In particular, if $\Omega :=int(S)=\emptyset $, then,
clearly, $S=\partial S$, so that we can put $\sigma =\mu$.
Also, since $\Omega =\emptyset $, the term $\Vert
f|\,_{\Omega }\Vert _{W_{p}^{1}(\Omega )}$ in \rf{Felafel}
can be omitted. Thus we have
\begin{theorem} \lbl{EXT-THIN} Let $p\in (n,\infty )$ and
let $\mu$ be a measure with $\operatorname{supp}\mu
=\partial S$ satisfying conditions \rf{DN} and \rf{B1}.
Assume that $int(S)=\emptyset $ and that $S$ satisfies the
ball condition. Then, for every $f\in C(S)$,
\bel{EE} \Vert f\Vert _{W_{p}^{1}({\bf R}^{n})|_{S}}\sim
\Vert f\Vert _{L_{p}(\mu )}+\left( \,\,\iint\limits_{\Vert
x-y\Vert <1}\frac{|f(x)-f(y)|^{p}}{\Vert x-y\Vert ^{p-n}\mu
(Q(x,\Vert x-y\Vert ))^{2}} \,d\mu (x)d\mu (y)\right)
^{\frac{1}{p}} \ee
with constants of equivalence depending only on $n,p,C_{\mu
},C_{\mu }^{\prime },C_{\mu }^{\prime \prime },$ and $\beta
_{S}$.
\end{theorem}
\par Norms similar to the one on the right hand side of
\rf{EE} have been used by Dynkin \cite{D1} to prove results
about the trace to the boundary of the unit circle of
analytic functions whose derivatives are in the Hardy space
$H^{p}$. Note also that similar norms have also been used
by Jonsson \cite{J} to characterize the restrictions of
Besov spaces to closed subsets of ${\bf R}^{n}$. See
Section 6 for more details about this.
\par Theorem \reff{EXT-THIN} implies several known results
about traces of Sobolev space to different classes of
``thin'' (or ``porous'') subsets of ${\bf R}^{n}$. For
instance, let $S={\bf R}^{d}\subset {\bf R}^{n}$ where
$0<d<n$. As we have seen in Remark \reff{SB-TRN}, in this
case, for $p>n$, the right hand side of \rf{EE} coincides
with one of the norms for the Besov space
$B_{p,p}^{1-\frac{n-d}{p}}({\bf R}^{d})$. This implies the
classical trace description of Sobolev spaces to ${\bf
R}^{d}$,
\bel{S-RD} W_{p}^{1}({\bf R}^{n})|_{{\bf
R}^{d}}=B_{p,p}^{1-\frac{n-d}{p}}({\bf R}^{d}), \ee
see, e.g. \cite{BIN}. (Note that this isomorphism holds
{\it for all} $p\ge 1,0<q\le \infty $.)
\par Jonsson and Wallin \cite{J1,JW} generalized
isomorphism \rf{S-RD} to the family of $d$-sets in $\RN$.
We recall that a Borel set $S$ in $\RN$ is called a
$d$-set, $0 < d\le n$, if the $d$-dimensional Hausdorff
measure $\Hc_d$ on $S$, $\mu:=\Hc_d|_S$, satisfies the
following condition: there exist constants $C_1,C_2>0,$
such that
\bel{DEF-D} C_1 r^d\le \mu(Q(x,r))\le C_2 r^d~~~~{\rm
for~all~~} x\in S ~~ {\rm and}~~ 0< r\le 1. \ee
Given a $d$-set $S\subset\RN$, the Besov space
$B^s_{p,p}(S), 0<\alpha< 1, 1\le p<\infty,$ is defined as
the family of all functions $f\in L_p(S;\mu)$ such that the
following norm
$$ \|f\|_{B^s_{p,p}(S)}:=\|f\|_{L_p(S;\mu)}+
\left(~\iint\limits_{\|x-y\|<1}
\frac{|f(x)-f(y)|^p}{\|x-y\|^{d+sp}}
\,d\mu(x)d\mu(y)\right)^{\frac1p} $$
is finite, see \cite{JW}, p. 103. Jonsson and Wallin proved
that for every $d$-set $S\subset\RN$ and every
$1<p<\infty$, $0<n-d<p$,
\bel{SOB-DSET} \WP|_S=B^{1-\frac{n-d}{p}}_{p,p}(S). \ee
\par It can be readily seen that every $d$-set with $d<n$
satisfies the ball condition and its interior is empty, so
that the isomorphism \rf{SOB-DSET} for $p>n$ also follows
from Theorem \reff{EXT-THIN}.
\par Thus Theorem \reff{EXT-THIN} extends both of the
isomorphisms \rf{S-RD} and \rf{SOB-DSET} (in the case
$p>n$) to all of the family of sets $S\subset {\bf R}^{n}$
which have empty interior and satisfy the ball condition.
In Section 7 we present an example (see Example
\reff{EXND}) of a constructive description of the trace
space $W_{p}^{1}({\bf R}^{2})|_{S}$, $ p>2$, for such a set
$S\subset {\bf R}^{2}$ which is not a $d$-set for any
$d>0$.
\medskip
\par {\bf 1.3 A Whitney-type extension theorem for
$W_{p}^{k}(R^{n})$-functions, $p>n$.} The proof of the
result which we describe in this subsection does not in
fact appear in this paper. Nevertheless, we mention it
here, in view of its close connection with the topics of
this paper. Its proof will appear in the forthcoming paper
\cite{S5} which will be devoted to a generalization of the
classical Whitney extension theorem.
\par As usual we let $L_{p}^{k}({\bf R}^{n})$ denote the
(homogeneous) Sobolev space of all functions defined on
${\bf R}^{n}$ whose distributional partial derivatives of
order $k$ belong to $L_{p}({\bf R}^{n})$. $L_{p}^{k}({\bf
R} ^{n})$ is normed by
$$ \Vert f\Vert _{L_{p}^{k}({\bf R}^{n})}:=\sum \{\Vert
D^{\alpha }f\Vert _{L_{p}({\bf R}^{n})}:|\alpha |=k\}. $$
By the Sobolev imbedding theorem, see e.g., \cite{M}, p.\
60, every $f\in L_{p}^{k}({\bf R}^{n}),p>n,$ can be
redefined, if necessary, in a set of Lebesgue measure zero
so that it belongs to the space $C^{k-1,\alpha }({\bf
R}^{n})$ of functions whose partial derivatives of order
$k-1$ satisfy the H\"{o}lder condition of order $\alpha
:=1-\frac{n}{p}$:
$$ |D^{\beta }f(x)-D^{\beta }f(y)| \le C(n,p)\Vert f\Vert
_{L_{p}^{k}({\bf R} ^{n})}\Vert x-y\Vert ^{1-\frac{n}{p}},
~~~x,y,\in {\bf R}^{n},|\beta |=k-1. $$
\par This means that, for $p>n$, we can identify each
element $f\in L_{p}^{k}({\bf R}^{n})$ with its unique
$C^{k-1}$-representative, and this enables us to define the
trace space $L_{p}^{k}({\bf R}^{n})|_{S},$ $p>n,$ of all
restrictions of $L_{p}^{k}({\bf R}^{n})$-functions to an
arbitrary closed set $S$,
$$ L_{p}^{k}({\bf R}^{n})|_{S}:=\{f:S\to {\bf R}:{\rm
there~is~}\ F\in L_{p}^{k}({\bf R}^{n})\cap C^{k-1}({\bf
R}^{n})~{\rm such~that~}F|_{S}=f\}. $$
We equip this space with the standard trace norm
$$ \Vert f\Vert _{L_{p}^{k}({\bf R}^{n})|_{S}}
:=\inf \{\Vert F\Vert _{L_{p}^{k}(%
{\bf R}^{n})}:F\in L_{p}^{k}({\bf R}^{n})
\cap C^{k-1}({\bf R}^{n})~~{\rm %
and~~~}F|_{S}=f\}.
$$
\par Let $\PK$ be the space of polynomials of degree at
most $k$ defined on $\RN$. Given a $k$-times differentiable
function $F$ and a point $x\in\RN,$ we let $T_{x}^{k}(F)$
denote the Taylor polynomial of $F$ at $x$ of degree at
most $k$:
$$ T_{x}^{k}(F)(y):=\sum_{|\alpha|\leq
k}\frac{1}{\alpha!}(D^{\alpha}F)(x)(y-x)^{\alpha}~,~~y\in
\RN. $$
\par Recall that, by Whitney's theorem \cite{W1}, given a
family of polynomials $\{P_x\in\Pc_{k-1}:x\in S\}$, there
exists a function $F\in L^k_\infty(\RN)$ such that
$T_{x}^{k-1}(F)=P_{x}$ for all $x\in S,$ if and only if
there exists a constant $\lambda>0$ such that, for every
$\alpha, |\alpha|\le k-1,$ we have
$$ |D^\alpha(P_x-P_y)(x)| \le \lambda\,
\|x-y\|^{k-|\alpha|},~~~~x,y\in S. $$
\par We have obtained the following generalization of this
theorem to the case $p>n$.
\begin{theorem}\lbl{EXTK}(\cite{S5}) Given positive
integers $k$ and $n$, a number $p\in (n,\infty)$ and a
family of polynomials $\{P_x\in\Pc_{k-1}:x\in S\}$, there
exists a $(k-1)$-continuously differentiable function
$F\in\LKP$ such that $T_{x}^{k-1}(F)=P_{x}$ for every $x\in
S$ if and only if there exists a constant $\lambda>0$ such
that for every finite family $\{Q_i:i=1,...,m\}$ of
disjoint cubes in $\RN$, every $x_i,y_i\in (11Q_i)\cap S,$
and every multi-index $\alpha, |\alpha|\le k-1,$ the
following inequality
$$
\sum_{i=1}^m\frac{|D^\alpha(P_{x_i}-P_{y_i})(x_i)|^p}
{(\diam Q_i)^{(k-|\alpha|)p-n}} \le\lambda
$$
holds. Moreover, $\|f\|_{\LKP|_S}\sim\inf
\lambda^{\frac{1}{p}}$ with constants of equivalence
depending only on $k,n$ and $p$.
\end{theorem}
\par Analogously to the role played by the Whitney
extension theorem in the study of the case $p=\infty$, we
expect that Theorem \reff{EXTK} will turn out to be useful
for the study of the following variant of the classical
{\it Whitney extension problem} \cite{W1,W2} for the space
$\LKP$, $p>n$:
\par {\it Given $p\in(n,\infty]$, a positive integer $k$,
and an arbitrary function $f:S\to\R$, what is a necessary
and sufficient condition for $f$ to be the restriction to
$S$ of a $(k-1)$-times continuously differentiable function
$F\in \LKP$?}
\par We observe that the Whitney problem for the spaces
$L^k_\infty$, $C^{k,\alpha}$, $C^k$, and their
generalizations has attracted a lot of attention in recent
years. We refer the reader to \cite{BS1}-\cite{BS3},
\cite{F1}-\cite{F4}, \cite{BM,BMP,S6,Zob1} and references
therein for numerous results in this direction, and for a
variety of techniques for obtaining them.
\par We also note that, as in the case $k=1$ (see Theorem
\reff{EXT-SIMP}), the classical Whitney extension method
provides an almost optimal extension operator for all the
scale of $\LKP$-spaces, $p>n$.
\medskip
\par {\bf Acknowledgement.} I am very thankful to Yu.
Brudnyi, M. Cwikel and V. Maz'ya for helpful discussions
and valuable advice.
\SECT{2. Traces of Sobolev and Besov spaces and local
oscillations.}{2}
\indent
\par {\bf 2.1 Local oscillations and $A_p$-functionals.}
\par The proofs of the results formulated in subsections
1.1 and 1.2 are based on estimates of local oscillations
and corresponding oscillation functionals
($A_p$-functionals) of functions and their extensions. The
$A_p$-functionals, see Definitions \reff{DEF-AP} and
\reff{PR-AP}, play a role of certain generalized moduli of
smoothness on a closed set providing one more approach to
description of traces of Sobolev spaces. We present these
results in Theorems \reff{EXT1} and \reff{EXT2} below.
\par Let us fix additional notation. Given a bounded
function $f$ defined on a set $U\subset\RN$ we let
$\Ec(f;U)$ denote its {\it oscillation} on $U$, i.e., the
quantity
\bel{OSC} \Ec(f;U):=\sup_{x,y\in U} |f(x)-f(y|. \ee
We put $\Ec(f;U):=0$ whenever $U=\emp$.
\par We call a finite family of disjoint cubes a {\it
packing}.
\par The point of departure for our approach is a
description of Sobolev spaces in terms of local
oscillations due to Brudnyi \cite{Br1}: a continuous
function $f\in\LOP, p>n,$ if and only if there is a
constant $\lambda>0$ such that for every $t>0$ we have
\bel{SPACK}
\sup_{\pi}\{\sum_{Q\in\pi}|Q|\Ec(f;Q)^p\}^{\frac{1}{p}}\le
\lambda t. \ee
Here the supremum is taken over all packings $\pi=\{Q\}$ of
equal cubes with $\diam Q\le t$. Moreover,
$\|f\|_{\LOP}\sim\inf\lambda$ with constants of equivalence
depending only on $n$ and $p$.
\par Inequality \rf{SPACK} motivates us to introduce the
following definition: given a locally bounded function
$f:\RN\to\R$ we let $\AP(\cdot;f)$ denote the function
\bel{DI} \AP(f;t):=
\sup_{\pi}\{\sum_{Q\in\pi}|Q|\Ec(f;Q)^p\}^{\frac{1}{p}} \ee
where $\pi$ runs over all packings  of equal cubes with
diameter at most $t$. Now, by \rf{SPACK}, for every $p>n$
we have
\bel{EM} \|f\|_{\LOP}\sim  \sup_{t>0}{\AP(f;t)}/{t}. \ee
\par Observe that, an approximation functional, similar to
$A_p$\,, has been introduced by Yu. Brudnyi in \cite{Br1}.
Let us recall its definition.
\par Given a cube $Q$ and a function $f\in L_q(Q), q>0,$ we
let $\Ec_1(f;Q)_{L_q}$ denote the so-called {\it normalized
local best approximation} of $f$ on $Q$ in $L_q$-norm:
\bel{ESP} \Ec_1(f;Q)_{L_q}:=|Q|^{-\frac{1}{q}}\inf_{c\in\R}
\|f-c\|_{L_q(Q)}=\inf_{c\in\R}\left(\frac{1}{|Q|}
\intl_Q|f-c|^q\,dx\right)^{\frac{1}{q}}. \ee
Following \cite{Br1}, given  a function $f\in
L_{q,\,loc}(\RN),~q,p\in(0,\infty],$ we define its {\it
$(1,p)$-modulus of continuity in $L_q$} by letting
\bel{KPU} \omp(f;t)_{L_q}:=\sup_{\pi}\left\{\sum_{Q\in \pi}
|Q|\,\Ec_1(f;Q)^p_{L_q}\right\}^{\frac{1}{p}}, ~~~t>0. \ee
Here $\pi$ runs over all packings of equal cubes in $\RN$
with diameter at most $t$.
\par There is an integral form of the $(1,p)$-modulus of
continuity,
\bel{IN1} \omp(f;t)_{L_q}\sim
\|\Ec_1(f;Q(\cdot;t))_{L_q}\|_{L_p(\RN)}, \ee
where constants of equivalence depend only on $n,p,$ and
$q$, see \cite{Br2}, Chapter 3, and also \cite{S4}. It is
also clear that, $\Ec_1(f;Q)_{L_\infty}\sim \Ec(f;Q)$ so
that
\bel{OMAP} \omp(f;\cdot)_{L_\infty}\sim\AP(f;\cdot) \ee
(with absolute constants). Therefore, by \rf{IN1}, for
every locally bounded function $f$ we have
\bel{AP-INT} \AP(f;t)\sim
\|\Ec(f;Q(\cdot;t))\|_{L_p(\RN)},~~~~t>0, \ee
with constants of equivalence depending only on $n$ and
$p$.
\par Let us also note an important property of the
$(1,p)$-modulus of continuity, see, e.g. \cite{Br1} or
\cite{DL, DS1}: for every $f\in L_{p,\,loc}(\RN)$ and $1\le
p\le\infty$
\bel{EQM}
\omega(f;t)_{L_p}\sim \omp(f;t)_{L_p},~~~~t>0,
\ee
with absolute constants of equivalence. (Recall that
$\omega(f;\cdot)_{L_p}$ , see \rf{DOM}, denotes the modulus
of smoothness of $f$ in $L_p$.)
\par Let us extend definition \rf{DI} of the
$A_p$-functional to an arbitrary closed set $S\subset\RN$.
\begin{definition}\lbl{DEF-AP}{\em Let $f:S\to\R$ be a
locally bounded function. We define the functional
$\AP(f;\cdot:S)$ by the formula
\bel{OM} \AP(f;t:S):=\sup_{\pi}\,\, \{\sum_{Q\in\pi}
|Q|\Ec(f;Q\cap S)^p\}^{\frac{1}{p}}, ~~~~t>0, \ee
where the supremum is taken over all packings $\pi=\{Q\}$
of {\it equal cubes  centered in $S$ with diameter at most
$t$.}
}\end{definition}
\par Thus, $\AP(F;t)=\AP(F;t:\RN), t>0$.
\par Observe that equivalence \rf{EM} provides a necessary
condition for a function $f:S\to\R$ to belong to the trace
space $\LOP|_S, p>n$. In fact, since
$$
\AP(f;t:S)\le \AP(F;t:\RN)=\AP(F;t)
$$
for every extension $F\in\LOP$ of $f$, by \rf{EM},
\bel{N0} \sup_{t>0} \AP(f;t:S)/t\le C\|f\|_{\LOP|_S}, \ee
so that
$$
f\in\LOP|_S~~\Longrightarrow~~\sup_{t>0}
\frac{\AP(f;t:S)}{t}<\infty.
$$
However, simple examples (for instance, a straight line in
$\RN$,) show that, in general, this condition is not
sufficient. Below we present an additional condition which
will enable us to describe the restrictions of Sobolev
functions via $A_p$-functionals.
\par To formulate this new condition we are needed the
following definitions and notation.
\begin{definition}\lbl{PR-Q}{\em Let $0\le\alpha\le 1$. A
cube $Q\subset\RN$ is said to be {\it $\alpha$-porous} with
respect to $S$ if there exists a cube $Q'$ such that
$$
Q'\subset Q\setminus S~~{\rm and ~~} \diam Q'\ge \alpha
\diam Q.
$$
We say that a packing $\pi$ is $\alpha$-porous with respect
to $S$ if every $Q\in\pi$ is $\alpha$-porous (with respect
to $S$).
}\end{definition}
\par Let us introduce a corresponding $A_p$-functional over
the family of all $\alpha$-porous packings.
\begin{definition}\lbl{PR-AP}{\em Given $0\le\alpha\le 1$
and a locally bounded function $f:S\to R$ we define the
functional $\PRM(f;t:S)$ by letting
\bel{OMP} \PRM(f;t:S):=\sup_{\pi}\,\,
\{\sum_{Q\in\pi}|Q|\Ec(f;Q\cap \partial S)^p\}
^{\frac{1}{p}}. \ee
Here $\pi=\{Q\}$ runs over all {\it $\alpha$-porous (with
respect to $S$) packings} of equal cubes with centers in
$\partial S$ and diameter at most $t$.
}\end{definition}
\par Note that
\bel{APCM} \PRM(f;t:S) \le\AP(f;t:S),~~~~ 0<\alpha<1,~ t>0,
\ee
and
$$
\Ac_{p,\alpha_2}(f;t:S)\le \Ac_{p,\alpha_1}(f;t:S),
~~~~0<\alpha_1<\alpha_2<1,~t>0.
$$
\begin{theorem}\lbl{EXT1} Let $\alpha\in(0,3/20]$ and let $
p\in(n,\infty)$. Then for every function $f\in C(S)$ the
following equivalence
$$
\|f\|_{\LOP|_S}\sim \sup_{0<t\le 2\,\diam S}
\frac{\AP(f;t:S)}{t}+ \left\{\intl_0^{\diam
S}\PRM(f;t:S)^p\,\frac{dt}{t^{p+1}}\right\}^{\frac{1}{p}}
$$
holds with constants depending only on $n,p$ and $\alpha$.
\end{theorem}
\par We turn to characterization of the traces of
$W^1_p$-functions via $A_p$- and
$\Ac_{p,\alpha}$-functionals. Let $\ve>0$ and let
$T:S_\ve\to S$ be a measurable mapping satisfying
inequality \rf{TR}. In addition, we assume that
\bel{ADT}
T(x)\in\partial S,~~~~x\in S_\ve\setminus S.
\ee
(In particular, formula \rf{TWQ} provides an example of
such a mapping).
\begin{theorem}\lbl{EXT2} Let
$\alpha\in(0,3/(10+10\theta)]$ and let $p\in(n,\infty)$.
Then for every $f\in C(S)$ we have
$$
\|f\|_{\WP|_S}\sim \|f\circ T\|_{L_p(S_\ve)}+
\sup_{0<t\le \ve} \frac{\AP(f;t:S)}{t}+
\left\{\intl_0^\ve\PRM(f;t:S)^p\,
\frac{dt}{t^{p+1}}\right\}^{\frac{1}{p}}
$$
with constants of equivalence depending only on
$n,p,\alpha,\theta$ and $\ve$.
\end{theorem}
\par Observe that, by choosing a more delicate Whitney
covering of $\RN\setminus S$, see Remark  \reff{DEL-W}, the
intervals $(0,3/20]$ and $(0,3/(10+10\theta)]$ in Theorems
\reff{EXT1} and \reff{EXT2} can be replaced with slightly
smaller intervals $(0,1/4)$ and $(0,1/(2+2\theta))$
respectively.
\par {\bf 2.1 Besov spaces on closed subsets of $\RN$.}
\par A. Jonsson \cite{J} characterized the trace of the
Besov space $\BS,\, s\in(n/p\,,1)$, to an arbitrary closed
subset $S\subset\RN$ via doubling measures supported on
$S$, see Theorem \reff{EXT-BMU}. In Section 6 we present
another intrinsic description of the space $\BS|_S$ which
involves only the values of a function defined on $S$ and
does not use doubling measures.
\par First we observe that every function $f\in \BS,$ $s\in
(n/p,1),$ can be redefined in a set of measure zero to
satisfy the H\"{o}lder condition of order $s-n/p$\,, see,
e.g., \cite{BIN} or \cite{N}, Chapter 6. Therefore, similar
to the case of Sobolev spaces, in what follows we will
identify $f$ with its continuous refinement. As usual,
given a closed subset $S\subset\RN$, we define the trace
space
$$ \BS|_S:=\{f:S\to\R: {\rm there~is~a~continuous}\ \  F\in
\BS\ \ {\rm such~that}\ \ F|_S=f\}, $$
and equip this space with the standard quotient space norm.
\par The next theorem provides an intrinsic
characterization of the restrictions of the Besov space
$\BS,$ $n/p<s<1,$ to an arbitrary closed subset
$S\subset\RN$. In its formulation given $\ve>0$ and
$\theta\ge 1$, we again let $T$ denote a measurable mapping
from $S_\ve$ into $S$ satisfying on $S_\ve$ inequality
\rf{TR}.
\begin{theorem}\lbl{EXT-BS-INTRO} A function $f\in C(S)$
can be extended to a (continuous) function $F\in\BS,
n/p<s<1,$ if and only if $f\circ T\in L_p(S_\ve)$ and
$$
\intl_0^{\ve}\AP(f;t:S)^q\,
\frac{dt}{t^{sq+1}}<\infty.
$$
Moreover,
\bel{BSTRN} \|f\|_{\BS|_S}\sim \|f\circ T\|_{L_p(S_\ve)}+
\left\{\intl_0^\ve\AP(f;t:S)^q\,
\frac{dt}{t^{sq+1}}\right\}^{\frac{1}{q}} \ee
with constants of equivalence depending only on
$n,s,p,q,\theta$ and $\ve$.
\end{theorem}
\par Similar to the case of the Sobolev space, Whitney's
method works for all the scale of the Besov spaces $\BS,$
$n/p<s<1,$ providing an almost optimal extension operator.
For Whitney-type extension theorems for the Besov spaces we
refer the reader to \cite{JW,Gud}.
\SECTLONG{3. The Whitney covering and}{local approximation
properties of  $\WP$-functions.}{3}
\indent
\par Throughout the paper $C,C_1,C_2,...$ will be generic
positive constants which depend only on parameters
determining sets (say, $n$, doubling or ``porosity"
constants, etc.) or function spaces ($p,q,s,$ etc). These
constants can change even in a single string of estimates.
The dependence of a constant on certain parameters is
expressed, for example, by the notation $C=C(n,p)$. We
write $A\sim B$ if there is a constant $C\ge 1$ such that
$A/C\le B\le CA$.
\par Recall that $|A|$ denotes the Lebesgue measure of $A$.
Given a Borel measure $\mu$ on $\RN$, a $\mu$-measurable
set $A\subset\RN$ and $p\in(0,\infty]$ we let $L_p(A:\mu)$
denote the $L_p$-space (with respect to the measure $\mu$)
of functions on $A$. We write $L_p(\mu)$ whenever $A=\RN$,
and $L_p(A)$ whenever $\mu$ is the Lebesgue measure.
\par If $f\in L_1(A:\mu)$ and $\mu(A)<\infty$, we let
$f_{A,\,\mu}:=\mu(A)^{-1}\int_Af\,d\mu$ denote the
$\mu$-average of $f$ on $A$. By $f_A$ we denote the average
of $f$ on $A$ with respect to the Lebesgue measure.
\par Recall that $Q=Q(x,r)$ denotes the closed cube in
$\RN$ centered at $x\in \RN$ with side length $2r$ whose
sides are parallel to the coordinate axes. We write $x_Q=x$
and $r_Q=r$ so that $Q=Q(x_Q,r_Q)$. Also recall that
$\lambda>0$ and $Q=Q(x,r)$ we let $\lambda Q$ denote the
cube $Q(x,\lambda r)$. By $Q^*$ we denote the cube
$Q^*:=\frac{9}{8}Q$.
\par It will be convenient for us to measure distances in
$\RN$ in the uniform norm
\bel{UNINORM}
\|x\|:=\max\{|x_i|:~i=1,...,n\}, \ \ \
x=(x_1,...,x_n)\in\RN.
\ee
Thus every cube
$$
Q=Q(x,r):=\{y\in\RN:\|y-x\|\le r\}
$$
is a ``ball" in $\|\cdot\|$-norm  of ``radius" $r$ centered
at $x$.
\par Given subsets $A,B\subset \RN$,  we put
$ \diam A:=\sup\{\|a-a'\|:~a,a'\in A\} $
and
$$
\dist(A,B):=\inf\{\|a-b\|:~a\in A, b\in B\}.
$$
For $x\in \RN$ we also set $\dist(x,A):=\dist(\{x\},A)$. By
$int(A)$ we denote the interior of $A$ and by $\chi_A$ the
characteristic function of $A$; we put $\chi_A\equiv 0$
whenever $A=\emptyset$.
\par As usual, $\nabla F(x):=\left(\frac{\partial
F}{\partial x_1}(x),...,\frac{\partial F}{\partial
x_n}(x)\right)$ stands for the gradient of a function $F$
at $x$; thus, by \rf{UNINORM}, $\|\nabla
F(x)\|=\max\left\{\left|\frac{\partial F}{\partial
x_i}(x)\right|:~i=1,...,n\right\}.$
\par Let $\pi=\{Q\}$ be a family of cubes in $\RN$. By
$M_\pi$ we denote its {\it covering multiplicity}, i.e.,
the minimal positive integer $M$ such that every point
$x\in\RN$ is covered by at most $M$ cubes from $\pi$. Thus
$$
M_{\pi}:=\sup_{x\in\RN}\sum_{Q\in \pi}\chi_Q(x).
$$
\par Everywhere in this paper $S$ denotes a {\it closed}
subset of $\RN$. Since the set $\RN\setminus S$ is open, it
admits a Whitney covering $\TW$ with the following
properties.
\begin{theorem}\lbl{Wcov}(see, e.g. \cite{St}, or \cite{G})
$\TW=\{Q_k\}$ is a countable family of cubes  such that
\par (i). $\RN\setminus S=\cup\{Q:Q\in \TW\}$;
\par (ii). For every cube $Q\in \TW$
\bel{DKQ}
\diam Q\le \dist(Q,S)\le 4\diam Q;
\ee
\par (iii). The covering multiplicity $M_{\TW}$ of $\TW$ is
bounded by a constant $N=N(n)$. Thus every point of
$\RN\setminus S$ is covered by at most $N$ cubes from
$\TW$.
\end{theorem}
\par We are also needed certain additional properties of
Whitney's cubes which we present in the next lemma. These
properties easily follow from Theorem \reff{Wcov}.
\begin{lemma}\lbl{Wadd}
(1). If $Q,K\in \TW$ and $Q^*\cap K^*\ne\emptyset$, then
$$
\frac{1}{4}\diam Q\le \diam K\le 4\diam Q.
$$
(Recall that $Q^*:=\frac{9}{8}Q$.)
\par (2). For every cube $K\in \TW$ there are at most
$N=N(n)$ cubes from the family
$ {\mathbb W}^*(S):=\{Q^*:Q\in \TW\} $
which intersect $K^*$.
\par (3). If $Q,K\in \TW$, then $Q^*\cap K^*\ne\emptyset$
if and only if  $Q\cap K\ne\emptyset$.
\end{lemma}
\par Given $\ve>0$ we put
\bel{DQEP} \TWB{\ve}:=\{Q\in\TW:~\diam Q\le\ve\}. \ee
\par The proofs of the necessity part of several results
stated in Sections 1 and 2 are based on a series of
auxiliary statements which we present in this section.
First of them is the following combinatorial
\begin{theorem}(\cite{BrK,Dol})\lbl{TFM} Let $\Ac=\{Q\}$ be
a collection of cubes in $\RN$ with the covering
multiplicity $ M_{\Ac}<\infty$. Then $\Ac$ can be
partitioned into at most $N=2^{n-1}(M_{\Ac}-1)+1$ families
of disjoint cubes.
\end{theorem}
\begin{lemma}\lbl{FM} Let $0<\alpha\le 1\le\beta\le
\lambda$ and let $t>0$. Let $\pi=\{Q\}$ be a family of
cubes in $\RN$ with the finite covering multiplicity such
that $\alpha t\le\diam Q\le\beta t$ for every $Q\in\pi$.
\par Let us assign every $Q\in\pi$ a cube $K_Q$ of diameter
at most $\lambda t$ and such that $K_Q\supset Q$. Then the
family of cubes $ \Ac=\{K_Q:~Q\in\pi\}$ can be partitioned
into at most $N=N(n,\alpha,\beta,\lambda,M_\pi)$ families
of disjoint cubes.
\end{lemma}
\par {\it Proof.} One can easily show that the covering
multiplicity of the family $\Ac$ is bounded by a constant
depending only on $n,\alpha,\beta,\lambda,$ and $M_\pi$. It
remains to make use of Theorem \reff{TFM} and the lemma
follows.\bx
\par The following proposition presents a classical Sobolev
imbedding inequality for the case $p>n$, see, e.g.
\cite{M}, p. 61, or \cite{MP}, p. 55. This inequality is
also known in the literature as Sobolev-Poin\'care
inequality (for $p>n$).
\begin{proposition}\lbl{EDIFF} Let $F\in L^1_p(\Omega)$ be
a continuous function defined on a domain
$\Omega\subset\RN$ and let $n<q\le p<\infty$. Then for
every cube $Q\subset\Omega$ and every $x,y\in Q$ the
following inequality
$$
|F(x)-F(y)|\le C(n,q)\,
r_Q\left(\frac{1}{|Q|}\intl_Q \|\nabla F(z)\|^q\,dz\right)
^{\frac{1}{q}}
$$
holds.
\end{proposition}
\par This inequality enables us to prove the following
useful property of the $A_p$-functional, see Definition
\reff{DEF-AP}.
\begin{proposition}\lbl{AP-DS} Let $p>n$ and let
$\Omega:=int(S)$. Assume that $f\in C(S)$ and
$f|_\Omega\in\LOPO$. Then
$$
\AP(f;t:S)\le
C(n,p)\{\AP(f;4t:\partial S)+t\|f|_\Omega\|_{\LOPO}\},
~~~~t>0.
$$
\end{proposition}
\par {\it Proof.} Let $t>0$ and let $\pi=\{Q\}$ be a
packing of equal cubes of diameter $\diam Q=2\tau\le t$
centered in $S$. If $Q\subset\Omega$ for every $Q\in\pi$,
then, by Proposition \reff{EDIFF},
\be J(f;Q)&:=&\sum_{Q\in\pi} |Q|\Ec(f;Q\cap
S)^p=\sum_{Q\in\pi} |Q|\Ec(f;Q)^p\nn\\
&\le& C\sum_{Q\in\pi} r_Q^p\intl_{Q} \|\nabla
f(y)\|^p\,dy\le Ct^p\intl_{\Omega} \|\nabla
f(y)\|^p\,dy\nn\ee
so that
\bel{EIN} J(f;Q)\le Ct^p\|f|_\Omega\|_{\LOPO}^p. \ee
\par Now assume that $Q\cap \partial S\ne\emp$ for each
$Q\in\pi$; thus there exists a point $b_Q\in Q\cap \partial
S$, $Q\in\pi$. Put $Q':=Q(b_Q,2\tau)=Q(b_Q,2r_Q)$. Then,
clearly, $Q\subset Q'\subset 3Q$.
\par We have
$$ \Ec(f;Q\cap S)=\sup_{x,y\in Q\cap S}|f(x)-f(y)|\le
2\sup_{x\in Q\cap S}|f(x)-f(b_Q)| $$
so that
$$ \Ec(f;Q\cap S)\le 2\max\{\sup_{x\in Q\cap \partial
S}|f(x)-f(b_Q)|, \sup_{x\in Q\cap \Omega}|f(x)-f(b_Q)|\}.
$$
Hence,
\bel{DFG1} \Ec(f;Q\cap S)\le  2\{\Ec(f;Q'\cap \partial S)+
\sup_{x\in Q\cap \Omega}|f(x)-f(b_Q)|\}. \ee
Given $x\in Q\cap \Omega$ we let $a_x$ denote a point
nearest to $x$ on $\partial S$. Thus
$$ \|x-a_x\|=\dist(x,\partial S)\le\|x-b_Q\|\le r_Q=\tau.
$$
Put $Q^{(x)}:=Q(x,\dist(x,\partial S))$. Then, clearly,
$int(Q^{(x)})\subset\Omega$ and $a_x\in Q^{(x)}\cap\partial
S$ so that for every $y\in Q^{(x)}$ we have
$$ \|y-x_Q\|\le \|y-x\|+\|x-x_Q\|\le \dist(x,\partial
S)+r_Q\le 2r_Q $$
proving that $Q^{(x)}\subset  2Q\subset 2Q'.$
\par By Lemma \reff{FM}, we can assume that the family of
cubes $2\pi':=\{2Q':~Q\in\pi\}$ is a packing. We have
$$ |f(x)-f(b_Q)|\le |f(x)-f(a_x)|+|f(a_x)-f(b_Q)| $$
so that
\bel{GKL} |f(x)-f(b_Q)|\le
\Ec(f,Q^{(x)})+\Ec(f,(2Q')\cap\partial S). \ee
Since $f$ is continuous on $S$ and
$int(Q^{(x)})\subset\Omega$, there exists a cube
$\tQ^{(x)}\subset int(Q^{(x)})$ such that
$\Ec(f,Q^{(x)})\le 2\Ec(f,\tQ^{(x)})$. Also, since
$\tQ^{(x)}\subset Q^{(x)} \subset 2Q$, we have
$r_{\tQ^{(x)}}\le 2r_Q=2\tau$.
\par Then, by Proposition \reff{EDIFF}, for any choice of a
point $x=x(Q)\in Q\cap\Omega$, we have
\be \sum_{Q\in\pi} |Q|\Ec(f;Q^{(x)})^p &\le& C
\sum_{Q\in\pi}
|Q|r^{p-n}_{\tQ^{(x)}}\intl_{\tQ^{(x)}}\|\nabla
f(y)\|^p\,dy\nn\\
&\le& C \tau^n(2\tau)^{p-n} \sum_{Q\in\pi}
\intl_{\tQ^{(x)}}\|\nabla f(y)\|^p\,dy\le C\tau^p
\intl_{\Omega}\|\nabla f(y)\|^p\,dy. \nn \ee
Hence, by \rf{GKL} and Definition \reff{DEF-AP},
\be \sum_{Q\in\pi} |Q|\sup_{x\in Q\cap\Omega}
|f(x)-f(b_Q)|^p&\le& C\{\sum_{Q\in\pi}|Q|\sup_{x\in
Q\cap\Omega}\Ec(f;Q^{(x)})^p\nn\\
&+&\sum_{Q\in\pi}|2Q'|\Ec(f;(2Q')\cap\partial S)^p\}\nn\\
&\le&C\{t^p \intl_{\Omega}\|\nabla f(y)\|^p\,dy
+\AP(f;4t:\partial S)^p\}.\nn \ee
\par Finally, by \rf{DFG1} and Definition \reff{DEF-AP},
\be J(f;\pi)&\le& C\{\sum_{Q\in\pi}
|Q'|\Ec(f;Q'\cap\partial S)^p+ \sum_{Q\in\pi}|Q|\sup_{x\in
Q\cap\Omega} |f(x)-f(b_Q)|^p\}\nn\\
&\le& C\{\AP(f;2t:\partial S)^p+t^p \|f|_\Omega\|_{\LOPO}^p
+\AP(f;4t:\partial S)^p\}. \nn \ee
Combining this inequality with \rf{EIN}, we obtain
$$
J(f;\pi)\le C\{\AP(f;4t:\partial S)^p+t^p
\|f|_\Omega\|_{\LOPO}^p \}
$$
provided $\pi$ is an {\it arbitrary} packing of equal cubes
of diameter at most $t$ centered in $S$. It remains to take
the supremum over all such packings $\pi$, and the
proposition follows.\bx
\begin{lemma}\lbl{LM-FMP} Let $\beta>1,\gamma\ge 1,$ and
let $\pi$ be a family of cubes with the covering
multiplicity $M_\pi<\infty.$ Then for every function $G\in
L_\beta(\RN)$ we have
$$
\{\sum_{Q\in\pi} |Q|\,|G_{\gamma
Q}|^\beta\}^{\frac{1}{\beta}}\le C\|G\|_{L_\beta(\RN)}.
$$
Here $C$ is a constant depending only on $n,\beta,$ and
$M_\pi$.
\end{lemma}
\par {\it Proof.} Let $Q\in\pi$ and let $x\in Q$. Then
$\gamma Q\subset Q_x$ where $Q_x:=Q(x,(\gamma+1)r_Q)$. It
is also clear that $Q_x\subset 3\gamma Q$. Hence
$$
|G_{\gamma Q}|\le\frac{1}{|\gamma Q|}
\intl_{\gamma Q}|G|\,dx\le
C(n)\frac{1}{|Q_x|}\intl_{Q_x}|G|\,dx\le \Mc G(x), ~~~x\in Q.
$$
Here $\Mc$ denotes the Hardy-Littlewood maximal function.
Integrating this inequality over $Q$ we obtain
$$
|Q|\,|G_{\gamma Q}|^\beta\le C\intl_{Q}\Mc G(x)^\beta dx
$$
with $C=C(n,\beta)$. Hence
$$
I:=\sum_{Q\in\pi} |Q|\,|G_{\gamma Q}|^\beta\le
C\sum_{Q\in\pi}\intl_{Q}\Mc G(x)^\beta dx=C
\intl_{\RN}\left(\sum_{Q\in\pi}\chi_Q(x)
\right)\Mc G(x)^\beta dx.
$$
This inequality and the Hardy-Littlewood maximal theorem
yield
$$
I\le C \,M_\pi\intl_{\RN}\Mc G(x)^\beta dx\le C\,M_\pi
\intl_{\RN}|G(x)|^\beta dx,
$$
proving the required inequality $I\le C \|G\|_{L_\beta
(\RN)}^\beta.$\bx
\par This lemma and Theorem \reff{Wcov} imply the following
\begin{lemma}\lbl{LMF} Let $\beta>1$,$\gamma\ge 1$, and let
$A$ be a closed subset of $\RN$. Then for every $G\in
L_\beta(\RN)$ the following inequality
$$
\{\sum_{Q\in\TWU} |Q|\,|G_{\gamma
Q}|^\beta\}^{\frac{1}{\beta}}\le C\|G\|_{L_\beta(\RN)}.
$$
hods. Here $C$ is a constant depending only on $n$ and
$\beta$.
\end{lemma}
\begin{lemma}\lbl{NT1} Let $n\in(p,\infty),$ $\gamma\ge 1,$
and let $\pi=\{Q\}$ be a family of cubes with the covering
multiplicity $M_\pi<\infty.$ Then for every continuous
function $F\in\LOP$ the following inequality
$$
\sum_{Q\in\pi}|Q|\left(\frac{\Ec(F;\gamma Q)}{\diam
Q}\right)^p\,\le C\intl_{\RN}\|\nabla F\|^p\,dx
$$
holds. Here  $C$ is a constant depending only on
$n,p,\gamma,$ and $M_\pi$.
\end{lemma}
\par {\it Proof.} We put $q:=(p+n)/2$ so that $n<q<p.$ By
Proposition \reff{EDIFF},
$$
\Ec(F;\gamma Q)=\sup_{x,y\in\gamma Q}|F(x)-F(y)|\le C(\gamma
r_Q)\left(\frac{1}{|\gamma Q|}\intl_{\gamma Q}
\|\nabla F\|^q\,dx\right)^{\frac1q}
$$
with $C=C(n,p)$. Put $G:=\|\nabla F\|^q$ and $\beta:=p/q.$
Then
$$
I:=\sum_{Q\in\pi}\,|Q|\left(\frac{\Ec(F;\gamma Q)}{\diam
Q}\right)^p\le C\sum_{Q\in\pi}\,|Q|\left(\frac{1}{|\gamma
Q|}\intl_{\gamma Q}\|\nabla F\|^q\,dx\right)
^{\frac{p}{q}}=
C\sum_{Q\in\pi}|Q|\,|G_{\gamma Q}|^\beta.
$$
Since $\beta>1$ and $M_\pi<\infty$, by Lemma \reff{LM-FMP},
$$
I\le C\|G\|^\beta_{L_\beta(\RN)}= C\intl_{\RN}
(\|\nabla F\|^q)^{\frac{p}{q}}\,dx=C\intl_{\RN}
\|\nabla F\|^p\,dx,
$$
proving the lemma.\bx
\par Now, let us replace in Definition \reff{PR-AP} the
boundary $\partial S$ of $S$ with the set $S$ itself. We
denote a functional obtained by $\SPR(f;\cdot:S)$.
\begin{definition}\lbl{PR-APS}{\em Let $f:S\to \R$ be a
locally bounded function and let $\alpha\in (0,1].$ We
define the functional $\SPR(f;\cdot:S)$ by the formula
\bel{OMPS} \SPR(f;t:S):= \sup_{\pi}\,\, \{\sum_{Q\in\pi}
|\,Q|\,\Ec(f;Q\cap S)^p\}^{\frac{1}{p}},~~~t>0. \ee
Here the supremum is taken over all $\alpha$-porous (with
respect to $S$) packings $\pi=\{Q\}$ of equal cubes with
centers in $S$ and diameter at most $t$.
}\end{definition}
\par Clearly, for all $t>0$ and $\alpha\in(0,1]$ we have
\bel{ADS} \PRM(f;t:S)\le \SPR(f;t:S). \ee
\begin{proposition}\lbl{NL} Let $\alpha\in(0,1]$ and let
$p\in(n,\infty)$. Then for every continuous function
$F\in\LOP$ we have
$$
\intl_0^{\infty}\SPR(F|_S;t:S)^p
\, \frac{dt}{t^{p+1}}\le
C(n,p,\alpha)\|F\|_{\LOP}^p.
$$
\end{proposition}
\par {\it Proof.} Put $f:=F|_S$. Since $\SPR(f;\cdot;S)$ is
a non-decreasing function, we have
\be I:=\intl_0^{\infty} \SPR(f;t:S)^p\, \frac{dt}{t^{p+1}}
&=& \sum^{\infty}_{m=-\infty}\,
\intl_{2^{m-1}}^{2^m}\SPR(f;t:S)^p\, \frac{dt}{t^{p+1}}\nn\\
&\le& C\sum^{\infty}_{m=-\infty} 2^{-mp}\SPR(f;2^m:S)^p.
\nn\ee
\par By \rf{OMPS}, for every integer $m$ such that
$\SPR(f;2^m:S)>0$ there exists an $\alpha$-porous packing
$\pi_m=\{Q\}$ which consists of equal cubes centered in $S$
with $r_Q=t_m\le 2^m$ and satisfies the following
inequality:
\bel{AUM} \SPR(f;2^m:S)^p\le 2 \sum_{Q\in\pi_m}
|Q|\Ec(f;Q\cap S)^p. \ee
\par Now fix an integer $m$. Since $\pi_m$ is
$\alpha$-porous, for each $Q\in\pi$ there exists a cube
$Q'\subset Q\setminus S$ such that $r_{Q'}\ge\alpha
r_Q=\alpha t_m$. Since $x_{Q'}\in Q\setminus S$, there is a
Whitney's cube $K_Q\in\TW$ which contains $x_{Q'}$.
\par It can be readily seen that inequality \rf{DKQ}
implies the following two properties of the cube $K_Q$: $
\diam K_Q\le 2^m$ and $\gamma K_Q\supset Q$ for some
$\gamma=\gamma(\alpha)>0$.
\par Hence $\Ec(f;Q\cap S)\le \Ec(F;\gamma K_Q)$ so that by
\rf{AUM}
$$
I\le C\sum^{\infty}_{m=-\infty}
\sum_{Q\in\pi_m} 2^{-mp}|Q|\Ec(F;\gamma K_Q)^p.
$$
Since $\diam K_Q\le 2^m$, we obtain
\be I&\le& C\sum^{\infty}_{m=-\infty} \sum_{Q\in\pi_m}
2^{-m(p-n)}\Ec(F;\gamma K_Q)^p\nn\\
&\le& C\sum_{K\in\TW}\Ec(F;\gamma K)^p
\sum_{\{m:~\diam K\le 2^m\}}2^{-m(p-n)}\nn\\
&\le& C\sum_{K\in\TW} \frac{\Ec(F;\gamma K)^p}{(\diam
K)^{p-n}} =2^n C\sum_{K\in\TW} |K| \left(\frac{\Ec(F;\gamma
K)}{\diam K}\right)^p. \nn\ee
Since $M_{\TW}\le N(n)$, by Lemma \reff{NT1}, $I\le
C\|F\|_{\LOP}^p$, and the proof is finished.\bx
\par Given $\ve>0$ and  $\theta\ge 1$ we let $T:\Se\to S$
denote a measurable mapping satisfying inequality \rf{TR}.
Also recall that $\TWE$ is a subfamily of $\TW$ of all
cubes with diameter at most $\ve$, see \rf{DQEP}.
\begin{lemma}\lbl{NR} For every continuous function $F\in
L_p(\RN),$ $p\in(1,\infty)$, we have
$$
\|F\circ T\|_{L_p(\Se)}^p\le
C(\|F\|_{L_p(\RN)}^p+
\sum_{Q\in\TWE}|\,Q|\,\Ec(F;\eta\,Q)^p\,)
$$
where $\eta:=10\theta+1$ and $C$  is a constant depending
only on $n,p$ and $\theta$. \end{lemma}
\par {\it Proof.} Put
$$
A_\ve:=\{Q\in\TW:~Q\cap S_\ve\ne\emp\}.
$$
Then $\dist(Q,S)\le \ve$ for every $Q\in A_\ve$ so that, by
part (ii) of Theorem \reff{Wcov}, $\diam
Q=2r_Q\le\dist(Q,S)\le\ve, $ proving that
$A_\ve\subset\TWE$.
\par Prove that
\bel{ETAQ}
T(Q)\subset \eta Q.
\ee
In fact, by \rf{TR}, for each $x\in Q$ we have
$$
\|T(x)-x_Q\|\le\|T(x)-x\|+\|x-x_Q\|\le\theta\dist(x,S)+r_Q.
$$
But
$ \dist(x,S)\le \dist(Q,S)+\diam Q, $
so that, by (ii), Theorem \reff{Wcov}, we have
$ \dist(x,S)\le 5\diam Q=10r_Q. $
Hence
$$
\|T(x)-x_Q\|\le (10\theta+1)r_Q=\eta r_Q
$$
proving \rf{ETAQ}.
\par Put $\tQ:=\eta Q$. We have
\be \intl_{Q\cap S_\ve}|F(T(x))|^p\,dx&\le& 2^p\intl_{Q\cap
S_\ve}(|F(T(x))-F_{\tQ}|^p
+|F_{\tQ}|^p)\,dx\nn\\
&\le& 2^p\intl_{Q\cap S_\ve}|F(T(x))-F_{\tQ}|^p\,dx +
2^p|Q||F_{\tQ}|^p. \nn\ee
Since $T(x)\in \tQ$, we obtain
$$
\intl_{Q\cap S_\ve}|F(T(x))-F_{\tQ}|^p\,dx\le |Q|
\sup_{x\in\tQ} |F(x)-F_{\tQ}|^p\le |Q|\Ec(F;\tQ)^p,
$$
so that
$$ \intl_{Q\cap S_\ve}|F(T(x))|^p\,dx\le 2^p(|Q|
\Ec(F;\tQ)^p +|Q||F_{\tQ}|^p). $$
This inequality and Lemma \reff{LMF} yield
\be \intl_{\Se\setminus S}|F(T(x))|^p\,dx&\le& \sum_{Q\in
A_\ve}~\intl_{Q\cap S_\ve}|F(T(x))|^p\,dx\nn\\&\le&
2^p\left(\sum_{Q\in\TWE}|Q|\, \Ec(F;\tQ)^p
+\sum_{Q\in\TW}|Q||F_{\tQ}|^p\right) \nn \\&\le&
C\left(\sum_{Q\in\TWE}|Q|\, \Ec(F;\tQ)^p
+\|F\|^p_{L_p(\RN)}\right).\nn \ee
It remains to note that $T(x)=x, x\in S,$ so that
\be \intl_{\Se}|F(T(x))|^p\,dx&=& \intl_{S}|F(T(x))|^p\,dx+
\intl_{\Se\setminus S}|F(T(x))|^p\,dx\nn\\&=&
\intl_{S}|F(x)|^p\,dx+ \intl_{\Se\setminus
S}|F(T(x))|^p\,dx\nn\ee
and the lemma follows.\bx
\begin{proposition}\lbl{N2} For every continuous function
$F\in\WP, p>n,$ we have
$$
\|F\circ T\|_{L_p(\Se)}\le C\|F\|_{\WP}
$$
where $C$  is a constant depending only on $n,p$ and
$\theta$.
\end{proposition}
\par {\it Proof.} Recall that $\diam Q\le \ve$ for every
$Q\in\TWE$. Hence, by Lemma \reff{NT1}, we have
$$
\sum_{Q\in\TWE}|Q|\,\Ec(F;\eta Q)^p\le \ve^p
\sum_{Q\in\TWE}|Q|\,\left(\frac{\Ec(F;\eta Q)}{\diam Q}
\right)^p\le C\|F\|_{\LOP},
$$
so that by Lemma \reff{NR}
$$
\|F\circ T\|_{L_p(\Se)}^p\le
C(\|F\|_{L_p(\RN)}^p+\sum_{Q\in\TWE}|Q|\Ec(F;\eta Q)^p)
\le C(\|F\|_{L_p(\RN)}^p+\|F\|_{\LOP})
$$
proving the proposition.\bx

\SECT{4. Local approximation properties of the extension
operator.}{4}
\indent
\par Let $\dlt>0, \bc\in\R,$ and let $f$ be a function
defined on a closed set $S\subset\RN$. We assign every
Whitney's cube $Q\in\TW$ a point $a_Q\in S$ which is
nearest to $Q$ on $S$ in the metric $\|\cdot\|$.  We put
\bel{CQ} c_Q:=\left \{
\begin{array}{ll}
f(a_Q),& \diam Q\le 2\dlt,\\\\
\bc,& \diam Q>2\dlt,
\end{array}
\right. \ee
and define a linear extension operator $\Es$ by letting
\bel{ExtOp} \Es f(x):=\left \{
\begin{array}{ll}
f(x),& x\in S,\\\\
\sum\limits_{Q\in \TW} c_Q\varphi_Q(x),& x\in \RN\setminus
S.
\end{array}
\right. \ee
\par Here $\{\varphi_Q:Q\in \TW\}$ is a smooth partition of
unity subordinated to the Whitney decomposition $\TW$, see,
e.g. \cite{St}. Recall main properties of this partition:
\par (a). $\phi_Q\in C^\infty(\RN)$ and $0\le\varphi_Q\le
1$ for every $Q\in \TW$;\smallskip
\par (b). $\supp \varphi_Q\subset Q^*(:=\frac{9}{8}Q),$
$Q\in \TW$;\smallskip
\par (c). $\sum\{\varphi_Q(x):~Q\in \TW\}=1$ for every
$x\in\RN\setminus S$;\smallskip
\par (d). for every cube $Q\in \TW$ we have
$$
\sup_{\RN}\|\nabla\varphi_Q\| \le C/\diam Q,
$$
where $C$ is a constant depending only on $n$.\smallskip
\par We fix $\dlt$ and $\bc$ and put
$$
\tf:=\Es f.
$$
\par In this section we present estimates of local
oscillations of $\tf$ via those of the function $f$. These
estimates are based on a series of auxiliary lemmas. To
formulate the first of them, given a cube $K\subset\RN$ we
define two families of Whitney's cubes:
$$
\Qc_1(K):=\{Q\in \TW:~Q\cap K\ne\emp\}
$$
and
\bel{AA} \Qc_2(K):=\{Q\in \TW:\exists\, Q'\in\Qc_1(K)~~{\rm
such~that~~}Q'\cap Q^*\ne\emp\}. \ee
\begin{lemma}\lbl{Cubes} Let $K$ be a cube centered
in $S$. Then for every $Q\in\Qc_2(K)$  we have $\diam Q\le
2\diam K$ and $\|x_K-x_Q\|\le \frac{5}{2}\diam K.$
\end{lemma}
\par For the proof see, e.g. \cite{S4}.
\begin{lemma}\lbl{L1} For every cube $K$ in $\RN$ and
every $c_0\in\R$ we have
$$
\sup_{K\setminus S} |\tf-c_0|\le
\sup_{Q\in\Qc_2(K)}|c_Q-c_0|.
$$
\end{lemma}
\par {\it Proof.} Clearly,
$K\setminus S\subset \cup\{Q:~Q\in\Qc_1(K)\}$ so that
$$
\sup_{K\setminus S} |\tf-c_0|\le \sup_{Q\in\Qc_1(K)}
\sup_{Q}|\tf-c_0|.
$$
\par Let $Q\in\Qc_1(K)$ and let
$ V(Q):=\{Q'\in \TW: (Q')^*\cap Q\ne\emp\}. $
By \rf{AA}, $V(Q)\subset \Qc_2(K)$.
\par Properties (a)-(c) of the partition of unity and
formula \rf{ExtOp} imply the following inequality:
\be \sup_{Q}|\tf-c_0|&=& \sup_{x\in Q}|\sum_{Q'\in \TW}
(c_{Q'}-c_0)\varphi_{Q'}(x)|\nn\\&=& \sup_{x\in Q}
|\sum_{Q'\in V(Q)} (c_{Q'}-c_0)\varphi_{Q'}(x)|\nn\\&\le&
(\max_{Q'\in V(Q)} |c_{Q'}-c_0|) (\sup_{x\in Q} \sum_{Q'\in
V(Q)} \varphi_{Q'}(x))\le \max_{Q'\in V(Q)}
|c_{Q'}-c_0|.\nn\ee
Hence,
$$
\sup_{K\setminus S} |\tf-c_0|\le \sup_{Q\in\Qc_1(K)}
\max_{Q'\in V(Q)}
|c_{Q'}-c_0|\le\sup_{Q\in\Qc_2(K)}|c_Q-c_0|.\BX
$$
\begin{lemma}\lbl{LP-N} We have
$$
\|\tf\|_{L_p(\RN\setminus S)}^p\le
C\sum_{Q\in\TWB{2\dlt}}|f(a_Q)|^p\,|Q|.
$$
\end{lemma}
\par {\it Proof.} Let $Q\in\TW$. By Lemma \reff{L1} with
$c_0=0$,
$$\sup_{Q}|\tf|\le\max\{ |c_{Q'}|:~Q'\in\Qc_2(Q)\}$$
so that
$$
\|\tf\|_{L_p(\RN\setminus S)}^p\le
\sum_{Q\in\TW}\intl_Q|\tf|^p\,dx\le\sum_{Q\in\TW} |Q|\sup_Q
|\tf|^p\le \sum_{Q\in\TW} |Q|\max_{Q'\in\Qc_2(Q)}
|c_{Q'}|^p.
$$
By Lemma \reff{Wadd} $|Q|\sim |Q'|, Q'\in\Qc_2(Q),$ so that
$$ \|\tf\|_{L_p(\RN\setminus S)}^p\le C\sum_{Q\in\TW}
\max_{Q'\in\Qc_2(Q)} |c_{Q'}|^p\,|Q'|\le
C\sum_{Q\in\TW} |c_{Q}|^p\,|Q|\card J_Q$$
where
$$ J_Q:=\{Q'\in\TW:~\Qc_2(Q')\ni Q\}.$$
By \rf{AA}, $Q'\in\Qc_2(Q)\Leftrightarrow Q\in\Qc_2(Q')$ so
that $J_Q=\Qc_2(Q)$. Moreover, by \rf{AA} and Lemma
\reff{Wadd},
$$
\card J_Q=\card\Qc_2(Q)\le(\card\Qc_1(Q))^2\le N(n)^2.
$$
Hence,
$$
\|\tf\|_{L_p(\RN\setminus S)}^p\le C\,N(n)^2\sum_{Q\in\TW}
|c_{Q}|^p|Q|.
$$
But, by \rf{CQ}, $c_Q=\bc=0$ whenever $\diam Q=2r_Q>2\dlt,$
so that
$$
\|\tf\|_{L_p(\RN\setminus S)}^p\le
CN(n)^2\sum_{Q\in\TW,\diam Q\le 2\dlt} |c_{Q}|^p|Q|=
CN(n)^2\sum_{Q\in\TWB{2\dlt}}|f(a_Q)|^p\,|Q|.
$$
\par The lemma is proved.\bx
\begin{lemma}\lbl{EKS}
\par (a). For every cube $K$ centered in $S$ we have:
$$
\sup_{K} |\tf-\bc| \le \sup_{(14K)\cap S}|f-\bc|.
$$
\par (b). If $K$ is a cube centered in $S$ and
$\diam K\le \delta$, then for every $c_0\in\R$
$$
\sup_{K} |\tf-c_0| \le \sup_{(14K)\cap S}|f-c_0|.
$$
\end{lemma}
\par {\it Proof.} It can be readily seen that Theorem
\reff{Wcov} and Lemma \reff{Cubes} imply the following:
\bel{AQK}
a_Q\in (14K)\cap S, ~~~Q\in \Qc_2(K).
\ee
\par Prove (a). By Lemma \reff{L1} with $c_0=\bc$,
$$
\sup_{K\setminus S}|\tf-\bc|\le\sup_{Q\in\Qc_2(K)}
|c_Q-\bc|.
$$
By \rf{CQ}, $c_Q$ takes the value $f(a_Q)$ or $\bc$ so that
by \rf{AQK}
$$
\sup_{K\setminus S}|\tf-\bc|\le\sup_{Q\in\Qc_2(K)}
|f(a_Q)-\bc|\le \sup_{(14K)\cap S}|f-\bc|.
$$
Finally,
$$
\sup_{K} |\tf-\bc| \le \max\{\sup_{K\cap S}
|f-\bc|,\sup_{K\setminus S}|\tf-\bc|\}\le\sup_{(14K)\cap
S}|f-\bc|,
$$
proving part (a) of the lemma.
\par Prove (b). Again, applying Lemma \reff{L1}, we obtain
$$
\sup_{K} |\tf-c_0| \le \max\{\sup_{K\cap
S}|f-c_0|,\sup_{K\setminus S}|\tf-c_0|\}\le
\max\{\sup_{K\cap S}|f-c_0|,\sup_{Q\in\Qc_2(K)}|c_Q-c_0|\}.
$$
But,  by Lemma \reff{Cubes}, $ \diam Q\le 2\diam K\le
2\dlt$\, for every $Q\in\Qc_2(K)$ so that by \rf{CQ}
$c_Q=f(a_Q)$ for each $Q\in\Qc_2(K)$. Hence, by \rf{AQK},
$$
\sup_{Q\in\Qc_2(K)} |c_{Q}-c_0| =\sup_{Q\in\Qc_2(K)}
|f(a_Q)-c_0|\le \sup_{(14K)\cap S}|f-c_0|
$$
proving part (b) of the lemma.\bx
\par We turn to estimates of local oscillations of $\tf$
on cubes which are located rather far from the set $S$. Let
$K$ be a cube  satisfying the following inequality
\bel{KA} \diam K\le\dist(K,S)/40. \ee
We let $Q_K\in \TW$ denote a Whitney's cube which contains
center of $K$, the point $x_K$.
\par The proof of the next lemma easily follows from
condition \rf{KA} and Theorem \reff{Wcov}.
\begin{lemma}\lbl{QK} $K\subset Q^*_K$ and
$$(8/9)\diam Q_K\le \dist(K,S)\le 5\diam Q_K. $$
\end{lemma}
\par We present one more property of cubes satisfying
condition \rf{KA}.
\begin{lemma}\lbl{Ksmall} For every cube $K$ satisfying
\rf{KA} the following inequality
$$
\Ec(\tf;K)\le C(n)\frac{r_K}{r_{Q_K}}\max\{|c_Q-c_{Q_K}|:~
Q\in\TW,~Q^*\cap K\ne\emp\}
$$
holds.
\end{lemma}
\par {\it Proof.} Since $K\subset\RN\setminus S$, the
function $\tf_K\in C^{\infty}(\RN)$ so that for every
$x,y\in K$ we have
$$
|\tf(x)-\tf(y)|\le C(\sup_K \|\nabla\tf\|)\,\|x-y\|\le
Cr_K\sup_K \|\nabla\tf\| .
$$
Since $K\subset Q^*_K$, see Lemma \reff{QK}, by properties
(b)-(d) of partition of unity, for every $x\in K$ we have
\be
\|\nabla\tf(x)\|&=&
\|\nabla(\sum_{Q\in\TW}(c_Q-c_K)\varphi_Q(x))\| =
\|\sum_{Q^*\cap K\ne\emptyset,Q\in\TW}
(c_Q-c_K)\nabla\varphi_Q(x)\|\nn\\
&\le& \sum_{Q^*\cap K\ne\emptyset,Q\in\TW}
|c_Q-c_K|\,\|\nabla\varphi_Q(x)\|\le C\,\sum_{Q^*\cap
K\ne\emptyset,Q\in\TW} \frac{|c_Q-c_K|}{\diam Q}.
\nn\ee
Clearly, if  $Q^*\cap K\ne\emptyset$, then $Q^*\cap Q_K
\ne\emptyset$ as well, so that by Lemma \reff{Wadd} the
number of such cubes $Q\in\TW$ is bounded by $N(n)$ and
$\diam Q\sim \diam Q_K$. Hence
$$
\sup_K\|\nabla\tf\|\le C\,r_{Q_K}^{-1}\, \max
\{|c_Q-c_K|:~Q\in\TW,~Q^*\cap K\ne\emptyset\}
$$
proving the lemma.\bx
\par Let $\pi$ be a packing of cubes in $\RN$. We put
$$
I(\tf;\pi):=\left\{\sum_{K\in\pi}\,r_K^{n-p}\,
\Ec(\tf;K)^p\right\}^\frac1p.
$$
We also define a family of cubes $\Qc(\pi)$ by letting
\bel{QPI}
\Qc(\pi):=\{Q\in\TW:~{\rm there~is~} K\in\pi
~{\rm such~ that}~~x_K\in Q\}.
\ee
\begin{lemma}\lbl{OEQ} For every packing $\pi$ of cubes
satisfying inequality \rf{KA} we have
$$
I(\tf;\pi)^p\le
C(n,p)\sum_{Q\in\Qc(\pi)}r_{Q}^{n-p}\max\{|c_Q-c_{Q'}|^p:
Q'\in\TW, Q\cap Q'\ne\emptyset\}.
$$
\end{lemma}
\par {\it Proof.} Fix a cube $\tQ\in\TW$ and suppose that
all cubes $K\in\pi$ are centered in $\tQ$. By Lemma
\reff{Ksmall} for each $K\in\pi$ we have
$$
\Ec(\tf;K)\le
C\frac{r_K}{r_{\tQ}}\max\{|c_Q-c_{\tQ}|:~Q\in\TW, Q^*\cap
K\ne\emp\}.
$$
Since $x_K\in\tQ$, by Lemma \reff{QK}, $K\subset\tQ^*$ so
that $Q^*\cap \tQ^*\ne\emp$ for every $Q\in\TW$ such that
$Q^*\cap K\ne\emp$. Then, by part (3) of Lemma \reff{Wadd},
$Q\cap \tQ\ne\emp$ as well so that
$ \Ec(\tf;K)\le C r_K J(\tf;\tQ)/r_{\tQ} $
where
$$
J(\tf;\tQ):= \max\{|c_Q-c_{\tQ}|:~Q\in\TW,~Q\cap
\tQ\ne\emptyset\}.
$$
Hence
$$
I(\tf;\pi)^p\le
C\sum_{K\in\pi}\,r_{\tQ}^{-p}J(\tf;\tQ)^p
\,r_K^n=C\,r_{\tQ}^{-p}J(\tf;\tQ)^p\,
\sum_{K\in\pi}|K|\,.
$$
Since $\pi$ consists of disjoint cubes lying in $\tQ^*$, we
obtain
$$
\sum_{K\in\pi}|K|=|\bigcup_{K\in \pi} K|\le |\tQ^*|\le C\,
|\tQ|=C r^n_{\tQ}
$$
proving that in the case under consideration
$ I(\tf;\pi)^p\le C\,r_{\tQ}^{n-p}J(\tf;\tQ). $
\par Now let $\pi$ be an {\it arbitrary} packing of cubes
satisfying inequality \rf{KA}. Then
\be I(\tf;\pi)^p &\le& \sum_{Q\in\TW}~\sum_{K\in\pi, x_K\in
Q}r_K^{n-p} \Ec(\tf;K)^p\nn\\
&\le&
C\sum_{Q\in\Qc(\pi)}r_{Q}^{n-p}\max\{|c_Q-c_{Q'}|^p:
Q'\in\TW, Q\cap Q'\ne\emptyset\}.\BX \nn\ee
\par We let $\KRN$ denote the family of all cubes in $\RN$.
We put
\bel{K1}
\Kc_1:=\{K\in\KRN:~\dist(K,S)<40\diam K,\diam
K\le\dlt/82\},
\ee
\bel{K2}
\Kc_2:=\{K\in\KRN:~\dist(K,S)<40\diam K,\diam
K>\dlt/82\},
\ee
\bel{K3}
\Kc_3:=\{K\in\KRN:~40\diam
K\le\dist(K,S)\le(2/5)\dlt\},
\ee
\bel{K4}
\Kc_4:=\{K\in\KRN:~40\diam
K\le\dist(K,S),(2/5)\dlt<\dist(K,S)\le 40\dlt\},
\ee
and, finally,
\bel{K5}
\Kc_5:=\{K\in\KRN:~40\diam
K\le\dist(K,S),\dist(K,S)>40\dlt\}.
\ee
Clearly, subfamilies $\{\Kc_i,~i=1,...,5\}$ provide a
partition of the family $\KRN$. Let us  estimate
$I(\tf;\pi)$ for $\pi\subset \Kc_i, i=1,...,5.$
\par We begin with the estimates of the quantity
$I(\tf;\pi)$ for $\pi\subset\Kc_1$.
\begin{lemma}\lbl{P1} For every packing
$\pi\subset\Kc_1$ of equal cubes of diameter $t$ there
exists a packing $\tpi=\{Q\}$ of equal cubes centered in
$S$ of diameter $\gamma t$ such that
$$
I(\tf;\pi)^p\le C(n,p)\,\sum_{Q\in\tpi}\,r_Q^{n-p}
\Ec(f;Q\cap S)^p.
$$
Here $\gamma(=574)$ is an absolute constant.
\end{lemma}
\par {\it Proof.} For each $K\in\pi$ we let $K'$ denote
the cube
$$
K':=Q(a_K,\diam K+\dist(K,S)).
$$
(Recall that $a_K$ is a point nearest to $K$ on $S$.) It
can be easily seen that $ K'\supset K$. Observe that for
every cube $K\in\Kc_1$ we have
\bel{KW}
r_{K'}:=\diam K+\dist(K,S)\le \diam K+40\diam K
=41\diam K\le \dlt/2.
\ee
Let us apply part (b) of Lemma \reff{EKS} to $K'$ (with
$c_0:=\tf(a_K)=f(a_K)$). We obtain
\be \Ec(\tf;K')&:=&\sup_{x,y\in K'}|\tf(x)-\tf(y)|\le
2\sup_{K'}|\tf-\tf(a_K)|\nn\\&\le& 2\sup_{(14K')\cap S}
|f-f(a_K)|\le 2\,\Ec(f;(14K')\cap S). \nn \ee
Since $K\subset K'$, we have
\bel{EA} I(\tf;\pi)^p\le \sum_{K\in\pi}\,r_K^{n-p}
\Ec(\tf;K')^p\le C^p\sum_{K\in\pi}\,r_{K'}^{n-p} \Ec(f;(14
K')\cap S)^p. \ee
\par Given $K\in\pi$ we put $\tK:=Q(a_K,\gamma\,t)$ where
$\gamma=14\cdot 41=574$.  Clearly, by \rf{KW}, $\tK\supset
14\,K'\supset K.$
\par We also put $\tpi:=\{\tK:~K\in\pi\}.$ Since
$r_{\tK}\sim r_K$ and $\tK\supset K,~K\in\pi$, by Lemma
\reff{FM} the packing $\tpi$ can be represented as union of
at most $C(n)$ packings. Therefore, without loss of
generality, we may suppose that $\tpi$ is a {\it packing}.
\par Finally, \rf{EA} and inclusion $\tK\supset K$ imply
the following:
$$
I(\tf;\pi)^p\le C^p\sum_{K\in\pi}\,r_{K'}^{n-p} \Ec(f;(14
K')\cap S)^p\le C^p\sum_{K\in\pi}\,r_{\tK}^{n-p}
\Ec(f;\tK\cap S)^p.
$$
\par The lemma is proved.\bx
\par This lemma and Definition \reff{DEF-AP} imply the
following
\begin{lemma}\lbl{P1AP} For every packing
$\pi\subset\Kc_1$ of equal cubes of diameter $t$ we have
$$ I(\tf;\pi)\le C(n,p)\,\AP(f;\gamma t:S)/t. $$
Here  $\gamma(=574)$ is an absolute constant.
\end{lemma}
\par We turn to estimates of oscillations of $\tf$ over
cubes from the family $\Kc_3$. (We do not present here such
estimates for the family $\Kc_2$. As we shall see below,
these estimates are either elementary or can be easily
avoided.)
\begin{lemma}\lbl{P3} For every packing $\pi\subset\Kc_3$
of equal cubes of diameter $t$ there exists a finite family
$\pi'=\{Q\}\subset\TW$ of Whitney's cubes of diameter
$2t\le\diam Q\le 2\dlt$ such that
$$
I(\tf;\pi)^p\le C(n,p)\,\sum\{r_{Q}^{n-p}|f(a_Q)-f(a_{Q'})|^p:
Q,Q'\in\pi', Q\cap Q'\ne\emptyset,r_{Q'}\le r_Q\}.
$$
\end{lemma}
\par {\it Proof.} By Lemma \reff{OEQ}
\bel{EIQ} I(\tf;\pi)^p\le
C\sum_{Q\in\Qc(\pi)}r_{Q}^{n-p}\max\{|c_Q-c_{Q'}|^p:
Q'\in\TW, Q\cap Q'\ne\emptyset\}, \ee
where $\Qc(\pi)$ is defined by \rf{QPI}. Thus for every
$Q\in\Qc(\pi)$ there is a cube $K\in\pi$ such that $x_K\in
Q$. Since $K\in\Kc_3$, we have $40\diam K\le\dist(K,S)$ and
$\dist(K,S)\le(2/5)\dlt.$ Clearly,
$$
\dist(K,S)\le \diam Q+\dist(Q,S)\le 5\diam Q,
$$
so that
$$40t=40\diam K\le\dist(K,S)\le 5\diam Q.$$
Hence $8t\le\diam Q$. If $Q'\cap Q\ne\emp$, $Q'\in\TW$,
then, by Lemma \reff{Wadd}, $\diam Q\le 4\diam Q'$ so that
$\diam Q'\ge 2t$.
\par On the other hand, by Lemma \reff{QK}, $\diam
Q\le(9/8)\dist(K,S)$ so that
$$
\diam Q\le(9/8)(2/5)\dlt=(9/20)\dlt.
$$
Then, by Lemma \reff{Wadd}, for every $Q'\in\TW$ such that
$Q'\cap Q\ne\emp$ we have
$$ \diam Q'\le 4\diam Q\le(9/5)\dlt<2\dlt, $$
or, equivalently, $r_{Q'}\le \dlt.$ Therefore, by \rf{CQ},
$c_Q=f(a_Q),~~c_{Q'}=f(a_{Q'}).$ (Recall that $a_Q$ stands
for a point nearest to $Q$ on $S$.)
\par Finally, we put
$$
\pi':=\{Q\in\TW:~{\rm there~is}~~Q'\in\Qc(\pi)~~{\rm such~
that}~~ Q\cap Q'\ne\emptyset\}.
$$
As we have proved, $2t\le\diam Q\le 2\dlt$ for every
$Q\in\pi'$.
\par Now, by \rf{EIQ},
\be I(\tf;\pi)^p&\le&
C\sum_{Q\in\Qc(\pi)}r_{Q}^{n-p}\max\{f(a_Q)-f(a_{Q'})|^p:
Q'\in\TW, Q\cap Q'\ne\emptyset\}\nn\\&\le&
C\,\sum\{r_{Q}^{n-p}|f(a_Q)-f(a_{Q'})|^p: Q,Q'\in\pi',Q\cap
Q'\ne\emptyset\}\nn\\&\le&
2C\,\sum\{r_{Q}^{n-p}|f(a_Q)-f(a_{Q'})|^p:
Q,Q'\in\pi',Q\cap Q'\ne\emptyset,r_{Q'}\le r_Q\}\nn\ee
proving the lemma.\bx
\begin{lemma}\lbl{P3A} For every packing $\pi\subset\Kc_3$
there exists a packing $\pcw\subset\TW$ of Whitney's cubes
of diameter at most $2\dlt$ such that
\bel{LP3IN} I(\tf;\pi)^p\le
C(n,p)\sum_{Q\in\pcw}\,r_{Q}^{n-p} \Ec(f;(11Q)\cap S)^p.
\ee
\end{lemma}
\par {\it Proof.} By Lemma \reff{P3} there exists a finite
family $\pi'\subset\TW$ of Whitney's cubes of diameter at
most $2\dlt$ such that
$$
I(\tf;\pi)^p\le C\sum\{r_{Q}^{n-p}|f(a_Q)-f(a_{Q'})|^p:
Q,Q'\in\pi', Q\cap Q'\ne\emptyset,r_{Q'}\le r_Q\}.
$$
Observe that for every two cubes $Q,Q'\in\TW$ such that
$Q\cap Q'\ne\emptyset, r_{Q'}\le r_Q,$ we have
$a_Q,a_{Q'}\in (11Q)\cap S$. This property easily follows
from Theorem \reff{Wcov} and Lemma \reff{Wadd}. Hence,
$$ |f(a_Q)-f(a_{Q'})|\le \Ec(f;(11Q)\cap S) $$
so that
$$
I(\tf;\pi)^p\le C\sum_{Q\in\pi'}
\,r_{Q}^{n-p}\Ec(f;(11Q)\cap S)^p\,\card \{Q'\in\TW:~
Q'\cap Q\ne\emptyset\}.
$$
By Lemma \reff{Wadd}, (2), $\card \{Q'\in\TW:~ Q'\cap
Q\ne\emptyset\}\le N(n)$ so that
\bel{RTF} I(\tf;\pi)^p\le C\sum_{Q\in\pi'}
\,r_{Q}^{n-p}\Ec(f;(11Q)\cap S)^p. \ee
\par Since $\pi'\subset\TW$, its covering multiplicity is
bounded by a constant $N(n)$, see Theorem \reff{Wcov}, so
that, by Theorem \reff{TFM}, $\pi'$ can be partitioned into
$C(n)$ packings $\{\pi'_i:~i=1,...,C(n)\}$. We let $\pcw$
denote a packing $\pi'_i$ for which the sum
$$
\sum_{Q\in\pi'_i} r_{Q}^{n-p}\Ec(f;(11Q)\cap S)^p
$$
is maximal. Clearly, by inequality \rf{RTF}, $\pcw$
satisfies inequality \rf{LP3IN}, and the lemma follows. \bx
\par In the next three lemmas we present estimates of local
oscillation of $\tf$ over cubes from $\Kc_3$ via the
functional $\PRM(f;\cdot:S)$, see \rf{OMP}.
\begin{lemma}\lbl{APCQ} Let $0<\beta<1,$
$0<t<T,\varsigma\ge 1$, $0\le s\le 1$, and let
$\pi=\{Q\}\subset\TW$ be a collection of Whitney's cubes of
diameter $t\le\diam Q\le T$. Given $Q\in\pi$ let $\tQ$ be a
$\beta$-porous cube centered in $\partial S$ such that
\bel{XD}\|x_Q-x_{\tQ}\|\le\varsigma\diam Q ~~~~{\rm
and}~~~~ \diam Q\le\diam \tQ\le \varsigma\diam Q. \ee
\par Then for every $\alpha\in(0,\beta)$ and every $f\in
C(S)$ we have
$$
\sum_{Q\in\pi}r_{\tQ}^{n-sp} \Ec(f;\tQ\cap\partial S)^p \le
C\intl_{t}^{4\varsigma T}
\PRM(f;u:S)^p \,\frac{du}{u^{sp+1}}.
$$
Here $C$ is a constant depending only on
$n,p,\varsigma,\beta$ and $\alpha$.
\end{lemma}
\par {\it Proof.} Since $\PRM$ is a non-increasing function
of $\alpha$, it suffices to prove the statement of the
lemma for $\alpha\in(\beta/2,\beta)$. We put
$\tau:=\beta/\alpha.$ Since $\alpha<\beta$, $\tau>1$.
\par Given integer $m$ consider a family of cubes
$ \tpi_m:=\{Q\in\pi:~\tau^{m-1}<r_{\tQ}\le\tau^m\}. $
Since for each $Q\in\pi$
$$
\diam\tQ\le\varsigma \diam Q~~~{\rm and}~~~t\le\diam Q\le
T,
$$
we have
$$
t\le\diam\tQ=2r_{\tQ}\le \varsigma \diam Q\le\varsigma
T,~~~~ Q\in\pi.
$$
Therefore we may assume that $t\le 2\tau^m$ and
$2\tau^{m-1}\le \varsigma T$.  Thus $M_0\le m\le M_1$ where
$M_0,M_1$ are positive integers such that
$ 2\tau^{M_1-1}\le \varsigma T<2\tau^{M_1} $
and
$ 2\tau^{M_0-1}\le t\le 2\tau^{M_0}. $
\par Given $Q\in\tpi_m$ we put
$$
Q^\sharp:=(\tau^m/r_{\tQ})\tQ=Q(x_{\tQ},\tau^m).
$$
Then $\tQ\subset Q^\sharp\subset \tau\tQ=(\beta/\alpha)\tQ$
so that
\bel{QAV}
\tQ\subset Q^\sharp\subset 2\tQ.
\ee
\par Since $\tQ$ has porosity $\beta$, it contains a cube
of diameter at least $\beta\diam\tQ$ which lies in
$\RN\setminus S$. Put $Q^\sharp:=(\tau^m/r_{\tQ})\tQ$. Then
the set $Q^\sharp\setminus S$ contains a cube of diameter
at least
$$
(r_{\tQ}/\tau^m)\beta\diam
Q^\sharp>(\tau^{m-1}/\tau^m)\beta\diam Q^\sharp =\beta/\tau\diam
Q^\sharp=\alpha\diam Q^\sharp,
$$
so that  $Q^\sharp$ is an $\alpha$-porous cube.
\par Thus the family of cubes
$$
\pi^{\sharp}_m:=\{Q^\sharp: Q\in\tpi_m\}
$$
is a family of equal cubes of porosity $\alpha$. By \rf{XD}
and by Lemma \reff{FM}, the family $\pi^\sharp_m$ can be
partitioned into at most $C(n,\tau)=C(n,\alpha,\beta)$
packings. Therefore, without loss of generality, later on
we can assume that $\pi^\sharp_m$ is a {\it packing}.
\par Thus $\pi^\sharp_m$ is a packing of cubes centered in
$\partial S$ with diameter $2\tau^m$ and porosity $\alpha$.
Since $\tQ\subset Q^\sharp$, we have
\bel{ET-Q}
\Ec(f;\tQ\cap \partial S)\le
\Ec(f;Q^\sharp\cap\partial S),
\ee
so that
\be I&:=& \sum_{Q\in\pi}r_{\tQ}^{n-sp} \Ec(f;\tQ\cap
\partial S)^p\le C\,\sum_{Q\in\pi}r_{Q^\sharp}^{n-sp}
\Ec(f;Q^\sharp\cap \partial S)^p\nn\\
&=& C\,\sum_{m=M_0}^{M_1}\tau^{-spm}\sum_{Q\in\pi^\sharp_m}
|Q^\sharp|\,\Ec(f;Q^\sharp\cap\partial S)^p. \nn \ee
Hence, by \rf{OMP}, we have
$$
I\le C\,\sum_{m=M_0}^{M_1} \tau^{-spm}\PRM(f;2\tau^m:S)^p.
$$
Since $\PRM(f;\cdot:S)$ is a non-decreasing function, it
can  be  readily shown that
$$
I\le C\, \intl_{2\tau^{M_0}}^{2\tau^{M_1+1}}
\PRM(f;u:S)^p \,\frac{du}{u^{sp+1}}.
$$
It remains to note that $t\le 2\tau^{M_0},$ and
$$
2\tau^{M_1+1}\le \tau^2\varsigma
T=(\beta/\alpha)^2\varsigma T\le 4\varsigma T,
$$
and the lemma follows.\bx
\par The following lemma easily follows from Theorem
\reff{Wcov} and Lemma \reff{Wadd} and we leave the details
of its proof to the reader.
\begin{lemma}\lbl{QP} Let $Q,Q'\in\TW, Q\cap Q'\ne\emp$,
and let $r_{Q'}\le r_Q$. Then
\par (i). $ a_{Q'}\in\tQ:=Q(a_Q,10\diam Q);$
\smallskip
\par (ii). For every $\beta\in(0,3/20)$ and every $\eta\in
(0,1]$ the cube $\eta \tQ$ is $\beta$-porous;
\smallskip
\par (iii). $3Q\subset \tQ\setminus S.$
\end{lemma}
\begin{lemma}\lbl{AP3} Let $0<\alpha<3/20$. Then for every
packing $\pi\subset\Kc_3$ of equal cubes of diameter $t$
the following inequality
$$
I(\tf;\pi)^p\le C\intl_{t}^{\gamma\delta}
\PRM(f;u:S)^p \,\frac{du}{u^{p+1}}
$$
holds. Here $C=C(n,p,\alpha)$ and $\gamma=160$.
\end{lemma}
\par {\it Proof.} By Lemma \reff{P3}, there exists a finite
family $\pi'=\{Q\}\subset\TW$ of cubes of diameter
$2t\le\diam Q\le 2\dlt$ such that
$$
I(\tf;\pi)^p\le C^p\sum\{r_{Q}^{n-p}|f(a_Q)-f(a_{Q'})|^p:
Q,Q'\in\pi', Q\cap Q'\ne\emptyset,r_{Q'}\le r_Q\}.
$$
\par Fix $Q\in\pi'$. Then, by Lemma \reff{QP}, $a_{Q'}\in
\tQ:=Q(a_Q,10\diam Q)$ for every cube $Q'\in\pi'$ such that
$Q'\cap Q \ne\emp$ and  $r_{Q'}\le r_Q$. Clearly,
$a_{Q},a_{Q'}\in\partial S$ so that
$$
|f(a_Q)-f(a_{Q'})|\le \Ec(f;\tQ\cap \partial S).
$$
\par Hence
\be I(\tf;\pi)^p&\le& C^p\sum\{r_{Q}^{n-p}\Ec(f;\tQ\cap
\partial S)^p:
Q,Q'\in\pi', Q\cap Q'\ne\emptyset,r_{Q'}\le r_Q\}\nn\\
&\le& C^p\sum_{Q\in\pi'}r_{\tQ}^{n-p}\Ec(f;\tQ\cap \partial
S)^p \card\{Q': Q'\in\TW, Q'\cap Q\ne\emptyset\}, \nn\ee
so that by Lemma \reff{Wadd}
$$
I(\tf;\pi)^p\le C^p\sum_{Q\in\pi'}r_{\tQ}^{n-p}
\Ec(f;\tQ\cap\partial S)^p.
$$
\par Let $\beta:=\min\{2\alpha,(\alpha+3/20)/2\}$. Then
$\alpha<\beta<3/20.$ Since $20\beta<3$, by Lemma \reff{QP},
$ (20\beta)Q\subset \tQ\setminus S,$ so that the cube
$\tQ=Q(a_Q,20r_Q)$ is $\beta$-porous.
\par Now, applying Lemma \reff{APCQ} with $s=1$,
$T=2\delta$ and $\varsigma=20$, we have
$$
I(\tf;\pi)^p\le C\intl_{t}^{\gamma\delta}
\PRM(f;u:S)^p \,\frac{du}{u^{p+1}}
$$
with $\gamma=2\varsigma\beta^2/\alpha^2$. But
$\beta/\alpha\le 2$ so that we can put $\gamma=160$, and
the proof of the lemma is finished.\bx
\par We turn to the family $\Kc_4$.
\begin{lemma}\lbl{P4} For every packing $\pi\subset\Kc_4$
there exists a packing $\pi'=\{Q\}$ of Whitney's cubes with
$0.01\dlt\le r_Q\le 90\dlt$ such that
\bel{PWQ} I(\tf;\pi)^p\le
C(n,p)\,\sum_{Q\in\pi'}r_{Q}^{n-p}\,|f(a_Q)-\bc|^p. \ee
\end{lemma}
\par {\it Proof.} Let $Q\in\Qc(\pi)$, i.e., there exists a
cube $K\in\pi$ centered in $Q$. Since $K\in\Kc_4$, by
\rf{K4}, $40 \diam K\le\dist(K,S)$
and $(2/5)\dlt<\dist(K,S)\le 40\dlt$. Also, by Lemma
\reff{QK},
$$
(1/5)\dist(K,S)\le \diam Q \le (9/8)\dist(K,S)
$$
so that $(2/25)\dlt\le \diam Q \le 45\dlt.$ Therefore, by
Lemma \reff{Wadd}, for each $Q'\in\TW, Q'\cap Q\ne\emp,$ we
have $0.02\dlt\le \diam Q' \le 180\dlt,$ or, equivalently,
$0.01\dlt\le r_{Q'}\le 90\dlt.$
\par Now, by Lemma \reff{OEQ}, we obtain
$$
I(\tf;\pi)^p\le C\,\sum_{Q\in\Qc(\pi)}r_{Q}^{n-p}
\max\{|c_Q-c_{Q'}|^p: Q'\in\TW, Q\cap Q'\ne\emptyset\}.
$$
\par We put
$$
\pi':=\{Q\in\TW:~{\rm there~is}~~Q'\in\Qc(\pi)~~{\rm such~
that}~~ Q\cap Q'\ne\emptyset\}.
$$
Then $\pi'=\{Q\}$ is a finite family of cubes such that
$0.01\dlt\le r_Q\le 90\dlt$ for every $Q\in\pi'$. Moreover,
since $\pi'\subset\TW$, its covering multiplicity is
bounded by a constant $N(n)$. Therefore, by Theorem
\reff{TFM}, $\pi'$ can be partitioned into at most $C(n)$
packings so that, without loss of generality, we may assume
that $\pi'$ itself is a packing .
\par Hence,
$$
I(\tf;\pi)^p\le
C\,\sum\{(r_{Q}+r_{Q'})^{n-p}|c_Q-c_{Q'}|^p: Q,Q'\in\pi',
Q\cap Q'\ne\emp\}.
$$
Since $r_Q\sim r_{Q'}$ for every $Q,Q'\in\TW,Q\cap
Q'\ne\emp$ (see Lemma \reff{Wadd}), we have
$$ (r_{Q}+r_{Q'})^{n-p}|c_Q-c_{Q'}|^p\le
C(r_{Q}^{n-p}|c_Q-\bc|^p+r_{Q'}^{n-p}|c_{Q'}-\bc|^p)$$
so that
\be I(\tf;\pi)^p&\le& C\,\sum\{r_{Q}^{n-p}|c_Q-\bc|^p+
r_{Q'}^{n-p}|c_{Q'}-\bc|^p:Q,Q'\in\pi', Q\cap Q'\ne\emp\}
\nn\\
&\le& C\,\sum_{Q\in\pi'}(r_{Q}^{n-p}|c_Q-\bc|^p)
\card\{Q'\in\TW: Q'\cap Q\ne\emp\}. \nn \ee
Hence, by Lemma \reff{Wadd}, we obtain
$$
I(\tf;\pi)^p\le C\,\sum_{Q\in\pi'}r_{Q}^{n-p}|c_Q-\bc|^p.
$$
By \rf{CQ}, $c_Q\in \{f(a_Q),\bc\}$ for every $Q\in\TW$ so
that
$$
|c_Q-\bc|\le |f(a_Q)-\bc|,~~~~ Q\in\TW
$$
proving that $\pi'$ satisfies inequality \rf{PWQ}.  \bx
\par This lemma implies the following
\begin{lemma}\lbl{P4C0} Let $\bc=0$. Then for every packing
$\pi\subset\Kc_4$ there exists a packing $\pi'=\{Q\}\subset
\TW$ with $0.01\dlt\le \diam Q\le 90\dlt$ such that the
following inequality
$$
I(\tf;\pi)^p\le C(n,p)\sum_{Q\in\pi'}|f(a_Q)|^p|Q|
$$
holds.
\end{lemma}
\par It remains to consider the case $\pi\subset\Kc_5$.
\begin{lemma}\lbl{P5} For every packing $\pi\subset\Kc_5$
we have $I(\tf;\pi)=0.$
\end{lemma}
\par {\it Proof.} Let $Q\in\Qc(\pi)$, i.e., there exists a
cube $K\in\pi$ centered in $Q$. Since $K\in\Kc_5$, by
\rf{K5}, $ 40 \diam K\le\dist(K,S)$ and $\dist(K,S)>
40\dlt$. On the other hand, by Lemma \reff{QK}, $\diam Q\ge
(1/5)\dist(K,S)$ so that $ \diam Q>(40/5)\dlt=8\dlt.$
\par Therefore, by Lemma \reff{Wadd}, for each $Q'\in\TW,
Q'\cap Q\ne\emp,$ we have
$$
\diam Q'\ge (1/4)\diam Q> 2\dlt.
$$
Hence, by \rf{CQ}, $c_{Q'}=c_Q=\bc$ so that, by Lemma
\reff{OEQ}, $I(\tf;\pi)^p=0.$\bx
\par In the last lemma of this section we estimate the
$L_p$-norm of the extension $\tf$. For its formulation we
put
\bel{WSDEF} \WS=\TWB{180\dlt}:=\{Q\in\TW:~r_Q\le 90\dlt\}.
\ee
\par Let $\ve>0$ and let $T:S_\ve\to S$ be a mapping such
that $T(x)=x$ for every $x\in S$. Let $f$ be a function
defined on $S$ such that $f\circ T\in L_p(S_\ve)$.
\begin{lemma}\lbl{NLP1} For every family $\pi\subset \WS$
the following is true:
\par (i). If $\dlt\le 0.001\ve,$ then $Q\subset S_\ve$
for each $Q\in\pi$;
\par (ii). There exists a packing $\pi'\subset\pi$
such that
$$
\sum_{Q\in\pi} |f(a_Q)|^p|Q|\le C\{\|f\circ
T\|^p_{L_p(S_\ve)}+\sum_{Q\in \pi'}|Q|\sup_Q |f\circ
T-f(a_Q)|^p\}.
$$
\end{lemma}
\par {\it Proof.} For every $Q\in\pi$ we have
\be |f(a_Q)-\fq|&\le&\frac1{|Q|}\intl_Q
|f(T(x))-f(a_Q)|\,dx \nn\\&\le& \sup_{x\in
Q}|f(T(x))-f(a_Q)|=:J(f;Q) \nn\ee
so that
\be I&:=&\sum_{Q\in\pi} |f(a_Q)|^p|Q|\le
2^p\{\sum_{Q\in\pi} (|f(a_Q)-\fq|^p|Q|+
|\fq|^p|Q|)\}\nn\\&\le& 2^p\{\sum_{Q\in\pi} J(f;Q)^p
|Q|+\sum_{Q\in\pi} |\fq|^p|Q|\}=2^p\{I_1+I_2\}. \nn \ee
\par By Lemma \reff{QP}, $Q\subset Q(a_Q,10\,r_Q)$ for
every $Q\in\TW$. Since $a_Q\in S$, we have $Q\subset
S_{10r_Q}.$ But $r_Q\le 90\dlt$, so that
$$
10r_Q\le=900\dlt\le 900\cdot(0.001)\ve<\ve.
$$
Hence, $Q\subset S_\ve$ for every $Q\in \pi$ proving the
part (i) of the lemma.
\par Prove (ii). We can assume that $I_1<\infty$, otherwise
part (ii) of the lemma is trivial. In this case there
exists a {\it finite} subfamily $\pi'\subset\pi$ such that
$$
I_1:=\sum_{Q\in\pi}|Q|J(f;Q)^p\le 2\,\sum_{Q\in
\pi'}|Q|J(f;Q)^p.
$$
By Theorem \reff{TFM}, the family $\pi'$ can be partitioned
into at most $C(n)$ packings so that without loss of
generality we may assume that $\pi'$ is a {\it packing}.
\par Let us estimate the quantity $I_2$. We have
$$
I_2:= \sum_{Q\in\pi}|\fq|^p|Q|
=\sum_{Q\in\pi}|Q|\left|\frac1{|Q|}\intl_Q f(T(x))\,
dx\right|^p\le \sum_{Q\in\pi}\intl_Q |f(T(x))|^p\,dx.
$$
Put $U:=\cup\{Q:~Q\in\pi\}.$ Since $Q\subset S_\ve$,
$Q\in\pi$, we have $U\subset S_\ve$. Then
$$
I_2\le\sum_{Q\in\pi}\intl_Q |f(T(x))|^p\,dx= \intl_U
|f(T(x))|^p\sum_{Q\in\pi}\chi_Q(x)\,dx
$$
so that, by Theorem \reff{Wcov},
$$
I_2\le N(n)\intl_{U} |f(T(x))|^p\,dx\le
N(n) \intl_{S_\ve} |f(T(x))|^p\,dx=N(n)\|f\circ
T\|^p_{L_p(S_\ve)}.
$$
\par Finally, we have
$$
I\le 2^p\{I_2+I_1\} \le C\{\|f\circ
T\|^p_{L_p(S_\ve)}+\sum_{Q\in \pi'}|Q|J(f;Q)^p\}.\BX
$$
\SECT{5. Proofs of the main theorems.}{5}
\indent
\par {\bf Proof of Theorem \reff{EXT-SIMP}.} The {\it
necessity} part of the theorem follows from Lemma
\reff{NT1} so that we turn to the proof of
\par {\it Sufficiency.} Let $f$ be a function defined on
$S$. Theorem's hypothesis is equivalent to the next
statement: there exists a constant $\lambda>0$ such that
for every packing $\pi$ the following inequality
\bel{MENQ} \sum_{Q\in\pi}r_Q^{n-p}\,\Ec(f;(11Q)\cap S)^p
\le \lambda \ee
holds.
\par We have to prove that there exists a
continuous function $F\in\LOP$ such that $F|_S=f$ and
\bel{ESF}
\|F\|_{\LOP}\le C\lambda^{\frac1p}.
\ee
\par We define $F$ as follows. We fix a point $x_0\in S$
and set
$$ \delta:=\diam S~~{\rm and}~~\bc:=f(x_0). $$
Then we put
$$
F:=\Es f
$$
where $\Es$ is the extension operator defined by formulas
\rf{ExtOp}.
\par Let us prove that inequality \rf{ESF} is satisfied. By
criterion \rf{SPACK} it suffices to show that for every
packing $\pi=\{K\}$ of equal cubes the following inequality
\bel{CR}
I(F;\pi):=\left\{\sum_{K\in\pi}r_K^{n-p}\,\Ec(F;K)^p
\right\}^\frac1p \le C\lambda^\frac1p
\ee
holds.
\par We prove this inequality in $5$ steps.
\par Put $t:=\diam K, K\in\pi$.
\par {\it The first step: $\pi\subset\Kc_1$,} see\rf{K1}.
By Lemma \rf{P1}, there exist an absolute constant
$\gamma>1$ and a packing $\tpi=\{Q\}$ of equal cubes
centered in $S$ of diameter $\gamma t$ such that
$$
I(F;\pi)^p\le C\,\sum_{Q\in\tpi}\,r_Q^{n-p} \Ec(f;Q\cap
S)^p
$$
so that, by \rf{MENQ}, $I(F;\pi)\le C\lambda^\frac1p$
proving \rf{CR}.\smallskip
\par {\it The second step: $\pi\subset \Kc_2$,} see
\rf{K2}. In this case
$$t=\diam K>\dlt/82=\diam S/82,~~~~K\in\pi$$
so that we may assume that $\dlt=\diam S<\infty.$
\par We put $\tK:=Q(x_0,\diam S)$. The reader can easily
prove that in this case
\bel{IMN2}
I(F;\pi)^p\le C\,r_{\tK}^{n-p}\,\Ec(f,\tK\cap
S)^p
\ee
so that applying \rf{MENQ} with $\pi=\{\tK\}$ we obtain the
required inequality \rf{CR}.\smallskip
\par {\it The third step: $\pi\subset\Kc_3$}, see \rf{K3}.
In this case, by Lemma \reff{P3A}, there exists a packing
$\pcw\subset\TW$ of Whitney's cubes of diameter at most
$2\dlt$ such that
$$
I(F;\pi)^p\le C\sum_{Q\in\pcw}\,r_{Q}^{n-p}
\Ec(f;(11Q)\cap S)^p,
$$
so that, by \rf{MENQ}, $I(F;\pi)^p\le C\lambda$ proving
\rf{CR}.\smallskip
\par {\it The forth step. $\pi\subset\Kc_4$}, see \rf{K4}.
In this case
$$
40\diam K\le\dist(K,S)~~~{\rm and}~~~
(2/5)\dlt<\dist(K,S)\le 40\dlt,~~~~K\in\pi,
$$
so that we may suppose that $\dlt=\diam S<\infty.$
\par  By Lemma \reff{P4}, there exists a packing
$\pi'=\{Q\}$ with $r_Q\sim \dlt=\diam S$ such that
$$ I(F;\pi)^p\le
C\sum_{Q\in\pi'}r_{Q}^{n-p}\,|f(a_Q)-\bc|^p.
$$
Recall that $\bc=f(x_0)$ so that
\be
I(F;\pi)^p&\le&
C\sum_{Q\in\pi'}r_{Q}^{n-p}\,|f(a_Q)-f(x_0)|^p\le
C\sum_{Q\in\pi'}r_{Q}^{n-p}\,\Ec(f;S)^p\nn\\ &\le& C(\diam
S)^{n-p}\,\Ec(f;S)^p.\nn
\ee
Hence
\bel{IN4}
I(F;\pi)^p\le C\,r_{\tK}^{n-p} \,\Ec(f;\tK\cap
S)^p
\ee
where $\tK:=Q(x_0,\diam S).$ It remains to apply \rf{MENQ}
with $\pi=\{\tK\}$ and the required inequality
$I(F;\pi)^p\le C\lambda$ follows. \smallskip
\par {\it The fifth step.} $\pi\subset\Kc_5$. In this case,
by Lemma \reff{P5}, $I(F;\pi)=0,$ so that \rf{CR} trivially
holds.
\par Theorem \reff{EXT-SIMP} is completely proved.\bx
\smallskip
\begin{remark}\lbl{DEL-W}{\em Given $\ve\in(0,1]$ the
factor $\gamma=11$ in the formulation of Theorem
\reff{EXT-SIMP} can be decreased to $\gamma=3+\ve$. It can
be done by a slight modification of the Whitney extension
formula \rf{ExtOp} where the Whitney covering $\TW$ is
replaced with the covering  $\TWM=\{Q\}$ of $\RN\setminus
S$ satisfying the following conditions: \par (a). $\ve\diam
Q\le\dist(Q,S)\le 4\ve\diam Q,~~Q\in\TWM$;
\par (b). Every point of $\RN\setminus S$ is covered by at
most $N=N(n,\ve)$ cubes from $\TWM$.
\par Also we put $Q^*:=(1+\ve/8)Q$. (The existence of the
family $\TWM$ readily follows from Besicovitch covering
theorem, see, e.g. \cite{G}). \par We observe that the
factor $\gamma=11$ first appears in the proof of Lemma
\reff{P3A}. Ele\-men\-tary changes in this proof show that
the corresponding factor for the modified Whitney's method
does not exceed the quantity $\gamma=3+9\ve$. \par Observe
that in this case the constants of the equivalence
$\|f\|_{\LOP|_S}\sim \inf \lambda^{\frac1p}$ depend on
$n,p$ and $\ve$.} \end{remark}
\medskip
\par {\bf Proof of Theorem \reff{EXT1}}. The {\it
necessity} part of the theorem follows from inequality
\rf{N0}, Proposition \reff{NL} and inequality \rf{ADS}.
\par Let us prove the {\it sufficiency}. Let $\lambda>0$,
$0<\alpha<3/20,$ and let $f:S\to\R$ be a function
satisfying the following inequality:
\bel{HT} \sup_{0<t\le 2\,\diam S}
\left(\frac{\AP(f;t:S)}{t}\right)^p+ \intl_0^{\diam
S}\left(\frac{\PRM(f;t:S)} {t}\right)^p\, \frac{dt}{t}\le
\lambda.
\ee
\par We have to prove the existence of a continuous
function $F\in\LOP$ such that $F|_S=f$ and
$\|F\|_{\LOP|}\le C\lambda^\frac1p$ where $C$ is a constant
depending only on $n,p$ and $\alpha$.
\par We use the same extension operator as in the proof of
Theorem \reff{EXT-SIMP}, i.e., we fix a point $x_0\in S$
and put
$$
F:=\Es f
$$
where $\delta:=\diam S,$ $\bc:=f(x_0),$ and $\Es$ is the
extension operator defined by formulas \rf{ExtOp}.
\par We are needed the following lemma which easily follows
from definition \rf{OM} of the $A_p$-functional.
\begin{lemma}\lbl{BIGD} Let $p>n$. Then
$$
u^{-\frac{n}{p}}\AP(f;u:S)\le (2\diam)^{-\frac{n}{p}}
\AP(f;2\diam S:S),~~~~~u\ge 2\diam S.
$$
\end{lemma}
\par Our aim is to show that inequality \rf{CR} is
satisfied. We prove this inequality following the same
scheme as in the prove of Theorem \reff{EXT-SIMP}, i.e., in
5 steps.
\par {\it The first step: $\pi\subset\Kc_1$.} Let
$\pi\subset\Kc_1$ be a packing of equal cubes of diameter
$t$. Then,  by Lemma \reff{P1AP},
$ I(F;\pi)\le Ct^{-1}\AP(f;\gamma t:S). $
On the other hand, by Lemma \reff{BIGD},
$$
\frac{\AP(f;u:S)}{u}\le u^{\frac{n}{p}-1}
\frac{\AP(f;2\diam S:S)}{(2\diam S)^\frac{n}{p}}\le
\frac{\AP(f;2\diam S:S)}{2\diam S},
~~~u\ge 2\diam S,
$$
so that, by \rf{HT},
$$
I(F;\pi)\le C\sup_{0<u\le 2\diam S} \AP(f;u:S)/u
\le C\lambda^{\frac1p}
$$
proving \rf{CR}.
\smallskip
\par {\it The second ($\pi\subset\Kc_2$), forth
($\pi\subset\Kc_4$) and fifth ($\pi\subset\Kc_5$) steps.}
Suppose that $\pi\subset\Kc_2$ or $\pi\subset\Kc_4$. In
these cases inequality \rf{CR} easy follows from
inequalities \rf{IMN2} and \rf{IN4}. In fact, by these
inequalities in both cases
$$
I(F;\pi)^p\le C^p\,r_{\tK}^{n-p}\,\Ec(f,\tK\cap S)^p
$$
where $\tK$ is a cube centered in $S$ with $\diam
\tK=2\diam S$. Hence, by \rf{HT}
$$
I(F;\pi)\le C\,\AP(f,2\diam S:S)/(2\diam S)
\le C\,\lambda^{\frac1p}
$$
so that inequality \rf{CR} is satisfied.
\par The fifth step ($\pi\subset\Kc_5$) is trivial because
in this case $I(F;\pi)=0$.
\smallskip
\par {\it The third step: $\pi\subset\Kc_3$.} By Lemma
\reff{AP3}, for every $0<\alpha<3/20$ we have
$$
I(F;\pi)^p\le C^p\intl_{0}^{\infty}
\left(\frac{\PRM(f;u:\partial S, \alpha)}{u}\right)^p
\,\frac{du}{u}.
$$
Elementary calculations show that, this inequality, Lemma
\reff{BIGD} and \rf{HT} imply the required inequality $
I(F;\pi)\le C\lambda^{\frac1p}.$ We leave the technical
details of the proof to the reader.
\par Theorem \reff{EXT1}  is completely proved.\bx
\medskip
\par {\bf Proof of Theorem \reff{EXT-NORM}.} The {\it
necessity} part of the theorem follows from Lemma
\reff{NT1} and Proposition \reff{N2}.
\par Let us prove the {\it sufficiency}. We will be needed
the following
\begin{lemma} \lbl{SMALLQ} Let $\varsigma>0$ be a positive
number and let $F$ be a continuous function from $L_p(\RN),
p>1.$ Suppose that there exists a constant $\lambda>0$ such
that for every packing $\pi$ of equal cubes of diameter at
most $\varsigma$ the following inequality
\bel{ISM}
\sum_{Q\in\pi}r_Q^{n-p}\,\Ec(F;Q)^p \le \lambda
\ee
holds. Then $F\in\LOP$ and
$$
\|F\|_{\LOP}\le
C(\lambda^\frac1p+\varsigma^{-1}\|F\|_{L_p(\RN)}).
$$
\end{lemma}
\par {\it Proof.} Let  $\pi$ be a packing of
equal cubes. Put
$$
J(F;\pi):=\left\{\sum_{Q\in\pi}r_Q^{n-p}\,
\Ec_1(F;Q)_{L_p}^p\right\}^\frac1p
$$
where $\Ec_1(F;Q)_{L_p}$ is the normalized best
approximation in $L_p$ defined by \rf{ESP}.
\par If $\diam Q\le\varsigma$ for every $Q\in\pi$, then by
\rf{ISM}
$$
J(F;\pi)^p\le
\sum_{Q\in\pi}r_Q^{n-p}\,\Ec_1(F;Q)_{L_\infty}^p \le
\sum_{Q\in\pi}r_Q^{n-p}\,\Ec(F;Q)^p\le\lambda.
$$
\par Suppose that for every $Q\in\pi$ we have
$\diam Q=2r_Q>\varsigma$. Then
$$
J(F;\pi)^p\le
\sum_{Q\in\pi}r_Q^{n-p}\,\left(\frac{1}{|Q|}
\intl_Q|F|^p\,dx\right) =\sum_{Q\in\pi}r_Q^{-p}\,
\intl_Q|F|^p\,dx,
$$
so that
$$
J(F;\pi)^p\le
2^p\varsigma^{-p}\,\sum_{Q\in\pi}\intl_Q|F|^p\,dx \le
2^p\varsigma^{-p}\|F\|_{L_p(\RN)}^p.
$$
\par Thus
$$
\sup_{\pi}J(F;\pi)^p\le
\lambda+2^p\varsigma^{-p}\|F\|_{L_p(\RN)}^p<\infty,
$$
where $\pi$ runs over all packings of equal cubes in $\RN$.
Hence, by a result of Brudnyi \cite{Br1}, $F\in\LOP$ and
$$
\|F\|_{\LOP}\le C\sup_{\pi}J(F;\pi)\le
C(\lambda^\frac1p+\varsigma^{-1}\|F\|_{L_p(\RN)}).
$$
\par The lemma is proved.\bx
\par Theorem's hypothesis can be formulated in the
following equivalent way, c.f., \rf{MENQ}.
\par Let $\theta\ge 1$ and let $T:S_\ve\to S$ be a
measurable mapping satisfying the following condition:
\bel{TR1}  \|T(x)-x\|\le \theta \dist(x,S),~~~~~ x\in
S_\ve. \ee
\par Suppose that there exists is a constant $\lambda>0$
such that:
\par (a).
\bel{MENQN1}
\|f\circ T\|_{L_p(S_\ve)}\le\lambda^\frac1p.
\ee
\par (b).  For every packing $\pi$ of cubes lying in
$S_\ve$, we have
\bel{MENQN} \sum_{Q\in\pi}r_Q^{n-p}\,\Ec(f;(\eta Q)\cap
S)^p \le \lambda, \ee
where $\eta:=10\,\theta+1$.
\par We have to prove that there exists a continuous
function $F\in\WP$ such that $F|_S=f$ and $\|F\|_{\WP}\le
C\lambda^{\frac{1}{p}}$ where the constant $C$ depends only
on $n,p,\ve$ and $\theta$.
\par We define the function $F:=\Es f$ by the formula
\rf{ExtOp} where the parameter $\bc:=0$ and the parameter
$$ \dlt:=0.001\ve. $$
\par First we prove that
\bel{ESNF} \|F\|_{L_p(\RN)}\le C\lambda^{\frac1p}. \ee
\begin{lemma}\lbl{NLP2} Under conditions \rf{MENQN}
and \rf{MENQN1} the following inequality
$$
\sum_{Q\in\WS} |f(a_Q)|^p|Q|\le C\lambda.
$$
holds. (Recall that  $\WS=\TWB{180\dlt}$, see \rf{WSDEF}.)
\end{lemma}
\par {\it Proof.} We put $\pi=\WS$ in Lemma \reff{NLP1}.
Then, by this lemma, there exists a packing $W'\subset\WS$
of cubes lying in $S_\ve$ such that
$$
I:=\sum_{Q\in\WS} |f(a_Q)|^p|Q|\le C\{\|f\circ
T\|^p_{L_p(S_\ve)}+\sum_{Q\in W'}|Q|\sup_Q |f\circ
T-f(a_Q)|^p\}.
$$
Prove that
\bel{AQF1}
\sup_{x\in Q} |f(T(x))-f(a_Q)|\le \Ec(f;(\eta
Q)\cap S).
\ee
In fact, by \rf{ETAQ}, $T(x)\in \eta Q$. Also, it can be
readily seen that $a_Q$ (a point nearest to $Q$ on $S$)
lies in the cube $9 Q\subset \eta Q$. Thus $a_Q,T(x)\in
\eta Q$ so that
$$ |f(T(x))-f(a_Q)|\le\Ec(f;(\eta Q)\cap S),$$
proving \rf{AQF1}.
\par Since $Q\in\WS$, $r_Q\le 90\dlt$ so that
\be
I'&:=& \sum_{Q\in W'}|Q|\sup_Q |f\circ T-f(a_Q)|^p\le
\sum_{Q\in W'}|Q|\Ec(f;(\eta Q)\cap S)^p\nn\\
&\le& (90\dlt)^p \sum_{Q\in W'}r_Q^{n-p} \Ec(f;(\eta Q)
\cap S)^p. \ee
Since $W'$ is a packing of  cubes lying in $S_\ve$, by
\rf{MENQN}, $I'\le C^p\,\lambda.$
\par Combining this inequality with \rf{MENQN1}, we finally
obtain
$$
I\le C\{\|f\circ T\|^p_{L_p(S_\ve)}+I'\}\le C\lambda.\BX
$$
\par We finish the proof of inequality \rf{ESNF} as
follows. By Lemma \reff{NLP2} and Lemma \reff{LP-N} we have
$\|F\|_{L_p(\RN\setminus S)}^p \le C\lambda.$ Since
$T(x)=x$ on $S$, by inequality \rf{MENQN1},
$\|F\|_{L_p(S)}^p=\|f\|_{L_p(S)}^p\le\lambda$. Hence,
$$
\|F\|_{L_p(\RN)}^p=\|F\|_{L_p(S)}^p+
\|F\|_{L_p(\RN\setminus S)}^p\le \lambda+C\lambda,
$$
proving \rf{ESNF}.
\par Let us prove that $\|F\|_{\LOP}\le
C\lambda^{\frac1p}.$ We will follow the same scheme as in
the proof of Theorem \reff{EXT-SIMP}.
\par We put $\varsigma:=\delta/90$.
By Lemma \reff{SMALLQ} it suffices to show that for every
packing $\pi=\{K\}$ of equal cubes of diameter
\bel{DSG} \diam K=t\le \varsigma \ee
the following inequality
\bel{CR1}
I(F;\pi):=\left\{\sum_{K\in\pi}r_K^{n-p}\,\Ec(F;K)^p
\right\}^\frac1p \le C\lambda^\frac1p
\ee
holds. Similarly to the proof of inequality \rf{CR}, we
prove \rf{CR1} in $5$ steps.
\par {\it The first step: $\pi\subset\Kc_1$.} Recall that
by \rf{K1}, $\diam K\le\dlt/82$ for every cube $K\in\pi$.
By Lemma \rf{P1}, there exist a packing $\tpi=\{Q\}$ of
equal cubes centered in $S$ of diameter $\gamma t$ such
that
$$
I(F;\pi)^p\le C^p\sum_{Q\in\tpi}\,r_Q^{n-p} \Ec(f;Q\cap
S)^p.
$$
Here $\gamma=574$. Since $t=\diam K\le\dlt/82, K\in\pi,$
and $\diam Q=\gamma t, Q\in\tpi,$ we have
$$
r_Q=\diam Q/2=\gamma t/2\le
\gamma\dlt/164=574\dlt/164<4\dlt.
$$
Since $\dlt=0.001\ve$, we have $r_Q<\ve$. In addition,
every cube $Q\in\tpi$ is centered in $S$ so that $Q\subset
S_\ve$.
\par Hence, by assumption \rf{MENQN}, $I(F;\pi)^p\le
C^p\lambda$ proving \rf{CR1}.
\smallskip
\par There is no necessity to consider the case of packings
$\pi\subset \Kc_2$. In fact, in this case $\diam K>\dlt/82$
for every $K\in\Kc_2$, see \rf{K2}, so that
$$\diam K>\dlt/82>\dlt/90=\varsigma$$
which contradicts \rf{DSG}.\smallskip
\par Therefore we turn to
\par {\it The third step: $\pi\subset\Kc_3$.} This case can
be treated in the same way as in the proof of Theorem
\reff{EXT-SIMP}.
\par By Lemma \reff{P3A} there exists a packing
$\pcw\subset\TW$ of Whitney's cubes of diameter at most
$2\dlt$ satisfying inequality \rf{LP3IN}. Since
$\eta:=10\theta+1\ge 11$, we obtain
\bel{IPK}
I(F;\pi)^p\le C\sum_{Q\in\pcw}\,r_{Q}^{n-p}
\Ec(f;(\eta Q)\cap S)^p.
\ee
\par By Lemma \reff{QP}, every cube $Q\in\pcw$ lies in
$5\diam Q$-neighborhood of $S$ so that $Q\subset
S_{10\dlt}.$ Since $\dlt=0.001\ve$, the cube $Q\subset
S_{\ve}.$
\par  Combining \rf{IPK} with \rf{MENQN}, we finally obtain
the required inequality $I(F;\pi)^p\le C\lambda$.
\smallskip
\par {\it The forth step: $\pi\subset\Kc_4$.} In this case
the inequality $I(F;\pi)^p\le C\lambda$ follows from Lemma
\reff{P4C0} and Lemma \reff{NLP2}.
\smallskip
\par The remaining case $\pi\subset\Kc_5$ is trivial
because in this case $I(\tf;\pi)=0$. Thus \rf{CR1} holds
for every packing $\pi$ of equal cubes proving the
theorem.\bx
\medskip
\par{\bf Proof of Theorem \reff{MF-C}.}{(\it Necessity.)}
\par (i). Let $F\in\LOP$ be a continuous function such that
$F|_S=f$ and
\bel{HT1}
\|F\|_{\LOP}\le 2\|f\|_{\LOP|_S}.
\ee
Fix $q\in(n,p)$. Then by Sobolev-Poin\'{c}are inequality
\rf{EDIFF}, for every cube $Q=Q(x,r)$ and every $y,z\in
Q\cap S$ we have
$$
|f(y)-f(z)|/r=|F(y)-F(z)|/r\le C\left(\frac{1}{|Q|}
\intl_Q\|\nabla F\|^q(u)\,du\right)^{\frac1q}
\le C M(\|\nabla F\|^q)^{\frac1q}(x).
$$
Hence,
$$
f^\sharp_{\infty,S}(x):=\sup\{|f(y)-f(z)|/r:~r>0,\,y,z\in
Q\cap S\}\le C M(\|\nabla F\|^q)^{\frac1q}(x),~~~x\in\RN,
$$
so that, by the Hardy-Littlewood maximal theorem,
$$
\|f^\sharp_{\infty,S}\|_{L_p(\RN)}\le C
\|M(\|\nabla F\|^q)^{\frac1q}\|_{L_p(\RN)}\le
C\|\nabla F\|_{L_p(\RN)},
$$
proving that $\|f^\sharp_{\infty,S}\|_{L_p(\RN)}\le
C\|F\|_{\LOP}.$
\par Combining this inequality with \rf{HT1}, we obtain
\bel{NECA}
\|f^\sharp_{\infty,S}\|_{L_p(\RN)}\le
C\|f\|_{\LOP|_S}.
\ee
\par (ii). The required inequality
$$
\|f\circ T\|_{L_p(S_\ve)}+
\|f^\sharp_{\infty,S}\|_{L_p(S_\ve)}\le C\|f\|_{\WP|_S}
$$
immediately follows from equivalence \rf{TRN} and
inequality \rf{NECA}.
\par {\it (Sufficiency)}. \par (i). Given $\gamma\ge 1$
consider a packing $\{Q_i:~i=1,...m\}$ and points
$x_i,y_i\in(\gamma Q_i)\cap S$. Let $z\in Q_i$ and let
$Q_z:=Q(z,(\gamma+1)r_{Q_i})$. Then, clearly, $Q_z\supset
\gamma Q_i$ so that $x_i,y_i\in Q_z\cap S.$
\par Since $r_{Q_z}\sim\diam Q_i$, we have
$$
\frac{|f(x_i)-f(y_i)|^p}{(\diam Q_i)^{p-n}}
\le C \,|Q_i|\,(|f(x_i)-f(y_i)|/r_{Q_z})^p
\le C\, |Q_i|\,(f^\sharp_{\infty,S}(z))^p,~~~~z\in Q_i.
$$
Integrating this inequality (with respect to $z$) on cube
$Q_i$ we obtain
$$
\frac{|f(x_i)-f(y_i)|^p}{(\diam Q_i)^{p-n}}
\le C \intl_{Q_i}(f^\sharp_{\infty,S})^p(z)\,dz
$$
so that
\bel{MF-SF} I:=\sum_{i=1}^m\frac{|f(x_i)-f(y_i)|^p}{(\diam
Q_i)^{p-n}} \le C
\intl_{U_\pi}(f^\sharp_{\infty,S})^p(z)\,dz
\ee
where $U_\pi:=\cup_{i=1}^m Q_i.$ Hence $I\le C
\|f^\sharp_{\infty,S} \|_{L_p(\RN)}^p.$
\par We put $\gamma=11$ and obtain inequality \rf{CONSTR}
with $\lambda:=C(n,p)\|f^\sharp_{\infty,S}\|_{L_p(\RN)}^p$.
Then, by \rf{EQVF},
$$
\|f\|_{\LOP|_S}\le C\lambda^{\frac1p}
\le C\|f^\sharp_{\infty,S}\|_{L_p(\RN)}
$$
proving the sufficiency part of the theorem for the case
(i).
\par (ii). By inequality \rf{MF-SF} with $\gamma=\eta$, for
every packing $\{Q_i:~i=1,...m\}$ lying in $S_\ve$ and
every two points $x_i,y_i\in(\eta Q_i)\cap S$ we have
$$
I:=\sum_{i=1}^m\frac{|f(x_i)-f(y_i)|^p}{(\diam
Q_i)^{p-n}} \le C \intl_{U_\pi}(f^\sharp_{\infty,S})^p(z)\,dz
\le C \|f^\sharp_{\infty,S}\|_{L_p(S_\ve)}^p=:\lambda.
$$
Hence, by Theorem \reff{EXT-NORM},
$$
\|f\|_{\WP|_S}\le C(\|f\circ T\|_{L_p(S_\ve)}
+\lambda^{\frac1p})\le C(\|f\circ T\|_{L_p(S_\ve)}+
\|f^\sharp_{\infty,S}\|_{L_p(S_\ve)}).
$$
\par Theorem \reff{MF-C} is completely proved.\bx
\medskip
\par{\bf Proof of Theorem \reff{EXT2}.} The {\it necessity}
part of the theorem follows from inequality \rf{N0},
Proposition \reff{NL}, inequality \rf{ADS} and Proposition
\reff{N2}.
\smallskip
\par Prove the {\it sufficiency}.
 Let $\theta\ge 1$ and let $T:S_\ve\to S$ be a measurable
mapping such that
\bel{TS}
T(x)\in\partial S,~~~~~x\in S_\ve\setminus S,
\ee
and inequality \rf{TR1} is satisfied.
\par We put
\bel{DLT2} \bc:=0~~~~~{\rm and} ~~~~~~
\dlt:=\frac{\ve}{7200(1+\theta)}, \ee
and define the function
\bel{DEFEXT}
F:=\Es f
\ee
by the formula \rf{ExtOp}.
\par In the remaining part of the section $C$ will denote a
positive constant depending only on $n,p,\ve,\theta$ and
$\alpha$.
\begin{lemma}\lbl{NORMPL} Let $t>0$ and let
$\pi\subset\WS$ be a collection of Whitney's cubes of
diameter at least $t$. Then for every continuous function
$f$ on $S$ and every $0<\alpha<0.3/(1+\theta)$ the
following inequality
$$
\sum_{Q\in\pi}|f(a_Q)|^p|Q|\le C\left\{\|f\circ
T\|_{L_p(S_\ve)}^p+
\intl_t^{\ve/2}\PRM(f;t:S)^p\,\frac{dt}{t}\right\}
$$
holds.
\end{lemma}
\par {\it Proof.} By Lemma \reff{NLP1} each cube
$Q\in\pi$ lies in $S_\ve$ and
\bel{IP1} I:=\sum_{Q\in\pi}|f(a_Q)|^p|Q|\le C^p\{\|f\circ
T\|^p_{L_p(S_\ve)}+\sum_{Q\in \pi}|Q|\sup_Q |f\circ
T-f(a_Q)|^p\}. \ee
\par Given $Q\in\pi$, we put
$ \tQ:=Q(a_Q,10(1+\theta)r_Q). $
The reader can readily see that inequality \rf{TR1} and
property \rf{DKQ} of Whitney's cubes imply the following:
$$
T(x)\in\tQ~~{\rm for~every}~~x\in Q.
$$
\par Since $a_Q\in\partial S$ and $T(x)\in\partial S, x\in
Q$, see \rf{TS}, we have
$$
|f(a_Q)-f(T(x))|\le \Ec(f;\tQ\cap\partial S),~~~~x\in Q,
$$
so that
$$
I':=\sum_{Q\in \pi}|Q|\sup_Q |f\circ T-f(a_Q)|^p \le
\sum_{Q\in \pi}|\tQ|\Ec(f;\tQ\cap\partial S)^p.
$$
\par Let $0<\beta<3/(10+10\theta)$. Prove that $\tQ$ has
porosity $\beta$. We put $\varsigma:=10(1+\theta)\beta$ so
that $0<\varsigma<3$. Therefore by Lemma \reff{QP}
$$
Q_\varsigma :=\varsigma Q\subset Q(a_Q,10\diam Q)\setminus
S.
$$
Since $\theta\ge 1$,
$$
Q(a_Q,10\diam Q)\subset \tQ:=Q(a_Q,5(1+\theta)\diam Q),
$$
so that $Q_\varsigma \subset \tQ\setminus S.$ But
$$
\beta\diam \tQ=10\beta (1+\theta)\diam Q=\varsigma \diam
Q=\diam Q_\varsigma,
$$
proving that $\tQ$ is a $\beta$-porous cube.
\par Now put
\bel{BT} \beta:=\min\{2\alpha,(\alpha+3/(10+10\theta))/2\}.
\ee
Then $0<\alpha<\beta<3/(10+10\theta)$. In addition, the
collection $\{\tQ:~Q\in\pi\}$ is a family of $\beta$-porous
cubes  centered in $S$. We apply Lemma \reff{APCQ} with
$s=0$, $T:=90\dlt$, $\varsigma:=10(1+\theta)$ and parameter
$\beta$ defined by \rf{BT}, and obtain the following
inequality:
\bel{IP2} I'\le \sum_{Q\in \pi}r_{\tQ}^{n}
\Ec(f;\tQ\cap\partial S)^p\le C\intl_{t}^{\gamma
T}\PRM(f;u:S)^p \,\frac{du}{u}. \ee
Here $\gamma=\beta^2\varsigma/\alpha^2$. Since
$\beta/\alpha\le 2$, we have
$$
\gamma T\le 4\cdot 10(1+\theta)\cdot
90\dlt=3600(1+\theta)\dlt,
$$
so that by \rf{DLT2} $ \gamma T\le\ve/2.$ Combining this
inequality, \rf{IP1} and \rf{IP2}, we obtain the statement
of the lemma.\bx
\begin{proposition}\lbl{NORMP} Let $f\in C(S)$ and let
$\alpha\in(0,0.3/(1+\theta)).$ Then
$$
\|F\|^p_{L_p(\RN)}\le C\left\{\|f\circ
T\|_{L_p(S_\ve)}^p+ \intl_0^{\ve/2}\PRM(f;t:S)^p\,
\frac{dt}{t}\right\}.
$$
\end{proposition}
\par {\it Proof.} By Lemma \reff{LP-N} and Lemma
\reff{NORMPL},
\be \|F\|_{L_p(\RN\setminus S)}^p&\le&
\sum_{Q\in\TW,r_Q\le\dlt} |f(a_Q)|^p|Q|\le\sum_{Q\in\WS}
|f(a_Q)|^p|Q|\nn\\&\le& C\left\{\|f\circ
T\|_{L_p(S_\ve)}^p+
\intl_0^{\ve/2}\PRM(f;t:S)^p\,\frac{dt}{t}\right\}. \nn\ee
\par Since $F|_S=f$ and $T(x)=x$ on $S$, we have
$$
\|F\|_{L_p(S)}=\|f\|_{L_p(S)}=\|f\circ T\|_{L_p(S)}\le
\|f\circ T\|^p_{L_p(S_\ve)}.
$$
Finally, these estimates and equality $\|F\|_{L_p(\RN)}^p=
\|F\|_{L_p(S)}^p+\|F\|_{L_p(\RN\setminus S)}^p$ imply the
statement of the proposition.\bx
\par The next proposition provides an estimate of the
$A_p$-functional of the extension $F$ of $f$ via
corresponding $A_p$- and $\Ac_{p,\alpha}$-functionals of
$f$. To formulate the result we put
$\ve_0:=\dlt/100=\ve/(7200\cdot 100(1+\theta))$, see
\rf{DLT2}.
\begin{proposition}\lbl{EXT-AP} Let $n<p<\infty$ and let
$0<\alpha<0.3/(1+\theta)$. Then for every $f\in C(S)$ and
every $0<t\le\ve_0$ we have
$$
\AP(F;t)\le C\left\{\AP(f;\gamma t:S)+
t\left(\intl_t^{\ve}\PRM(f;u:S)^p\,
\frac{du}{u^{p+1}}\right)^{\frac1p}+t\|f\circ
T\|_{L_p(S_\ve)}\right\}.
$$
Here $\gamma$ is an absolute constant ($\gamma\le 574$).
\end{proposition}
\par {\it Proof.} Let $\pi$ be a packing of equal cubes of
diameter $s>0$. Put
$$
J(F;\pi):=\left\{\sum_{Q\in\pi}|Q|\Ec(F;Q)^p\right\}
^{\frac1p}
$$
and
$$
H(f;t):=\AP(f;\gamma t:S)+
t\left(\intl_t^{\ve}\PRM(f;u:S)^p\,
\frac{du}{u^{p+1}}\right)^{\frac1p}+t\|f\circ
T\|_{L_p(S_\ve)}
$$
(with $\gamma=574$). By \rf{DI}, we have to proof that
\bel{JHP}
J(F;\pi)\le C\,H(f;t),~~~~0<s\le t\le\ve_0.
\ee
\par Recall that
$$
I(F;\pi):=\left\{\sum_{Q\in\pi}r_Q^{n-p}\Ec(F;Q)^p\right\}
^{\frac1p}
$$
so that $J(F;\pi)=s I(F;\pi)/2.$
\par We prove \rf{JHP} in 5 steps.
\par {\it The first step: $\pi\subset\Kc_1$.} By Lemma
\reff{P1AP} in this case
$$
J(F;\pi)=s I(F;\pi)/2\le C s(Cs^{-1}\AP(f;\gamma s:S)\le
C\AP(f;\gamma s:S)
$$
with $\gamma\le 574,$  so that \rf{JHP} is satisfied.
\par {\it The second step: $\pi\subset\Kc_2$.} Recall that
$\diam Q>\delta/82$ for every $Q\in\Kc_2$. But $\diam
Q=s\le\delta/100$ for each $Q\in\pi$ so that this case can
be omitted.
\par {\it The third step: $\pi\subset\Kc_3$.} Since
$\theta\ge 1$, the parameter $\alpha\in(0,3/20)$ so that by
Lemma \reff{AP3} and inequality $\gamma\dlt\le\ve$, see
\rf{DLT2},
\be
J(F;\pi)^p&=&(s I(F;\pi)/2)^p\le C s^p\intl_{s}^{\gamma
\delta} \PRM(f;u:S)^p \,\frac{du}{u^{p+1}}\nn\\
&\le& C s^p\intl_{s}^{\ve}
\PRM(f;u:S)^p \,\frac{du}{u^{p+1}}\nn\\
&\le& C \left(s^p\intl_{s}^{t} \PRM(f;u:S)^p
\,\frac{du}{u^{p+1}}+ t^p\intl_{t}^{\ve} \PRM(f;u:S)^p
\,\frac{du}{u^{p+1}}\right)=C(I_1+I_2). \nn \ee
Since $\PRM(f;\cdot:S)$ is a non-decreasing function and
$\PRM(f;\cdot:S)\le\AP(f;\cdot:S)$, see \rf{APCM}, the item
$I_1$ is bounded by $(1/p)\AP(f;t: S)^p$ proving \rf{JHP}
in the case under consideration.
\par {\it The forth step: $\pi\subset\Kc_4$.} By Lemma
\reff{P4C0} there exists a packing of Whitney's cubes
$\pi'=\{Q\}$ with $0.01\dlt\le\diam Q\le 90\dlt$ such that
$$
I(F;\pi)^p\le C\sum_{Q\in\pi'}|f(a_Q)|^p|Q|.
$$
Since $\pi'$ is a subfamily of $\WS:=\{Q\in\TW:~\diam Q\le
90\dlt\}$ of cubes of diameter at least $\dlt/100=\ve_0$,
and $0<t\le \ve_0$, by Lemma \reff{NORMPL}
$$
I(F;\pi)^p\le C\left(\|f\circ T\|_{L_p(S_\ve)}^p+
\intl_{t}^{\ve}\PRM(f;u:S)^p\,\frac{du}{u^{p+1}}\right).
$$
\par Clearly, $J(F;\pi)=0$ whenever $\pi\subset\Kc_5$, see
Lemma \reff{P5}, and the proof is finished.\bx
\par Let us finish the proof of the sufficiency part of
Theorem \reff{EXT2}. We let $\lambda$ denote the right-hand
side of the theorem's equivalence.
\par Then by Proposition \reff{NORMP} $\|F\|_{L_p(\RN)}\le
C\lambda$. On the other hand, by Proposition \reff{EXT-AP}
$\AP(F;t)\le C\lambda\,t$ for every $0<t\le\ve_0$.
Therefore, by \rf{DI}, for every packing of equal cubes of
diameter at most $\ve_0$ the inequality \rf{ISM} of Lemma
\reff{SMALLQ} holds (with constant $C\lambda^p$ instead of
$\lambda$). Then by this lemma $F\in\LOP$ and
$$
\|F\|_{\LOP}\le C(\lambda+\ve_0^{-1}\|F\|_{L_p(\RN)})\le
C(\lambda+C\lambda)\le C\lambda.
$$
\par Theorem \reff{EXT2} is completely proved.
\par This theorem enables us to prove the following
\begin{proposition}\lbl{WDS} Let $p>n$ and let
$\Omega:=int(S)$. Then for every function $f\in C(S)$ the
following equivalence
$$ \|f\|_{\WP|_S}\sim \|f|_\Omega\|_{\WPO}+ \|f|_{\partial
S}\|_{\WP|_{\partial S}} $$
holds with constants depending only on $p$ and $n$.
\end{proposition}
\par {\it Proof.} We put
\bel{FYT} \lambda:=\|f|_\Omega\|_{\WPO}+ \|f|_{\partial
S}\|_{\WP|_{\partial S}}. \ee
\par Observe that inequality $\lambda\le 2\|f\|_{\WP|_S}$
is trivial. Prove the converse inequality.
\par Fix $\ve>0$ and a mapping $T:S_\ve\to S$ satisfying
conditions \rf{TR} and \rf{ADT}. By Theorem \reff{EXT2}, it
suffices to show that for some
$\alpha\in(0,3/(10+10\theta)]$ the following inequalities
\bel{IT1} \|f\circ T\|_{L_p(S_\ve)}\le C\lambda, \ee
\bel{IT2} \sup_{0<t\le \ve} \frac{\AP(f;t:S)}{t}\le
C\lambda \ee
and
\bel{IT3} \left\{\intl_0^\ve\PRM(f;t:S)^p\,
\frac{dt}{t^{p+1}}\right\}^{\frac{1}{p}} \le C\lambda \ee
hold with a constant $C$ depending only on $n$ and $p$.
(Recall that  $\theta$ is a constant from inequality
\rf{TR}.)
\par Let $T_S:\RN\to S$ is a mapping defined by the formula
\rf{TWQ}. Also, let us consider the corresponding mapping
$T_{\partial S}:\RN\to S$ defined by the same formula
\rf{TWQ} but for $\partial S$ rather than $S$. Clearly,
$T_{\partial S}|_{\RN\setminus \Omega}=T_S$.
\par Put
$$ I:=\|f\circ T_{\partial S}\|_{L_p((\partial S)_{4\ve})}+
\sup_{0<t\le 4\ve} \frac{\AP(f;t:\partial S)}{t}+
\left\{\intl_0^{4\ve}\PRM(f;t:\partial S)^p\,
\frac{dt}{t^{p+1}}\right\}^{\frac{1}{p}} $$
and apply Theorem \reff{EXT2} to the set $\partial S$. By
this theorem,
\bel{ELC} I\le C\|f|_{\partial S}\|_{\WP|_{\partial S}} \le
C\lambda \ee
provided
\bel{ALI} \alpha\in(0,3/(10+10\theta) \ee
where $\theta=\theta(n)$ is a constant from inequality
\rf{TR} corresponding to the mapping $T_{\partial S}$.
\par Prove that inequalities \rf{IT1}-\rf{IT3} hold for
$T=T_S$ and every $\alpha$ satisfying \rf{ALI}.
\par We begin with inequality \rf{IT1}. Since $T_{\partial
S}|_{\RN\setminus \Omega}=T_S$, we have
$$
\|f\circ T_{S}\|_{L_p(S_{\ve})}^p=
\|f\|_{L_p(\Omega)}^p+ \|f\circ T_{S}\|_{L_p((\partial
S)_\ve\setminus\Omega)}^p=
\|f\|_{L_p(\Omega)}^p+ \|f\circ T_{\partial
S}\|_{L_p((\partial S)_\ve \setminus\Omega)}^p,
$$
so that
\bel{FTD} \|f\circ T_{S}\|_{L_p(S_{\ve})}^p\le
\|f|_\Omega\|_{\WPO}^p+ \|f\circ T_{\partial
S}\|_{L_p((\partial S)_\ve)}^p.
\ee
Hence, by \rf{ELC} and \rf{FYT}, we have $\|f\circ
T_{S}\|_{L_p(S_{\ve})}\le \lambda^p+(C\lambda)^p$ proving
\rf{IT1}.
\par Prove \rf{IT2}. By Proposition \reff{AP-DS},
\bel{ADD}
\AP(f;t:S)\le C\{\AP(f;4t:\partial
S)+t\|f|_\Omega\|_{\LOPO}\}, ~~~~t>0
\ee
so that, by\rf{ELC} and \rf{FYT}, $\AP(f;t:S)\le C\lambda
t,$ $t\in(0,\ve),$ proving \rf{IT2}.
\par The remaining inequality \rf{IT3} is trivial. In fact,
by Definition \reff{PR-AP},
$$ \PRM(f;t:S)\le\PRM(f;t:\partial S),~~~t>0, $$
so that, by \rf{ELC},
$$
\intl_0^\ve\PRM(f;t:S)^p\,
\frac{dt}{t^{p+1}}
\le  \intl_0^\ve\PRM(f;t:\partial S)^p\,
\frac{dt}{t^{p+1}}\le C\lambda^p.
$$
\par The proposition is completely proved.\bx
\begin{remark}\lbl{EST-AFT} {\em Theorem \reff{EXT2} and
estimates \rf{FTD} and \rf{ADD} imply the following
equivalence
\be \|f\|_{\WP|_S}&\sim& \|f|_\Omega\|_{\WPO}+ \|f\circ
T_{\partial S}\|_{L_p((\partial S)_\ve)}+ \sup_{0<t\le \ve}
\frac{\AP(f;t:\partial S)}{t}\nn\\&+&
\left\{\intl_0^\ve\PRM(f;t:S)^p\,
\frac{dt}{t^{p+1}}\right\}^{\frac{1}{p}}.
\nn
\ee
Here $f$, $p$, $\alpha$, $\ve$ and constants of equivalence
are the same as in Theorem \reff{EXT2}.}
\end{remark}
\par We finish the section with an important consequence of
Proposition \reff{EXT-AP}.
\par Recall that by $\omega(F;t)_{L_p}$ we denote the
modulus of continuity of a function $F$ in $L_p$, see
\rf{DOM}.
\begin{theorem}\lbl{EXT-MS} Let $\ve>0,\theta\ge 1,$ and
let $T:S_\ve\to S$ be a measurable mapping satisfying
inequality \rf{TR1}. There exists a linear continuous
extension operator $E:C(S)\to C(\RN)$ such that for every
$f\in C(S)$, $n<p<\infty$, and every $t\in(0,\ve_0)$ we
have:
\bel{OAP} \omega(Ef;t)_{L_p}\le
Ct\left\{\left(\intl_t^{\ve}\AP(f;u:S)^p\,
\frac{du}{u^{p+1}}\right)^{\frac1p}+\|f\circ
T\|_{L_p(S_\ve)}\right\}. \ee
\par In addition, if the mapping $T$ satisfies condition
\rf{TS}, then for every $t\in(0,\ve_0)$ and every
$\alpha\in(0,0.3/(1+\theta))$ the following inequality
$$
\omega(Ef;t)_{L_p}\le C\left\{\AP(f;\gamma t:S)+
t\left(\intl_t^{\ve}\PRM(f;u:S)\,
\frac{du}{u^{p+1}}\right)^{\frac1p}+t\|f\circ
T\|_{L_p(S_\ve)}\right\}
$$
holds. Here $\gamma$ is an absolute constant, and $\ve_0$
is a positive constant depending only on $\ve$ and
$\theta$.
\end{theorem}
\par {\it Proof.} We let $E:=\Es$ denote the same extension
operator as in the proof of Theorem \reff{EXT2}, i.e., we
put $\bc:=0$ and define $\dlt$ by formula \rf{DLT2}. \par
Clearly, by \rf{KPU},
$ \omp(Ef;t)_{L_p}\le \,\omp(Ef;t)_{L_\infty}, ~t>0,$
so that by \rf{OMAP} and \rf{EQM}
$ \omega(Ef;t)_{L_p}\le C \AP(Ef;t),~t>0.$
This inequality and Proposition \reff{EXT-AP} imply the
second inequality of the theorem.
\par Prove the first inequality. Observe that
$\PRM(f;\cdot:S)\le \AP(f;\cdot:S)$, see \rf{APCM}, so that
in the second inequality the quantity $\PRM(f;\cdot:S)$ can
be replaced with $\AP(f;\cdot:S)$. Moreover, a trivial
modification of the proof (see, in particular, Lemma
\reff{NORMPL}) allows us to omit condition \rf{TS} in this
inequality.
Moreover, we can assume that $2\gamma\ve_0\le \ve$ so that
$2\gamma t\le \ve$ for every $t\in (0,\ve_0)$. Since
$\AP(f;\cdot:S)$ is non-decreasing, we have
$$
t^p\intl_{\gamma t}^{\ve}\AP(f;u:S)^p\,
\frac{du}{u^{p+1}}\ge
(2\gamma)^{-p}(2^p-1)\AP(f;\gamma t:S)^p,
$$
proving the theorem.\bx
\SECT{6. Restrictions of Besov  $\BS$-spaces, $n/p<s<1$, to
closed sets.}{6}
\indent
\par In this section we present an intrinsic
characterization of the restrictions of Besov spaces $\BS$,
$s\in(n/p,1),$ to closed subsets of $\RN$ via local
oscillations.
\par We recall that the space $\BS, 0<s<1, 1\le
p\le\infty,$ $0<q\le\infty,$ is defined by the finiteness
of the norm
$$ \|f\|_{\BS}:= \|f\|_{L_p(\RN)}+\left(~\intl_{\|h\|<1}
\|\triangle^1_h f\|_{L_p(\RN)}^q\,
\frac{dh}{\|h\|^{sq+n}}\right)^{\frac{1}{q}} $$
(with standard modification if $q=\infty$), see, e.g.
\cite{BIN,T2,T3}.
\par Here $\triangle^1_h f(x)=f(x+h)-f(x).$ In particular,
for $p=q$,
$$ \|f\|_{B^s_{p,p}(\RN)}=
\|f\|_{L_p(\RN)}+\left(~\,\iint\limits_{\|x-y\|<1}
\frac{|f(x)-f(y)|^p} {\|x-y\|^{sp+n}}\, dx\,dy
\right)^{\frac{1}{p}}.$$
\par There are many different expressions for equivalent
norms in the Besov space. We will be interested in two of
them. The first equivalence is one of the classical
definition of the Besov space norm: Let  $p\ge 1,q>0,$
$s\in (0,1)$ and $T>0$. Then for every $f\in L_p(\RN)$
\bel{EN1} \|f\|_{\BS}\sim \|f\|_{L_p(\RN)}+ \left(\intl_0^T
\omega(f;t)_{L_p}^q\,
\frac{dt}{t^{sq+1}}\right)^{\frac{1}{q}}. \ee
\par The second equivalence involves the notion of local
oscillation $\Ec(f,Q(x,t))$, see \rf{OSC}. This result, see
\cite{Br1} and \cite{T2}, p. 186, states that for every
$f\in L_p(\RN)$
$$ \|f\|_{\BS}\sim \|f\|_{L_p(\RN)}+ \left(\intl_0^T
\|\Ec(f;Q(\cdot;t))\|_{L_p(\RN)}^q\,
\frac{dt}{t^{sq+1}}\right)^{\frac{1}{q}} $$
provided $p\ge 1,q>0,s\in(0,1)$ and $sp>n$.
\par This equivalence and \rf{AP-INT} imply the following
description of Besov spaces via $A_p$-functionals: Let
$p\ge 1,q>0,s\in(0,1)$ and $sp>n$. Then for every $f\in
L_p(\RN)$ we have
\bel{ENAP}
\|f\|_{\BS}\sim \|f\|_{L_p(\RN)}+
\left(\intl_0^T \AP(f;t)^q\,
\frac{dt}{t^{sp+1}}\right)^{\frac{1}{q}}.
\ee
In all cases the constants of equivalence depend only on
$n,s,p,q$ and $T.$
\par We turn to
\medskip
\par {\bf Proof of Theorem \reff{EXT-BS-INTRO}.}
\par We are needed two technical lemmas. First of them is a
particular case of the well known Hardy inequality. Let
$0<\varsigma,\eta\le 1,$ $0<p,q\le\infty,$ $a>0$. Given
$g\in L_{1,loc}(\R_+)$ put
$$ U_\varsigma g(t):=t^\varsigma\left(\intl_0^t
\left(\frac{|g(u)|} {u^\varsigma}\right)^p\,
\frac{du}{u}\right)^{\frac{1}{p}}, ~~~~~~~~~~ V_\varsigma
g(t):=t^\varsigma\left(\intl_t^a \left(\frac{|g(u)|}
{u^\varsigma}\right)^p\, \frac{du}{u}\right)^{\frac{1}{p}},
$$
and
$$
\|g\|_{\eta,q}:=\left(\intl_0^a \left(\frac{|g(t)|}
{t^\eta}\right)^q\, \frac{dt}{t}\right)^{\frac{1}{q}},
$$
(with the corresponding modification for $p$ or
$q=\infty$).
\begin{lemma}\lbl{HR-IN} For every $g\in L_{1,loc}(\R_+)$
the following two inequalities
\bel{HE1} \|U_\varsigma g\|_{\eta,q}\le
C\|g\|_{\eta,q},~~~{\rm if}~~\varsigma<\eta, \ee
and
\bel{HE2} \|V_\varsigma g\|_{\eta,q}\le
C\|g\|_{\eta,q},~~~{\rm if}~~\varsigma>\eta,  \ee
hold with constant $C=C(\varsigma,\eta,p,q)$.
\end{lemma}
\par Let us formulate the second lemma (its proof is
elementary and can be left to the reader as an easy
exercise).
\begin{lemma}\lbl{INTPQ} Let $g=g(t),t>0,$ be a
non-decreasing function. Then for every $0<p,q\le \infty$
and every $0<s,a<\infty$ the following inequality
$$
\left(\intl_{0}^{a}g(t)^p\,\frac{dt}{t}\right)^{\frac1p}
\le C \left(\intl_{0}^{2a}\left(\frac{g(u)}{u^s}\right)^q
\,\frac{du}{u}\right)^{\frac1q}
$$
holds (with the standard modification for $p=\infty$ or
$q=\infty$). Here $C$ is a constant depending only on
$p,q,s,$ and $a$.
\end{lemma}
\medskip
\par We turn to the proof of Theorem \reff{EXT-BS-INTRO}.
\smallskip \par  {\it (Necessity.)} We put $\tQ:=11\theta
Q$ and
$$ I:=\sum_{Q\in\TWE}|Q|\Ec(F;\tQ)^p. $$
By $\pi_m$ we denote the collections of cubes
$$
\pi_m:=\{Q\in\TWE:~2^m<r_Q\le 2^{m+1}\}.
$$
Since $r_Q=(1/2)\diam Q\le\ve/2$ for each $Q\in\TWE$, we
may suppose that $m \le M$ where $M$ satisfies the
following inequality:
$2^M\le \ve/2\le 2^{M+1}.$ (Otherwise $\pi=\emp.$)
\par Given $Q\in\pi_m$ we put
$Q^\curlywedge:=Q(x_Q,11\theta\,2^{m+1})$. Then
$Q\subset\tQ\subset Q^\curlywedge$. By Lemma \reff{FM}, the
collection
$\pi^\curlywedge_m:=\{Q^\curlywedge:~Q\in\pi_m\}$ of equal
cubes (of diameter $22\theta\,2^{m+1}$) can be partitioned
into at most $N=N(n,\theta)$ packings. Therefore without
loss of generality we may assume that $ \pi^\curlywedge$ is
a {\it packing}.
\par Put $\tau_m:=11\theta\,2^{m+1}$. Then, by \rf{DI},
$$ \sum_{Q\in\pi_m}|Q|\Ec(F;\tQ)^p\le
\sum_{Q\in\pi_m}|Q|\Ec(F;Q^\curlywedge)^p\le
C\sum_{K\in\pi^\curlywedge_m}|K|\Ec(F;K)^p\le
C\AP(F;\tau_{m+1})^p. $$
Hence,
$$ I=\sum_{m=-\infty}^M\sum_{Q\in\pi_m}|Q|\Ec(F;\tQ)^p \le
C\sum_{m=-\infty}^M \AP(F;\tau_{m+1})^p. $$
Since $\AP(\cdot;t)$ is a non-decreasing function, it can
be readily seen that
$$
I\le C\sum_{m=-\infty}^M
\intl_{\tau_{m+1}}^{\tau_{m+2}}\AP(F;t)^p\,\frac{dt}{t}
=C\intl_{0}^{\tau_{M+2}}\AP(F;t)^p\,\frac{dt}{t}.
$$
Since $2^{M+1}\le\ve$, we have $\tau_{M+2}=11\theta
2^{M+2}\le 22\theta\ve=:\tve.$ Hence, by Lemma
\reff{INTPQ}, we obtain
$$
I^{\frac1p}\le C
\left(\intl_{0}^{2\tve}\left(\frac{\AP(F;t)}{t^s}\right)^q
\,\frac{dt}{t}\right)^{\frac1q}.
$$
\par On the other hand, by Lemma  \reff{NR},
$ \|F\circ T\|_{L_p(\Se)}\le
C\left\{\|F\|_{L_p(\RN)}+I^{\frac1p}\right\}. $
Combining this inequality, equivalence \rf{ENAP} and
obvious inequality
$$
\AP(F|_S;t:S)\le \AP(F;t),~~~~F\in C(\RN),
$$
we obtain the statement of the necessity part of the
theorem.
\smallskip
\par {\it (Sufficiency.)} The proof of the
sufficiency part follows the same scheme as the proof of
Theorem \reff{EXT2}. Given a function $f\in C(S)$ we let
$\lambda$ denote the right-hand side of the equivalence
\rf{BSTRN}.
\par As in Theorem \reff{EXT2}, we define
an extension $F$ of $f$ by formulas \rf{DLT2} and
\rf{DEFEXT}. Then we prove that $F\in\BS$ and
$\|F\|_{\BS}\le C\lambda.$
\par First we show that $\|F\|_{L_p(\RN)}\le C\lambda.$
Note that $\PRM(f;\cdot:S)\le \AP(f;\cdot:S)$, so that by
Lemma \reff{NORMP}
\be \|F\|_{L_p(\RN)}&\le& C\left\{\|f\circ
T\|_{L_p(S_\ve)}+ \left(\intl_0^{\ve/2}\PRM(f;t:S)^p\,
\frac{dt}{t}\right)^{\frac1p}\right\}\nn\\
&\le& C\left\{\|f\circ T\|_{L_p(S_\ve)}+
\left(\intl_0^{\ve/2}\AP(f;t: S)^p\,
\frac{dt}{t}\right)^{\frac1p}\right\}.\nn \ee
This inequality and Lemma \reff{INTPQ} yield
$$ \|F\|_{L_p(\RN)}\le C\left\{\|f\circ T\|_{L_p(S_\ve)}+
\left(\intl_0^{\ve} \left(\frac{\AP(f;t:S)}{t^s}\right)^q\,
\frac{dt}{t}\right)^{\frac1q}\right\}\le C\lambda.$$
\par On the other hand, inequality \rf{OAP} and the Hardy
inequality \rf{HE2} imply the following estimate
$$
\left(\intl_0^\ve \left(\frac{\omega(F;t)_{L_p}}
{t^s}\right)^q\, \frac{dt}{t}\right)^{\frac{1}{q}}\le C
\|f\circ T\|_{L_p(S_\ve)}+ C\left(\intl_0^\ve\left
(\frac{\AP(f;t:S)}{t^s}\right)^q\,
\frac{dt}{t}\right)^{\frac{1}{q}}\le C\lambda.
$$
Finally, by \rf{EN1}, $\|F\|_{\BS}\le C\lambda$.
\par Theorem \reff{EXT-BS-INTRO} is proved.\bx
\medskip
\par In the aforementioned paper \cite{J} A. Jonsson
describes the trace space $\BS|_S$, $n/p<s<1,$ via certain
$L_p$-oscillations of functions with respect to a doubling
measure supported on $S$. Let us recall this result.
\par Let $\mu$ be a doubling measure with support $S$
satisfying the condition $(D_n)$ \rf{DN} and the condition
\rf{B1}. Observe that for every cube $Q=Q(x,r)$ such that
$x\in S$ and $r\le 1$, by the $(D_n)$-condition,
$$
\mu (Q(x,1))\le C_\mu (1/r)^n\mu(Q(x,r))=2^n\,C_\mu
\mu(Q(x,r))/|Q(x,r)|.
$$
By \rf{B1}, $C'_\mu\le\mu(Q(x,1))$ so that
\bel{LEBMU} |Q(x,r)|\le C\mu(Q(x,r)),~~~~{\rm
for~every}~~x\in S,~r\le 1, \ee
with $C=2^nC_\mu/C'_\mu$.
\begin{theorem}\lbl{EXT-BMU} (\cite{J}) Let $n/p<s<1$,
$1\le p,q\le\infty,$ and let $\mu$ be a measure with
support $S$ satisfying conditions $(D_n)$ and \rf{B1}. Then
for every $f\in C(S)$
\be &&\|f\|_{\Bspq(\RN)|_S}\sim \|f\|_{L_p(\mu)}\nn\\&+&
\left\{\sum_{\nu=0}^\infty\left(2^{\nu(s-n/p)}
\left(~\iint\limits _{\|x-y\|<2^{-\nu}}
\frac{|f(x)-f(y)|^p}
{\mu(Q(x,2^{-\nu}))\mu(Q(y,2^{-\nu}))}\,d\mu(x)\,d\mu(y)
\right)^{\frac1p}\right)^q\right\}^{\frac{1}{q}} \nn \ee
(with the standard modification of the right-hand side
whenever $p=\infty$ or $q=\infty$). The constants of this
equivalence depend only on $n,p,q,s,$ and parameters
$C_\mu, C'_\mu, C''_\mu$.
\end{theorem}
\par It is shown in \cite{J} that for $p=q$
$$ \|f\|_{\Bspq(\RN)|_S}\sim\|f\|_{L_p(\mu)}+
\left(~\iint\limits_{\|x-y\|<1} \frac{|f(x)-f(y)|^p}
{\|x-y\|^{sp-n}\mu(Q(x,\|x-y\|))^2}\,d\mu(x)\,d\mu(y)
\right)^{\frac1p}. $$
\par Let us show that the second item on the right-hand
side of these equivalences may be expressed in terms of
certain $L_p(\mu)$-version of the functional $A_p$\,, see
Definition \reff{DEF-APLS} below.
\par Let $\mu$ be a doubling measure on $\RN$ with support
$\supp\mu=S$ and let $1\le q\le\infty$. Given a function
$f\in L_{q,loc}(\mu)$ and a cube $Q\subset\RN$, we define
$L_q$-oscillation of $f$ on $Q$ (with respect to $\mu$) by
letting
\bel{EMU} \EFQ_{L_q(\mu)}:=\left(\frac{1}{\mu(Q)^2}
\iint\limits_{Q\times Q}|f(x)-f(y)|^q\,
d\mu(x)d\mu(y)\right)^{\frac1q}, \ee
whenever $1\le q<\infty$, and
$$
\EFQ_{L_\infty(\mu)}:=\Ec(f;Q\cap S)=\sup_{x,y\in Q\cap
S}|f(x)-f(y)|.
$$
It can be readily seen that
\bel{EBST}
\EFQ_{L_q(\mu)}\sim \Ec_1(f;Q)_{L_q(\mu)}
\ee
where
$$
\Ec_1(f;Q)_{L_q(\mu)}:=
\inf_{c\in\R}\mu(Q)^{-\frac1q}\|f-c\|_{L_q(Q:\mu)}.
$$
In addition,
\bel{MO} \EFQ_{L_q(\mu)}\sim
\left(\frac{1}{\mu(Q)}\intl_Q|f(x)-f_{Q,\mu}|^q\,
d\mu(x)\right)^{\frac1q}. \ee
Recall that  $f_{Q,\mu}$ denotes the $\mu$-average of $f$
on $Q$.
\par Since $Q\subset Q(x,2r_Q)$ and $Q(x,r_Q)\subset 2Q,$
whenever $x\in Q$, by the doubling condition \rf{DOUBC},
\bel{EM1} \EFQ_{L_q(\mu)}\sim \left(~\iint\limits_{Q\times
Q}\frac{|f(x)-f(y)|^q}{\mu(Q(x,r_Q))^2}\,
d\mu(x)d\mu(y)\right)^{\frac1q} \ee
provided $r_Q\in(0,1/2]$. Here the constants of equivalence
depend only on $q$ and the doubling constant $c_\mu$.
\begin{definition}\lbl{DEF-APLS}
{\em Given $1\le p<\infty$,
$1\le q\le\infty$ and a function $f\in L_{q,loc}(\mu)$ we
define the functional $\AP(f,\cdot:S)_{L_q(\mu)}$ by
letting
$$
\AP(f;t:S)_{L_q(\mu)}:= \sup_{\pi}\{\sum_{Q\in\pi}|Q|
\EFQ_{L_q(\mu)}^p\}^{\frac{1}{p}}
$$
where $\pi$ runs over all packings  of equal cubes with
diameter at most $t$ whose centers belong to $S$.
}\end{definition}
\par Clearly, if  $S=\supp\mu$, then
$$
\AP(f;t:S)_{L_\infty(\mu)}=\AP(f;t:S), ~~~~t>0,
$$
see Definition \reff{DEF-AP}, so that for all $1\le
q\le\infty$ we have
\bel{AP-APH} \AP(f;t:S)_{L_q(\mu)}\le\AP(f;t:S), ~~~~t>0.
\ee
\par We present two important properties of the functional
$\AP(\cdot;\cdot:S)_{L_p(\mu)}$\,.
\begin{lemma}\lbl{APLP} For every function $f\in L_p(\mu)$,
we have
$$
\AP(f;t:S)_{L_p(\mu)}\le
C(1+t^{\frac{n}{p}})\|f\|_{L_p(\mu)}, ~~~~t>0,
$$
where $C$ is a constant depending only on $p,n, C_\mu$ and
$C'_\mu$.
\end{lemma}
\par {\it Proof.} Let $\pi=\{Q\}$ be a packing of equal
cubes of diameter $\diam Q=s\le t$ centered in $S$. By
\rf{EBST} for each $Q\in\pi$ we have
$$
\EFQ_{L_p(\mu)}^p\le C \inf_{c\in\R}
\frac{1}{\mu(Q)}\intl_Q |f-c|^p\,d\mu\le C
\frac{1}{\mu(Q)}\intl_Q|f|^p\,d\mu.
$$
By \rf{LEBMU}, $|Q|\le C\mu(Q)$ provided $r_Q\le 1$, so
that whenever $s\le 1$ we obtain
$$
|Q|\EFQ_{L_p(\mu)}^p\le C
\frac{|Q|}{\mu(Q)}\intl_Q|f|^p\,d\mu\le
C\,\intl_Q|f|^p\,d\mu.
$$
\par Let $1<s=\diam Q\le t$. Then
$\mu(Q)\ge\mu(Q(x_Q,1))\ge C'_\mu$ so that
$$
|Q|\EFQ_{L_p(\mu)}^p\le C
\frac{|Q|}{\mu(Q)}\intl_Q|f|^p\,d\mu\le
C\,\frac{|Q|}{C'_\mu}\intl_Q|f|^p\,d\mu\le
C\,\frac{t^n}{C'_\mu}\intl_Q|f|^p\,d\mu.
$$
Hence, for every $Q\in\pi$ of diameter $\diam Q=s\in(0,t]$
we have
$$
|Q|\EFQ_{L_p(\mu)}^p\le C(1+t^n) \intl_Q|f|^p\,d\mu,
$$
so that
$$
\sum_{Q\in\pi}|Q|\EFQ_{L_p(\mu)}^p\le C(1+t^n)
\sum_{Q\in\pi}\intl_Q|f|^p\,d\mu\le C(1+t^n)
\intl_{\RN}|f|^p\,d\mu.
$$
It remains to take the supremum over all packings $\pi$ of
equal cubes with diameter at most $t$ centered in $S$, and
the lemma follows.\bx
\begin{proposition}\lbl{FAPN}  Let $\mu$ be a doubling
measure and let $f\in L_{p,loc}(\mu)$, $1\le p<\infty$.
Then for every $t\in(0,1/4]$ we have
$$ \frac1C\AP(f;t/4:S)^p_{L_p(\mu)}\le
\iint\limits_{\|x-y\|<t} \frac{t^n\,|f(x)-f(y)|^p}
{\mu(Q(x,t))^2}\,d\mu(x)\,d\mu(y)\le
C\AP(f;4t:S)^p_{L_p(\mu)}.
$$
Here $C$ is a positive constant depending only on $n,p$ and
$c_\mu$.
\end{proposition}
\par {\it Proof.} By \rf{EBST}, for every two cubes
$Q,\tQ$, $Q\subset\tQ$, we have
$$
\EFQ_{L_p(\mu)}\sim\mu(Q)^{-\frac1p}\inf_{c\in\R}
\|f-c\|_{L_p(Q:\mu)}\le\mu(Q)^{-\frac1p}\inf_{c\in\R}
\|f-c\|_{L_p(\tQ:\mu)},
$$
so that
$$ \EFQ_{L_p(\mu)}\le C
\frac{\mu(\tQ)^{\frac1p}}{\mu(Q)^{\frac1p}}
\,\Ec(f;\tQ)_{L_p(\mu)}. $$
\par First we prove the second inequality of the
proposition. We have
$$
I:=\iint\limits_{\|x-y\|<t} \frac{t^n\,|f(x)-f(y)|^p}
{\mu(Q(x,t))^2}\,d\mu(y)\,d\mu(x)= t^n\, \intl_{S}
\intl_{Q(x,t)}\frac{|f(x)-f(y)|^p}{\mu(Q(x,t))^2}
\,d\mu(y)\,d\mu(x).
$$
\par Let
\bel{TPI-P} \tpi:=\{Q_i=Q(x_i,t/2):~i=1,2,...\} \ee
be a covering of $S$ by cubes centered in $S$ with the
covering multiplicity bounded by a constant $C=C(n)$. Then,
$$
I \le Ct^n\sum_{i=1}^{\infty}\intl_{Q_i}
\intl_{Q(x,t)}\frac{|f(x)-f(y)|^p}{\mu(Q(x,t))^2}
\,d\mu(y)\,d\mu(x).
$$
Clearly, for each $x\in Q_i$ we have $Q(x,t)\subset
3Q_i=Q(x_i,3t)\subset 3Q(x,t)$ so that, by the doubling
condition, $ \mu(Q(x,t))\sim \mu(3Q_i)\sim \mu(3Q(x,t))$
provided $t\le 1/3.$ Hence, by \rf{EMU},
\be I&\le& Ct^n\sum_{i=1}^{\infty}\intl_{3Q_i}
\intl_{3Q_i}\frac{|f(x)-f(y)|^p}{\mu(Q(x,t))^2}
\,d\mu(y)\,d\mu(x)\nn\\
&\le& Ct^n\sum_{i=1}^{\infty}\frac{1}{\mu(3Q_i)^2}
\intl_{3Q_i} \intl_{3Q_i}|f(x)-f(y)|^p
\,d\mu(y)\,d\mu(x)\nn\\
&\le& C\sum_{i=1}^{\infty} |3Q_i|\Ec(f;3Q_i)^p_{L_p(\mu)}.
\nn\ee
\par By Lemma \reff{FM}, the collection of cubes
$\{2Q_i:~i=1,2,...\}$ can be partitioned into at most
$N(n)$ packings so that, by Definition \reff{DEF-APLS},
$$
I\le C N(n)\AP(f;4t:S)^p_{L_p(\mu)}
$$
proving the second inequality of the proposition.
\par Prove the first inequality. Let $0<\tau\le t$ and let
$\pi$ be a packing of equal cubes of diameter $\tau$
centered in $S$. Put
$$
J:=\sum_{Q\in\pi} |Q|\Ec(f;Q)^p_{L_p(\mu)}.
$$
Then by \rf{EBST} and \rf{DN}
$$
J\le C\sum_{Q\in\pi} \frac{\tau^n}{\mu(Q)}\inf_{c\in\R}
\intl_Q|f-c|^p\,d\mu\le C\sum_{Q\in\pi}
\frac{t^n}{\mu(Q(x_Q,t))}\inf_{c\in\R}
\intl_Q|f-c|^p\,d\mu.
$$
\par Given a cube $Q_i\in\tpi$, see \rf{TPI-P}, we put
$\pi_i:=\{Q\in\pi:~x_Q\in Q_i\}.$ Since $\tau\le t$, each
cube $Q\in\pi_i$ lies in the cube $3Q_i$ so that
$$
\sum_{Q\in\pi_i} \inf_{c\in\R} \intl_Q|f-c|^p\,d\mu \le
\inf_{c\in\R}\sum_{Q\in\pi_i} \intl_Q|f-c|^p\,d\mu\le
\inf_{c\in\R} \intl_{3Q_i}|f-c|^p\,d\mu.
$$
On the other hand, since $x_Q\in Q_i$, the cube
$3Q_i\subset Q(x_Q,2t)$, so that, by the doubling
condition, $\mu(3Q_i)\sim \mu(Q(x_Q,t))$ provided $t\le
1/2.$ Hence,
\be J&\le& Ct^n\sum_{i=1}^\infty\sum_{Q\in\pi_i}
\frac{1}{\mu(Q(x_Q,t))}\inf_{c\in\R}
\intl_Q|f-c|^p\,d\mu\nn\\
&\le&
Ct^n\sum_{i=1}^\infty\frac{1}{\mu(3Q_i)}\sum_{Q\in\pi_i}
\inf_{c\in\R} \intl_Q|f-c|^p\,d\mu, \nn\\
&\le& Ct^n\sum_{i=1}^\infty
\frac{1}{\mu(3Q_i)}\inf_{c\in\R} \intl_{3Q_i}|f-c|^p\,d\mu,
\nn\ee
so that, by \rf{EBST},
$$
J\le Ct^n\sum_{i=1}^\infty\Ec(f;3Q_i)_{L_p(\mu)}.
$$
Hence, by \rf{EM1},
$$
J\le Ct^n\sum_{i=1}^\infty\intl_{3Q_i}\intl_{3Q_i}
\frac{|f(x)-f(y)|^p}{\mu(Q(x,3t/2))^2}\,
d\mu(y)d\mu(x).
$$
Note that $3Q_i\subset Q(x,3t)$ for each $x\in 3Q_i$. In
addition, by the doubling condition, $\mu(Q(x,4t))\sim
\mu(Q(x,3t/2))$, $x\in 3Q_i$, provided $t\le 1/4$. Hence,
$$
J\le Ct^n\sum_{i=1}^\infty\intl_{3Q_i}\intl_{Q(x,3t)}
\frac{|f(x)-f(y)|^p}{\mu(Q(x,4t))^2}\,
d\mu(y)d\mu(x).
$$
\par By Lemma \reff{FM}, without loss of generality one can
assume that the collection of cubes $\{3Q_i:~i=1,2,...\}$
is a packing. Hence
\be J&\le& Ct^n\intl_{S}\intl_{Q(x,3t)}
\frac{|f(x)-f(y)|^p}{\mu(Q(x,4t))^2}\,
d\mu(y)d\mu(x)\nn\\&\le& Ct^n\iint\limits_{\|x-y\|<4t}
\frac{|f(x)-f(y)|^p}{\mu(Q(x,4t))^2}\, d\mu(y)d\mu(x).
\nn\ee
Taking the supremum over all packings $\pi$ of equal cubes
centered in $S$ with diameter $\tau\le t$, we finally
obtain
$$
\AP(f;t:S)_{L_p(\mu)}^p \le Ct^n\iint\limits_{\|x-y\|<4t}
\frac{|f(x)-f(y)|^p}{\mu(Q(x,4t))^2}\,
d\mu(y)d\mu(x).
$$
\par The proposition is proved.\bx
\par Now, combining this proposition with Theorem
\reff{EXT-BMU}, we obtain
$$
\|f\|_{B^s_{p,q}(\RN)|_S}\sim \|f\|_{L_p(\mu)}+
\left\{\sum_{\nu=0}^\infty 2^{\nu sq}
\AP(f;2^{-\nu}:S)_{L_p(\mu)}^q\right\}^{\frac{1}{q}}.
$$
Finally, this equivalence and Lemma \reff{APLP} imply the
following formula for the trace norm in the Besov space:
\par {\it Let $\mu$ be a measure with support $S$
satisfying conditions $(D_n)$ and \rf{B1}. Let  $1\le
p,q\le\infty,\ve>0,$ and let $n/p<s<1$. Then for every
$f\in C(S)$ the following equivalence
$$
\|f\|_{B^s_{p,q}(\RN)|_S}\sim \|f\|_{L_p(\mu)}+
\left\{\intl_{0}^{\ve}
\left(\frac{\AP(f;t:S)_{L_p(\mu)}}{t^s}\right)^q\,
\frac{dt}{t}\right\} ^{\frac{1}{q}}
$$
holds. The constants of this equivalence depend on $\ve$
and the same parameters as the constants in Theorem
\reff{EXT-BMU}.}
\SECT{7. Restrictions of Sobolev spaces and doubling
measures.}{7}
\indent
\par In this section we prove the results stated in
subsection 1.2.
\par Let $\mu$ be a positive Borel regular measure
supported on the closed set $S$. We assume that $\mu$
satisfies condition $(D_n)$ \rf{DN} and condition \rf{B1}.
\begin{definition}\lbl{AP-PRMU} {\em Let $\alpha\in(0,1]$,
$p\in[1,\infty)$ and $q\in[1,\infty]$. Given $f\in
L_{q,loc}(\mu)$ we define the functional
$\PRM(f,\cdot:S)_{L_q(\mu)}$ by letting
\bel{PORMU}
\PRM(f;t:S)_{L_q(\mu)}:=
\sup_{\pi}\{\sum_{Q\in\pi}|Q|
\EFQ_{L_q(\mu)}^p\}^{\frac{1}{p}},~~~~t>0,
\ee
where $\pi=\{Q\}$ runs over all $\alpha$-porous (with
respect to $S$) packings of equal cubes with centers in
$\partial S$ and diameter at most $t$.
}\end{definition}
\par Clearly, by Definition \reff{PR-APS}, for every $1\le
q\le\infty$ and $\alpha\in(0,1]$
\bel{AMU-AP} \PRM(f;t:S)_{L_q(\mu)}\le\SPR(f,t:S),~~~~t>0.
\ee
\par It is also clear that
\bel{Y} \PRM(f;t:S)_{L_q(\mu)}\le
\AP(f;t:S)_{L_q(\mu)},~~~~t>0,\ee
see Definition \reff{DEF-APLS}. Also observe that the
functional $\PRM(f;\cdot:S)$, see Definition \reff{PR-AP},
equals the functional $\PRM(f;\cdot:S)_{L_\infty(\sigma)}$
where $\sigma$ is a measure whose support coincides with
$\partial S$.
\par The first main result of the section is the following
\begin{theorem}\lbl{EXT-MU-NORM} Let $p\in(n,\infty),$
$\alpha\in(0,1/7)$, and $\ve>0$. Then for every function
$f\in C(S)$ we have
\be
\|f\|_{\WP|_S}&\sim& \|f\|_{L_p(\mu)}+
\sup_{0<t\le \ve} \frac{\AP(f;t:S)_{L_p(\mu)}}{t}\nn\\
&+& \left\{\intl_0^\ve\PRM(f;t:S)_{L_p(\mu)}^p\,
\frac{dt}{t^{p+1}}\right\}^{\frac{1}{p}} \nn\ee
with constants of equivalence depending only on
$n,p,\alpha,\ve$ and the parameters $C_\mu,C'_\mu,C''_\mu$.
\end{theorem}
\par The proof is based on a series of auxiliary lemmas and
propositions. In all of these results we will assume that
$f\in C(S)$, $p\in(n,\infty), \ve>0,$ and $C$ is a constant
depending only on $n,p,\alpha,\ve$ and
$C_\mu,C'_\mu,C''_\mu$.
\begin{lemma}\lbl{MU-LP} We have $\|f\|_{L_p(\mu)}\le
C\,\|f\|_{\WP|_S}.$
\end{lemma}
\par {\it Proof.} Let $F\in \WP$ be a continuous function
such that $F|_S=f$ and
$$ \|F\|_{\WP}\le 2\,\|f\|_{\WP|_S}. $$
\par Let $\pi=\{Q_i=Q(x_i,1/2):~i=1,2,...\}$ be a covering
of $S$ by equal cubes of diameter $1$ centered in $S$ whose
covering multiplicity $M_\pi\le N(n).$ By condition
\rf{B1}, $\mu(Q_i)\le C''_\mu$ so that
$$
\|f\|_{L_p(\mu)}^p=\|F\|_{L_p(\mu)}^p\le \sum_{i=1}^\infty
\intl_{Q_i}|F|^p\,d\mu\le\sum_{i=1}^\infty
(\sup_{Q_i}|F|^p)\mu(Q_i)\le C''_\mu\sum_{i=1}^\infty
\sup_{Q_i}|F|^p.
$$
Hence
$$
\|f\|_{L_p(\mu)}^p\le C\{\sum_{i=1}^\infty
\sup_{Q_i}|F-F_{Q_i}|^p+\sum_{i=1}^\infty|F_{Q_i}|^p\} \le
C\{\sum_{i=1}^\infty
\Ec(F;Q_i)^p+
\sum_{i=1}^\infty|Q_i|^{-1}\intl_{Q_i}|F|^p\,dx\}.
$$
(Recall that $F_{Q}$ stands for the average of $F$ on $Q$.)
Since $|Q_i|=1$, we obtain
$$
\|f\|_{L_p(\mu)}^p\le C\{\sum_{i=1}^\infty
|Q_i|\Ec(F;Q_i)^p+\sum_{i=1}^\infty\intl_{Q_i}|F|^p\,dx\}.
$$
\par Since $M_\pi\le N(n),$ the second term on the
right-hand side of this inequality is bounded by
$N(n)\|F\|_{L_p(\RN)}^p$. Let us estimate the first term.
By Theorem \reff{TFM}, the collection $\pi$ can be
partitioned into at most $C(n)$ packings so that without
loss of generality we may assume that $\pi$ itself is a
packing. Then, by \rf{DI} and \rf{EM},
$$
\sum_{i=1}^\infty |Q_i|\Ec(F;Q_i)^p\le \AP(F;1)^p\le
C\|F\|_{\LOP}^p.
$$
Hence,
$$
\|f\|_{L_p(\mu)}^p\le
C(\|F\|_{\LOP}^p+\|F\|^p_{L_p(\RN)})\le C\|F\|_{\WP}^p\le
C\|f\|_{\WP|_S}^p,
$$
proving the lemma.\bx
\par We put
$$ \|f\|'_S= \|f\|_{L_p(\mu)}+ \sup_{0<t\le \ve}
\frac{\AP(f;t:S)_{L_p(\mu)}}{t} +
\left\{\intl_0^\ve\PRM(f;t:S)_{L_p(\mu)}^p\,
\frac{dt}{t^{p+1}}\right\}^{\frac{1}{p}}. $$
\par Combining Lemma \reff{MU-LP} with inequalities
\rf{AP-APH}, \reff{N0}, \rf{AMU-AP}, and Proposition
\reff{NL}, we obtain the required inequality
$ \|f\|'_S\le C\|f\|_{\WP|_S} $
provided $n<p<\infty, \ve>0$, and $0<\alpha\le 1$.
\par Let us prove that for every $p>n,$ $\ve>0$,
and $0<\alpha<1/7$ the opposite inequality
\bel{OPIN}
\|f\|_{\WP|_S}\le C\|f\|'_S.
\ee
holds as well. We begin the proof with the next auxiliary
\begin{lemma}\lbl{MUS} For every $\ve>0$ we have
$$
\|f\|_{L_p(S)}\le C\{\|f\|_{L_p(\mu)}+\AP(f;\ve:S)\}.
$$
\end{lemma}
\par {\it Proof.} It suffices to prove the lemma for
$\ve\in(0,1]$. We let $\pi:=\{Q_i=Q(x_i,\ve/2):~x_i\in
S,\,i\in I\}$ denote a covering of $S$ with the covering
multiplicity $M_\pi\le C(n)$. We have
$$
\|f\|_{L_p(S)}^p\le \sum_{i\in I}\intl_{Q_i\cap
S}|f|^p\,dx\le\sum_{i\in I}|Q_i|\sup_{Q_i\cap S}|f|^p.
$$
Hence,
\be \|f\|_{L_p(S)}^p&\le& 2^p\{\sum_{i\in
I}|Q_i|\sup_{Q_i\cap S} |f-f_{Q_i,\mu}|^p+\sum_{i\in I}
|Q_i||f_{Q_i,\mu}|^p\}\nn\\ &\le& 2^p\{\sum_{i\in
I}|Q_i|\Ec(f;Q_i\cap S)^p+\sum_{i\in I}|Q_i|
|f_{Q_i,\mu}|^p\}.\nn \ee
By \rf{LEBMU}), $|Q_i|\le C\mu(Q_i)$ so that
$$
\|f\|_{L_p(S)}^p\le C\{\sum_{i\in I}|Q_i|\Ec(f;Q_i\cap
S)^p+\sum_{i\in I}\intl_{Q_i} |f|^p\,d\mu\}.
$$
\par Again, by Theorem \reff{TFM}, we may suppose that
$\pi$ is a packing, so that, by \rf{DI},
$$
\|f\|_{L_p(S)}^p\le C\{\AP(f;\ve:S)^p+\intl_{\RN}
|f|^p\,d\mu.\}
$$
proving the lemma.\bx
\par Recall that given $\ve>0$ by $\TWE$ we denote the
family of all Whitney's cubes of diameter at most $\ve$,
see \rf{DQEP}. Given $Q\in\TW$ we put
$$\bQ:=Q(a_Q,r_Q).$$
\begin{lemma}\lbl{T-AV} For every $0<\ve\le 1/50$ we have
$$
\sum_{Q\in\TWE}|Q||f_{\bQ,\mu}|^p\le C\{\|f\|_{L_p(\mu)}^p
+\|f\|_{L_p(S)}^p\}.
$$
\end{lemma}
\par {\it Proof.} The reader can easily see that, by
Theorem \reff{Wcov}, for $Q\in\TWE$ and every $y\in
\bQ=Q(a_Q,r_Q)$ we have the following:
$$ \bQ\subset Q(x_Q,10r_Q)=10 Q, $$
\bel{QAQ} \bQ\subset 10 Q\subset Q(x,11Q),~~~x\in Q, \ee
and
$$ Q(x,11r_Q)\subset Q(a_Q,22r_Q)=22\bQ,~~~x\in Q. $$
This and \rf{QAQ} yield
\bel{MAQ} \mu(\bQ)\le \mu(Q(x,11Q))\le \mu(22\tQ), ~~~x\in
Q. \ee
\par Since $Q\in\TWE$, its diameter $\diam Q=2r_Q\le \ve$
so that
$$
\diam 22\bQ=44r_{\bQ}=44\,r_Q\le 44\ve\le 1.
$$
Therefore, by the doubling condition and \rf{MAQ},
\bel{EQ1}
\mu(Q(x,11Q))\sim \mu(\bQ), ~~~x\in Q.
\ee
In addition, by \rf{LEBMU},
\bel{M1}
|Q|=|\bQ|\le C\mu(\bQ).
\ee
\par We put $\tmu:=\mu+\lambda$ where $\lambda$ denotes the
Lebesgue measure on $\RN$. By $\fc$ we denote the extension
of $f$ by $0$ from $S$ to all of $\RN$. Thus $\fc(x)=f(x),
x\in S,$ and $\fc(x)=0$ on $\RN\setminus S$. Given $g\in
L_{1,loc}(\tmu)$ we let $\Mc_{\tmu}g$ denote the
Hardy--Littlewood maximal function of $g$ (with respect to
the measure $\tmu$):
$$
\Mc_{\tmu}g(x):=\sup_{r>0}\frac{1}{\tmu(Q(x,r))}
\intl_{Q(x,r)} |g|\,d\tmu.
$$
Now, by \rf{QAQ}, for every $x\in Q$ we have
$$
|f_{\bQ,\mu}|\le \frac{1}{\mu(\bQ)} \intl_{\bQ} |f|\,d\mu
\le \frac{1}{\mu(\bQ)} \intl_{Q(x,11r_Q)} |f|\,d\mu\le
\frac{1}{\mu(\bQ)} \intl_{Q(x,11r_Q)} |\fc|\,d\tmu,
$$
so that
$$
|f_{\bQ,\mu}|\le \frac{\tmu(Q(x,11r_Q))}{\mu(\bQ)}
\frac{1}{\tmu(Q(x,11r_Q))}\intl_{Q(x,11r_Q)} |\fc|\,d\tmu
\le \frac{\tmu(Q(x,11r_Q))}{\mu(\bQ)} [\Mc_{\tmu}(\fc)](x).
$$
But, by \rf{EQ1} and \rf{M1},
$$
\frac{\tmu(Q(x,11r_Q))}{\mu(\bQ)}=
\frac{\mu(Q(x,11r_Q))+|Q(x,11r_Q)|}{\mu(\bQ)}\le C
\frac{\mu(\bQ)+|Q|}{\mu(\bQ)}\le C,
$$
so that
$$
|f_{\bQ,\mu}|\le C\, [\Mc_{\tmu}(\fc)](x), ~~~~x\in Q.
$$
Since $|Q|\le |Q|+\mu(Q)=\tmu(Q),$ we obtain
$$
|Q|\,|f_{\bQ,\mu}|^p\le \tmu(Q)\,|f_{\bQ,\mu}|^p\le C\,
\intl_{Q}[\Mc_{\tmu}(\fc)]^p(x)\,d\tmu(x).
$$
Hence,
\be I&:=&\sum_{Q\in\TWE}|Q||f_{\bQ,\mu}|^p\le
C\sum_{Q\in\TWE}\intl_{Q}
[\Mc_{\tmu}(\fc)]^p(x)\,d\tmu(x)\nn\\
&\le&
C\,M_{\TW}\intl_{\RN}[\Mc_{\tmu}(\fc)]^p(x)\,d\tmu(x).
\nn\ee
(Recall that $M_{\TW}$ denotes the covering multiplicity of
the family $\TW$ of Whitney's cubes.) By Theorem
\reff{Wcov}, $M_{\TW}\le N(n)$ so that
$$
I\le C\intl_{\RN}[\Mc_{\tmu}(\fc)]^p(x)\,d\tmu(x)
=C\|\Mc_{\tmu}(\fc)\|^p_{L_p(\tmu)}.
$$
Observe that  $\tmu=\mu+\lambda$ is a positive Borel
regular measure on $\RN$. In this case the
Hardy--Littlewood maximal operator is a bounded sublinear
operator from $L_p(\tmu)$ into $L_p(\tmu), p>1,$ whose
operator norm does not exceed a constant depending only on
$p$ and $n$. See, e.g., \cite{GK}, Corollary 3.3. Thus,
$$
\|\Mc_{\tmu}(\fc)\|^p_{L_p(\tmu)}\le
C\|\fc\|^p_{L_p(\tmu)},
$$
so that
\be I&\le&
C\|\fc\|^p_{L_p(\tmu)}=C\intl_{\RN}|\fc|^p\,d\tmu=
C\{\intl_{\RN}|f|^p\,d\mu+\intl_{\RN}|\fc|^p\,dx\}\nn\\&=&
C\{\|f\|_{L_p(\mu)}^p +\|f\|_{L_p(S)}^p\},\nn\ee
proving the lemma.\bx
\begin{definition}\lbl{STR-APS} {\em Let $0\le\alpha\le 1$
and let $Q$ be a cube centered in $S$. We say that $Q$ is
{\it strongly $\alpha$-porous} if the cube $\eta Q$ is
$\alpha$-porous (with respect to $S$) for every $0<\eta\le
1$. \par A packing $\pi$ is said to be {\it strongly
$\alpha$-porous} if every cube $Q\in\pi$ is strongly
$\alpha$-porous.
}\end{definition}
\par Given a function $f\in C(S)$ and a cube $Q\subset\RN$,
we put
\bel{ECT}
\Ect(f;Q):=\frac{1}{\mu(Q)}\intl_Q|f-f(x_Q)|\,d\mu.
\ee
Let us introduce a corresponding analog of the functional
\rf{PORMU} over families of strongly porous cubes.
\begin{definition}\lbl{STR-F} {\em \par Let
$p\in[1,\infty),\alpha\in[0,1]$, and let $f\in C(S)$. We
define the functional $\TRS(f,\cdot:S,\mu)$ by the formula
$$
\TRS(f;t:S,\mu):= \sup_{\pi}\{\sum_{Q\in\pi}|Q|\,
\Ect(f;Q)^p\}^{\frac{1}{p}}, ~~~~t>0.
$$
Here $\pi=\{Q\}$ runs over all {\it strongly}
$\alpha$-porous (with respect to $S$) packings of equal
cubes of diameter at most $t$ centered in $\partial S$.
}\end{definition}
\begin{proposition}\lbl{EXT-J} For every $\alpha\in(0,1/7)$
we have
$$
\|f\|_{\WP|_S}\le C\,\left(\|f\|_{L_p(\mu)}+\sup_{0<t\le
\ve} \frac{\AP(f;t:S)}{t} + \left\{\intl_0^{\ve}
\TRS(f;t:S,\mu)^p\,
\frac{dt}{t^{p+1}}\right\}^{\frac{1}{p}}\right).
$$
\end{proposition}
\par {\it Proof.} Since the righthand side of the
proposition's inequality is a non-decreasing function of
$\ve$, it suffices to prove the proposition for
\bel{MU-EPS} \ve\in(0,0.02). \ee
Put
\bel{MU-DLT} \bc:=0~~~~~{\rm and} ~~~~~~
\dlt:=\frac{\ve}{2000}, \ee
and define the function
$$ F:=\Es f $$
by the formula \rf{ExtOp}.
Since $ \|f\|_{\WP|_S} \le \|F\|_{\WP}, $ our aim is to
show that
$$\|F\|_{\WP}:=\|F\|_{L_p(\RN)}+\|F\|_{\LOP}$$
is bounded by the righthand side of the proposition's
inequality.
\par To estimate the norm $\|F\|_{L_p(\RN)}$ we will be
needed the following auxiliary lemmas.
\begin{lemma}\lbl{MU-APC} Let $0<\beta<1,$
$\varsigma\ge 1$ and $0<T\le 1/\varsigma.$ Let
$\pi=\{Q\}\subset\TW$ be a collection of Whitney's cubes of
diameter at most $T$. Given $Q\in\pi$ let $\tQ$ be a
strongly $\beta$-porous cube centered in $\partial S$ such
that
$$\|x_Q-x_{\tQ}\|\le\varsigma\diam Q ~~~~{\rm and}~~~~
\diam Q\le\diam \tQ\le \varsigma\diam Q.
$$
\par Then for every $f\in C(S)$ and every
$\alpha\in(0,\beta)$ we have
$$
\sum_{Q\in\pi}r_{\tQ}^{n-p} \Ect(f;\tQ)^p \le
C\intl_{0}^{4\varsigma T}
\TRS(f;t:S,\mu)^p \,\frac{dt}{t^{p+1}}
$$
where $C$ is a constant depending on $n,p, \varsigma$ and
$c_\mu$.
\end{lemma}
\par {\it Proof.} The lemma can be proven precisely in the
same way as Lemma \reff{APCQ} (with $t=0$ in its
formulation). We need only to replace in this proof
\rf{ET-Q} with the inequality
\bel{NEW-ETQ} \Ect(f;\tQ\cap S)\le C\,\Ect(f;Q^\sharp\cap
S). \ee
To prove this let us recall that, by \rf{QAV}, $\tQ\subset
Q^\sharp\subset 2\tQ$ and $x_{\tQ}=x_{Q^\sharp}$ so that
$$
\Ect(f;\tQ\cap S)
=\frac{1}{\mu(\tQ)}\intl_{\tQ}|f-f(x_{\tQ})|\,d\mu \le
\frac{1}{\mu(\tQ)}
\intl_{Q^\sharp}|f-f(x_{Q^\sharp})|\,d\mu.
$$
Since
$\diam \tQ\le\varsigma\diam Q\le\varsigma T\le 1,$
by the doubling condition \rf{DOUBC}
$\mu(Q^\sharp) \le \mu(2\tQ) \le c_\mu \mu(\tQ) , $
proving \rf{NEW-ETQ} with $C=c_\mu$.\bx
\begin{lemma}\lbl{FQ-EST} For every $\alpha\in(0,3/20)$
$$
\sum_{Q\in\TWB{90\dlt}}|Q||f(a_Q)|^p\le
\,C\left\{\|f\|_{L_p(\mu)}+\AP(f;\ve:S)+
\left(\intl_{0}^{\ve}
\TRS(f;t:S,\mu)^p
\,\frac{dt}{t^{p+1}}\right)^{\frac1p}\right\}.
$$
\end{lemma}
\par {\it Proof.} Put
$$
I:=\sum_{Q\in\TWB{90\dlt}}|Q||f(a_Q)|^p.
$$
Then
\bel{FEL} I\le
2^p\{\sum_{Q\in\TWB{90\dlt}}|Q||f(a_Q)-f_{\bQ,\mu}|^p+
\sum_{Q\in\TWB{90\dlt}}|Q||f_{\bQ,\mu}|^p\}=2^p\{I_1+I_2\}.
\ee
Recall that $\bQ:=Q(a_Q,r_Q)$. By \rf{MU-EPS} and
\rf{MU-DLT}, we have $90\dlt\le\ve<1/50$ so that by Lemma
\reff{T-AV}
$$
I_2:=\sum_{Q\in\TWB{90\dlt}}|Q||f_{\bQ,\mu}|^p\le\,
C\{\|f\|_{L_p(\mu)}+\|f\|_{L_p(S)}\}.
$$
Applying again Lemma \reff{MUS}, we obtain
\bel{EI2} I_2\le\, C\{\|f\|_{L_p(\mu)}+\AP(f;\ve:S)\}. \ee
\par Let us estimate the quantity
$$
I_1:= \sum_{Q\in\TWB{90\dlt}}|Q||f(a_Q)-f_{\bQ,\mu}|^p.
$$
We have
\be I_1&\le&
\sum_{Q\in\TWB{90\dlt}}|Q|\left(\frac{1}{\mu(\bQ)}\intl_{\bQ}
|f-f(a_Q)|\,d\mu\right)^p\nn\\&=&
\sum_{Q\in\TWB{90\dlt}}|\bQ|\left(\frac{1}{\mu(\bQ)}\intl_{\bQ}
|f-f(x_{\bQ})|\,d\mu\right)^p \le
C\sum_{Q\in\TWB{90\dlt}}r_{\bQ}^{n-p}\Ect(f;\bQ)^p.\nn\ee
Observe that, by Definition \reff{STR-APS} and by Lemma
\reff{QP}, the cube $\tQ:=Q(a_Q,20r_Q)$ is {\it strongly}
$\beta$-porous for every $\beta\in(0,3/20)$. Since
$\bQ=Q(a_Q,r_Q)=(1/20)\,\tQ,$
the cube $\bQ$ is strongly $\beta$-porous for every
$\beta\in(0,3/20)$ as well.
\par It can be easily seen that the cubes
$\{\bQ:~Q\in\TWB{90\dlt}\}$ satisfy conditions of Lemma
\reff{MU-APC} with $\varsigma=5$. Since $90\dlt\le
1/\varsigma =1/5$, see \rf{MU-EPS} and \rf{MU-DLT}, in the
lemma's hypothesis one can put  $T:=90\dlt$. Observe also
that $4\varsigma T=72\dlt\le\ve$ so that, by Lemma
\reff{MU-APC}, for every $\alpha\in(0,3/20)$ we have
$$
I_1\le C\,\intl_{0}^{4\varsigma T}
\TRS(f;t:S,\mu)^p \,\frac{dt}{t^{p+1}}\le
C\,\intl_{0}^{\ve} \TRS(f;t:S,\mu)^p
\,\frac{dt}{t^{p+1}}.
$$
\par Combining \rf{FEL} with this estimate and inequality
\rf{EI2}, we finally obtain the statement of the lemma.\bx
\begin{lemma}\lbl{FLP-MU} For every $\alpha\in(0,3/20)$ we
have
$$
\|F\|_{L_p(\RN)}\le
\,C\left\{\|f\|_{L_p(\mu)}+\AP(f;\ve:S)+
\left(\intl_{0}^{\ve}
\TRS(f;t:S,\mu)^p
\,\frac{dt}{t^{p+1}}\right)^{\frac1p}\right\}.
$$
\end{lemma}
\par {\it Proof.} By Lemma \reff{MUS}
\bel{AAA} \|F\|_{L_p(S)}=\|f\|_{L_p(S)}\le \,C
\{\|f\|_{L_p(\mu)}+\AP(f;\ve:S)\} \ee
provided $\ve\in(0,1]$. In turn, by Lemma \reff{LP-N},
$$
\|F\|_{L_p(\RN\setminus S)}\le
\,C\,\sum_{Q\in\TWB{2\dlt}}|Q||f(a_Q)|^p,
$$
so that, by Lemma \reff{FQ-EST},
$$
\|F\|_{L_p(\RN\setminus S)}\le
\,C\left\{\|f\|_{L_p(\mu)}+\AP(f;\ve:S)+
\left(\intl_{0}^{\ve}
\TRS(f;t:S,\mu)^p
\,\frac{dt}{t^{p+1}}\right)^{\frac1p}\right\}.
$$
\par This inequality and \rf{AAA} imply the result of the
lemma.\bx
\par We turn to estimates of the norm $\|F\|_{\LOP}$. We
will mainly follow the scheme used in the proofs of the
sufficiency part of  Theorems \reff{EXT-NORM} and
\reff{EXT2}.
\par We put
$$
\lambda(f):=\|f\|_{L_p(\mu)}+ \sup_{0<t\le \ve}
\frac{\AP(f;t:S)}{t} + \left\{\intl_0^{\ve}
\TRS(f;t:S,\mu)^p\,
\frac{dt}{t^{p+1}}\right\}^{\frac{1}{p}}.
$$
Our goal is to show that
\bel{PRST} \|F\|_{\LOP|_S}\le C\lambda(f), \ee
provided $p>n,$ $0<\ve\le 1/50$, and $0<\alpha<1/7$.
\par To prove this inequality given a packing $\pi=\{K\}$
of equal cubes in $\RN$ we again put
$$
I(F;\pi):=\left\{\sum_{K\in\pi}r_K^{n-p}\,\Ec(F;K)^p
\right\}^\frac1p.
$$
\par We also put $\varsigma:=\dlt/90$. Let us see that it
suffices to show that for every packing $\pi=\{K\}$ of
equal cubes of diameter $t\in(0,\varsigma)$ the following
inequality
\bel{EI-LAST}
I(F;\pi)\le C\,\lambda(f).
\ee
holds.
\par In fact, in this case, by Lemma \reff{SMALLQ},
$F\in\LOP$ and
$$
\|F\|_{\LOP}\le
C(\lambda(f)+\varsigma^{-1}\|F\|_{L_p(\RN)}).
$$
In turn, by Lemma \reff{FLP-MU}, $\|F\|_{L_p(\RN)}\le
C\,\lambda(f)$ so that $\|F\|_{\LOP}\le C\lambda(f)$ as
well. Hence,
$$
\|F\|_{\WP}=
\|F\|_{L_p(\RN)}+\|F\|_{\LOP}\le C\lambda(f),
$$
proving \rf{PRST}.
\par As above, we prove inequality \rf{EI-LAST} in five
steps.
\smallskip
\par {\it The first step: $\pi\subset\Kc_1$,} see \rf{K1}.
By Lemma \reff{P1AP},
$$
I(F;\pi)\le Ct^{-1}\AP(f;\gamma t:S)
$$
with $\gamma\le 574$. Since $t\le\varsigma=\dlt/90$, we
have $\gamma t\le\gamma\dlt/90\le 7\dlt\le\ve$ (recall that
$\dlt=\ve/2000$). Hence,
$$
I(F;\pi)\le C\sup_{0<t\le\ve}\frac{\AP(f;t:S)}{t}\le
C\lambda(f)
$$
so that \rf{EI-LAST} is satisfied.
\par There is no necessity to consider the case
$\pi\subset\Kc_2$. In fact, by \rf{K2}, every cube
$K\in\Kc_2$ has the diameter $\diam K>\dlt/82$ which
contradicts our assumption $\diam K\le\dlt/90, K\in\pi$.
\smallskip
\par {\it The third step: $\pi\subset\Kc_3$.} The proof of
inequality \rf{EI-LAST} in this case is based on the
following
\begin{lemma}\lbl{K3-STR} Let $\delta\in(0,1]$. Then for
every packing $\pi\subset\Kc_3$ of equal cubes and every
$0<\alpha <1/7$ we have
$$
I(F;\pi)^p\le C\,\intl_0^{\delta}
\TRS(f;t:S,\mu)^p\,\frac{dt}{t^{p+1}}.
$$
\end{lemma}
\par {\it Proof.} By Lemma \reff{P3}, there exists a finite
family $\pi'=\{Q\}\subset\TW$ of cubes of diameter $\diam
Q\le 2\delta$ such that
$$
I(F;\pi)^p\le C\sum\{r_{Q}^{n-p}|f(a_Q)-f(a_{Q'})|^p:
Q,Q'\in\pi', Q\cap Q'\ne\emptyset,r_{Q'}\le r_Q\}.
$$
\par Given $Q\in\pi'$ we put
$$ \bQ:=Q(a_Q,r_Q),~~~~~\widehat{Q}:=Q(a_Q,21r_Q). $$
By Lemma \reff{QP}, for every $Q'\in\pi'$ such that $Q'\cap
Q\ne\emp$ and $r_{Q'}\le r_Q$ we have
\bel{AQPR} \|a_Q-a_{Q'}\|\le 20r_Q. \ee
In addition, the lemma states the following: Let
$\tQ:=Q(a_Q,20r_Q)$. Then the cube $\eta\tQ$ is
$\beta$-porous for every $\eta\in(0,1]$ and $\beta \in
(0,3/20)$. Since $\widehat{Q}=(21/20)\tQ$, the cube
$\eta\widehat{Q}$ is $\beta$-porous for every
$\eta\in(0,1]$ but with $\beta < (20/21)3/20=1/7.$
\par Thus, by Definition \reff{STR-APS}, the cube
$\widehat{Q}$ is {\it strongly} $\beta$-porous for every
$\beta\in(0,1/7).$ Now we have
\be |f(a_Q)-f(a_{Q'})|&\le& |f(a_Q)-f_{\bQ,\mu}|+
|f_{\bQ,\mu}-f(a_{Q'})|\nn\\
&\le& \frac{1}{\mu(\bQ)}\intl_{\bQ} |f-f(a_Q)|\,d\mu+
\frac{1}{\mu(\bQ)}\intl_{\bQ} |f-f(a_{Q'})|\,d\mu. \nn\ee
Observe that $\bQ\subset\widehat{Q}$, see \rf{AQPR}, so
that, by \rf{ECT},
\be |f(a_Q)-f(a_{Q'})| &\le&
\frac{\mu(\widehat{Q})}{\mu(\bQ)}
\frac{1}{\mu(\widehat{Q})}\intl_{\widehat{Q}}
|f-f(x_{\widehat{Q}})|\,d\mu+
\frac{\mu(\widehat{Q'})}{\mu(\bQ)}
\frac{1}{\mu(\widehat{Q'})}\intl_{\widehat{Q'}}
|f-f(x_{\widehat{Q'}})|\,d\mu \nn\\
&=&\frac{\mu(\widehat{Q})}{\mu(\bQ)} \Ect(f;\widehat{Q})
+\frac{\mu(\widehat{Q'})}{\mu(\bQ)} \Ect(f;\widehat{Q'})\nn
\ee
\par Recall that for every cube $Q\in\Kc_3$ we have $\diam
Q\le 0.01\dlt$ so that
$$
\diam\widehat{Q}=21\diam Q\le 0.21\dlt\le 1.
$$
Therefore, by the doubling condition,
$\mu(\widehat{Q})/\mu(\bQ)\le C.$ Also, by \rf{AQPR},
$$
\widehat{Q'}=Q(a_{Q'},21r_{Q'}) \subset
Q(a_{Q'},21r_{Q})\subset Q(a_{Q},41r_{Q})=41\bQ.
$$
Since
$$
\diam 41\bQ=41\diam Q\le0.41\dlt\le 1,
$$
by the doubling condition
$$
\frac{\mu(\widehat{Q'})}{\mu(\bQ)}\le
\frac{\mu(41\bQ)}{\mu(\bQ)}\le C.
$$
Hence
$$
|f(a_Q)-f(a_{Q'})| \le C(\Ect(f;\widehat{Q})
+\Ect(f;\widehat{Q'}))
$$
so that
$$
I(F;\pi)^p\le C\sum\{r_{Q}^{n-p}(\Ect(f;\widehat{Q})^p
+\Ect(f;\widehat{Q'})^p): Q,Q'\in\pi', Q\cap
Q'\ne\emptyset,r_{Q'}\le r_Q\}.
$$
Since $r_Q\sim r_{Q'}$ for every $Q,Q'\in\TW$ such that
$Q\cap Q'\ne\emp$, we obtain
\be I(F;\pi)^p&\le& C\sum\{r_{Q}^{n-p}\Ect(f;\widehat{Q})^p
+r_{Q'}^{n-p}\Ect(f;\widehat{Q'})^p:
Q,Q'\in\pi'\}\nn\\
&=& 2C \sum_{Q\in\pi'}r_{Q}^{n-p}\Ect(f;\widehat{Q})^p.\nn
\ee
We also have
$$
\diam Q\le \diam \widehat{Q}=21\diam Q.
$$
Put $\varsigma:=21, T:=0.01\dlt$. Since $\dlt\le 1$, we
have $T\le 1/\varsigma$.  We also put
$$
\beta:=\min\{2\alpha,(\alpha+1/7)/2\}.
$$
Since $\alpha\in(0,1/7)$, we have $0<\alpha<\beta \le
2\alpha$ and $\beta\in(0,1/7)$. We have also proved above
that for every $Q\in\pi'$ the cube $\widehat{Q}$ is {\it
strongly} $\beta$-porous. \par Thus all conditions of Lemma
\reff{MU-APC} are satisfied. By this lemma,
$$
I(F;\pi)^p\le C^p
\sum_{Q\in\pi'}r_{Q}^{n-p}\Ect(f;\widehat{Q})^p \le
\intl_0^{4\varsigma T}
\TRS(f;t:S,\mu)^p\,\frac{dt}{t^{p+1}}.
$$
But $4\varsigma T=0.84\dlt\le\dlt$, and the proof of the
lemma is finished.\bx
\par Since $\dlt\le\ve$, by this lemma,
$$
I(F;\pi)^p\le C\,\intl_0^{\ve}
\TRS(f;t:S,\mu)^p\,\frac{dt}{t^{p+1}}\le
C\lambda(f)^p,
$$
proving inequality \rf{EI-LAST} for $\pi\subset\Kc_3$.
\par {\it The forth step: $\pi\subset\Kc_4$.} By Lemma
\reff{P4C0}, there exists a finite collection
$\pi'=\{Q\}\subset\TW$ of Whitney's cubes of diameter
$\diam Q\le 90\dlt$ such that
$$
I(F;\pi)^p\le C\sum_{Q\in\pi'}|Q||f(a_Q)|^p.
$$
Therefore, by Lemma \reff{FQ-EST},
$$
I(F;\pi)^p\le C\sum_{Q\in\TWB{90\dlt}}|Q||f(a_Q)|^p\le
C\lambda(f)^p,
$$
proving \rf{EI-LAST} for $\pi\subset\Kc_4$.
\par {\it The fifth step: $\pi\subset\Kc_5$.} This case is
trivial because, by Lemma \reff{P5}, $I(F;\pi)=0$ for every
packing $\pi\subset\Kc_5.$
\par Proposition \reff{EXT-J} is completely proved.\bx
\par We are in a position to finish the proof of Theorem
\reff{EXT-MU-NORM}. We observe that, by Proposition
\reff{EXT-J}, it remains to replace the functional
$\AP(f;\cdot:S)$ and $\TRS(f;\cdot:S,\mu)$ in the righthand
side of the proposition's inequality with the corresponding
functionals $\AP(f,\cdot:S)_{L_p(\mu)}$, see Definition
\reff{DEF-AP}, and $\PRM(f;\cdot:S)_{L_p(\mu)}$, see
Definition \reff{STR-F}. \par We will make this with the
help of the two following lemmas.
\begin{lemma}\lbl{J-AP} Let $s\in(\frac{n}{p},1)$ and let
$\alpha\in(0,1)$. Then for every $t\in(0,1]$ we have
$$
\TRS(f;t:S,\mu)\le C\,t^{s}\left(\intl_0^{2t}\left(
\frac{\PRM(f;u:S)_{L_p(\mu)}}
{u^s}\right)^p\,\frac{du}{u}\right)^{\frac1p}.
$$
Here $C$ is a constant depending only on $n,p,s$ and
$c_\mu$.
\end{lemma}
\par {\it Proof.} Let $Q$ be a strongly $\alpha$-porous
cube of diameter at most $t$ centered in $\partial S$, see
Definition \reff{STR-APS}. Put
$$
Q_i:=2^{-i}Q=Q(x_Q,2^{-i}r_Q),~~~i=0,1,...
$$
Then, by \rf{MO},
\be \Ect(f;Q)&:=&\frac{1}{\mu(Q)}
\intl_Q|f-f(x_Q)|\,d\mu\le\frac{1}{\mu(Q)}
\intl_Q|f-f_{Q,\mu}|\,d\mu+|f_{Q,\mu}-f(x_Q)|
\nn\\
&\le& \left(\frac{1}{\mu(Q)}
\intl_Q|f-{Q,\mu}|^p\,d\mu\right)^{\frac1p}
+|f_{Q,\mu}-f(x_Q)|\nn\\
&\le& C\Ec(f;Q)_{L_p(\mu)}+|f_{Q,\mu}-f(x_Q)|. \nn \ee
Since $f$ is continuous on $S$,
$$
\lim_{r\to 0}\frac{1}{\mu(rQ)} \intl_{rQ}f\,d\mu=f(x_Q),
$$
so that
$$
|f_{Q,\mu}-f(x_Q)|=
|\sum_{i=0}^\infty(f_{Q_i,\mu}-f_{Q_{i+1},\mu})| \le
\sum_{i=0}^\infty|f_{Q_i,\mu}-f_{Q_{i+1},\mu}|.
$$
Since $\diam Q_i=2^{-i}\diam Q\le t\le 1$, by the doubling
condition,
\be |f_{Q_i,\mu}-f_{Q_{i+1},\mu}|&\le&
\frac{1}{\mu(Q_{i+1})}
\intl_{Q_{i+1}}|f-f_{Q_i,\mu}|\,d\mu\nn\\
&\le&\frac{\mu(Q_{i})}{\mu(Q_{i+1})} \frac{1}{\mu(Q_{i})}
\intl_{Q_{i}}|f-f_{Q_i,\mu}|\,d\mu\le c_\mu
\frac{1}{\mu(Q_{i})} \intl_{Q_{i}}|f-f_{Q_i,\mu}|\,d\mu,
\nn\ee
so that, by \rf{MO},
$$
|f_{Q_i,\mu}-f_{Q_{i+1},\mu}|\le c_\mu
\left(\frac{1}{\mu(Q_{i})}
\intl_{Q_{i}}|f-f_{Q_i,\mu}|^p\,d\mu\right)^{\frac1p}\le
C\Ec(f;Q_i)_{L_p(\mu)}.
$$
We obtain
$$
|f_{Q,\mu}-f(x_Q)|\le\, C \sum_{i=0}^\infty
\Ec(f;Q_i)_{L_p(\mu)}.
$$
Hence, by the H\"older inequality,
$$
|f_{Q,\mu}-f(x_Q)|\le\, C \left(\sum_{i=0}^\infty
r_{Q_i}^{(s-n/p)p^*}\right)^{\frac{1}{p^*}}
\left(\sum_{i=0}^\infty \left(\frac{\Ec(f;Q_i)_{L_p(\mu)}}
{r_{Q_i}^{s-n/p}}\right)^p\right)^{\frac{1}{p}},
~~~~\frac1p+\frac{1}{p^*}=1.
$$
Since $r_{Q_i}=2^{-i}r_Q$,
$$
 \left(\sum_{i=0}^\infty
r_{Q_i}^{(s-n/p)p^*}\right)^{\frac{1}{p^*}}\le\,
Cr_Q^{s-n/p}.
$$
Hence,
$$
|f_{Q,\mu}-f(x_Q)|\le\, C
r_Q^{s-n/p}\left(\sum_{i=0}^\infty
\left(\frac{\Ec(f;Q_i)_{L_p(\mu)}}
{r_{Q_i}^{s-n/p}}\right)^p\right)^{\frac{1}{p}},
$$
so that
\be \Ect(f;Q)^p&\le&
C\{\Ec(f;Q)_{L_p(\mu)}^p+r_Q^{sp-n}\sum_{i=0}^\infty
|Q_i|\left(\frac{\Ec(f;Q_i)_{L_p(\mu)}}
{r_{Q_i}^{s}}\right)^p\}\nn\\&\le&
C\,r_Q^{sp-n}\sum_{i=0}^\infty
|Q_i|\left(\frac{\Ec(f;Q_i)_{L_p(\mu)}}
{r_{Q_i}^{s}}\right)^p. \nn \ee
We obtain
\bel{EQE} |Q|\Ect(f;Q)^p\le C\,r_Q^{sp}\sum_{i=0}^\infty
|Q_i|\left(\frac{\Ec(f;Q_i)_{L_p(\mu)}}
{r_{Q_i}^{s}}\right)^p. \ee
\par Let $\pi=\{Q\}$ be a packing of strongly
$\alpha$-porous equal cubes of diameter $\diam Q=\tau\le t$
centered in $\partial S$. Then, by Definition \reff{ECT},
for every $Q\in\pi$ and every $i=0,1,...$ the cube
$Q_i=2^{-i}Q$ is $\alpha$-porous. Now, by \rf{EQE},
\be
I&:=&\sum_{Q\in\pi} |Q|\Ect(f;Q)^p\le
C\,\sum_{Q\in\pi}(\diam Q)^{sp}\sum_{i=0}^\infty
|Q_i|\left(\frac{\Ec(f;Q_i)_{L_p(\mu)}}
{r_{Q_i}^{s}}\right)^p\nn\\
&=& C\,\tau^{sp}\sum_{Q\in\pi} \sum_{i=0}^\infty
|Q_i|\left(\frac{\Ec(f;Q_i)_{L_p(\mu)}}
{r_{Q_i}^{s}}\right)^p\nn\\
&=&C\,\tau^{sp}\sum_{i=0}^\infty(2^{-i}\tau)^{-sp}
\sum_{Q\in\pi} |2^{-i}Q|\Ec(f;2^{-i}Q)_{L_p(\mu)}^p
\nn
\ee
so that, by Definition \reff{AP-PRMU},
$$
I\le C\,\tau^{sp}\sum_{i=0}^\infty(2^{-i}\tau)^{-sp}
\PRM(f;2^{-i}\tau:S)_{L_p(\mu)}^p.
$$
Since $\PRM(f;\cdot:S)_{L_p(\mu)}$ is a non-decreasing
function, we obtain
$$
I\le C\,\tau^{sp}\intl_0^{2\tau}\left(
\frac{\PRM(f;u:S)_{L_p(\mu)}}
{u^s}\right)^p\,\frac{du}{u}\le C\,t^{sp}\intl_0^{2t}\left(
\frac{\PRM(f;u:S)_{L_p(\mu)}} {u^s}\right)^p\,\frac{du}{u}.
$$
It remains to take the supremum over all packings $\pi$
satisfying conditions of Definition \reff{STR-F}, and the
result of the lemma follows.\bx
\par In a similar way one can prove the following estimate
of the functional $\AP(f;\cdot:S)$ via the corresponding
functional $\AP(f;\cdot:S)_{L_p(\mu)}$, see Definitions
\reff{DEF-AP} and \reff{DEF-APLS}. (Recall that $
\AP(f;t:S)= \AP(f;t:S)_{L_\infty(\mu)}$ provided
$S=\supp\mu$.)
\begin{lemma}\lbl{AP-APMU}For every  $s\in(\frac{n}{p},1)$
the following inequality
$$
\AP(f;t:S)^p\le C\,t^{s}\left(\intl\limits_0^{2t}\left(
\frac{\AP(f,u:S)_{L_p(\mu)}}
{u^s}\right)^p\,\frac{du}{u}\right)^{\frac1p},~~~~0<t\le 1,
$$
holds. Here $C$ is a constant depending only on $n,p,s$ and
$c_\mu$.
\end{lemma}
\par Now, we put $s:=(\frac{n}{p}+1)/2$ in Lemmas
\reff{J-AP} and \reff{AP-APMU} and then apply the Hardy
inequality \rf{HE1}. We obtain:
$$
\intl\limits_0^{\ve} \TRS(f;t:S,\mu)^p\,
\frac{dt}{t^{p+1}}\le C\,\intl\limits_0^\ve
\PRM(f;t:S)_{L_p(\mu)}^p\, \frac{dt}{t^{p+1}}
$$
and
$$
\sup_{0<t\le \ve} \frac{\AP(f;t:S)}{t}\le C\,\sup_{0<t\le
\ve} \frac{\AP(f;t:S)_{L_p(\mu)}}{t}
$$
with $C$ depending only on $n,p$ and $c_\mu$. Combining
these inequalities with the result of Proposition
\reff{EXT-J}, we obtain the required inequality \rf{OPIN}.
\par Theorem \reff{EXT-MU-NORM} is proved.\bx
\smallskip
\par This theorem leads us to a constructive description of
the trace space $\WP|_S$. To formulate this result given
$x,y\in S$ and $\alpha\in (0,1]$ we put
\bel{DXYA} \rxy:=\inf\{\diam Q:~Q\ni x,y, ~\alpha Q
\subset\RN\setminus S\}. \ee
We put $\rxy:=+\infty$ whenever there no exists a cube $Q$
such that $x,y\in Q$ and $\alpha Q\subset\RN\setminus S$.
\par Clearly, the function $\rho_S$ defined by \rf{DX15}
coincides with $\rho_{\alpha,S}$ where $\alpha=
\frac{1}{15}$.
\par We prove the following generalization of part (ii) of
Theorem \reff{EXT-DIFF-S15}.
\begin{theorem}\lbl{EXT-DIFF} Let $p\in(n,\infty),$
$\ve>0$, and let $\alpha\in(0,1/14)$. Then for every
function $f\in C(S)$ the following equivalence
\be \|f\|_{\WP|_S}&\sim& \|f\|_{L_p(\mu)}+ \sup_{0<t\le
\ve} \left(~\iint\limits_{\|x-y\|<t} \frac{|f(x)-f(y)|^p}
{t^{p-n}\mu(Q(x,t))\mu(Q(y,t))}
\,d\mu(x)\,d\mu(y)\right)^{\frac1p}\nn\\
&+& \left(~\iint\limits_{\rxy<\ve} \frac{|f(x)-f(y)|^p}
{\rxy^{p-n}\mu(Q(x,\rxy))^2} \,d\mu(x)\,d\mu(y)
\right)^{\frac1p} \nn\ee
holds with constants depending only on $n,p,\ve,\alpha$,
and parameters $C_\mu,C'_\mu,C''_\mu$.
\end{theorem}
\par The proof of the theorem relies on results of Sections
6 and 7 and the following two lemmas.
\begin{lemma}\lbl{WPN} Let $p\in(n,\infty),$
$\alpha\in(0,1)$ and $\ve\in(0,1/4]$. Then for every
function $f\in C(S)$ we have
\be
1/C\intl_0^{\alpha\ve}\PRMA{2\alpha}(f;t:S)_{L_p(\mu)}^p\,
\frac{dt}{t^{p+1}} &\le& \iint\limits_{\rxy<\ve}
\frac{|f(x)-f(y)|^p} {\rxy^{p-n}
\mu(Q(x,\rxy))^2}\,d\mu(x)\,d\mu(y)\nn\\
&\le& C\intl_0^{4\ve}\PRMA{\alpha/8}(f;t:S)_{L_p(\mu)}^p\,
\frac{dt}{t^{p+1}}. \nn\ee
Here $C$ is a positive constant depending only on
$n,p,\ve,\alpha,$ and the doubling constant $c_\mu$.
\end{lemma}
\par {\it Proof.} Let us prove the first inequality of the
lemma. We put $t_i:=2^{-i}\alpha\ve, i=0,1,...\,\,$. We
have
\be I_1&:=&\intl_0^{\alpha\ve}\left
(\frac{\PRMA{2\alpha}(f;t:S)_{L_p(\mu)}} {t}\right)^p\,
\frac{dt}{t}= \sum_{i=0}^\infty\intl_{t_{i+1}}^{t_i}\left
(\frac{\PRMA{2\alpha}(f;t:S)_{L_p(\mu)}} {t}\right)^p\,
\frac{dt}{t}\nn\\&\le& C\sum_{i=0}^\infty t_i^{-p}
\PRMA{2\alpha}(f;t_i:S)_{L_p(\mu)}^p. \nn \ee
\par By Definition \reff{AP-PRMU}, for every $i=0,1,...$
there exists a packing $\pi_i=\{Q\}$ of $2\alpha$-porous
equal cubes centered in $\partial S$ with diameter
\bel{TI}
\diam Q=\tau_i\le t_i=2^{-i}\alpha\ve
\ee
such that
$$
\PRMA{2\alpha}(f;t_i:S)_{L_p(\mu)}^p\le
2\sum_{Q\in\pi_i}|Q|\Ec(f;Q)^p_{L_p(\mu)}.
$$
We may assume that $\{\tau_i:~i=0,1,...\}$ is
non-decreasing sequence. By \rf{EMU},
\be I_1&\le& C\sum_{i=0}^\infty \sum_{Q\in\pi_i}
t_i^{-p}r_Q^{n}\Ec(f;Q)^p_{L_p(\mu)}\nn\\&=&
C\sum_{i=0}^\infty \sum_{Q\in\pi_i}\frac{
t_i^{-p}r_Q^{n}}{\mu(Q)^2} \iint\limits_{Q\times
Q}|f(x)-f(y)|^p\, d\mu(y)d\mu(x). \nn\ee
\par We let $H$ denote a subset of $\RN\times \RN$ defined
by the formula
$$
H:=\cup\{Q\times Q:~Q\in\pi_i,~ i=0,1,2,...\}.
$$
Also, given $(x,y)\in H, x\ne y$, we define a collection of
cubes $\Qc_{xy}$ by letting
$$
\Qc_{xy}:=\{Q:~\exists\, i\in\{0,1,2,...\} {\rm
~such~that~~} Q\in\pi_i~~{\rm~and~} Q\ni x,y \}.
$$
Then
$$
I_1 \le C\iint\limits_H |f(x)-f(y)|^p\,
g(x,y)d\mu(y)d\mu(x).
$$
where
$$
g(x,y):=\sum\{t_i^{-p}r_{Q_i}^{n}\mu(Q_i)^{-2}:~
Q_i\in\pi_i, Q_i\ni x,y,~i=0,1,...\}.
$$
\par Fix $(x,y)\in H$ and $Q\in\Qc_{xy}$; thus $Q\in\pi_i$
for some $i=0,1,...,$ and $x,y\in Q$. Recall that $x_Q\in
S$, $\diam Q\le 2^{-i}\ve$, and $Q$ is $2\alpha$-porous.
The latter implies the existence of a cube $\tQ\subset
Q\setminus S$ such that $\diam Q\le (2\alpha)^{-1}\diam
\tQ$. Then
$$
Q\subset 2(2\alpha)^{-1}\tQ=\alpha^{-1}\tQ,
$$
so that $x,y\in \alpha^{-1}\tQ$. Therefore , by \rf{DXYA},
$\rxy\le \alpha^{-1}\diam \tQ$. Note also that $x_Q\in S$
and $\tQ\subset Q\setminus S$ so that $\diam \tQ\le
\frac12\diam Q$. Hence,
\bel{QDXY} \rxy\le \alpha^{-1}\diam \tQ\le
\tfrac12\,\alpha^{-1}\diam Q\le 2^{-i-1}\ve. \ee
In particular, $\rxy<\ve$ so that
$$
I_1 \le C\iint\limits_{\rxy<\ve} |f(x)-f(y)|^p\,
g(x,y)\, d\mu(y)d\mu(x).
$$
\par Let us estimate $g(x,y)$. Since $\diam
Q_i\ge\|x-y\|>0$ for each $Q_i\ni x,y$,
$$
i_0:=\sup\{i:~ \exists \,Q_i\in\pi_i, Q_i\ni x,y\}<\infty.
$$
\par In addition, for each $ Q_i\in\pi_i, Q_i\ni x,y,$ we
have $\diam Q_i\ge\diam Q_{i_0}$ and $Q_i\cap
Q_{i_0}\ne\emp$ so that $2Q_i\supset Q_{i_0}$. On the other
hand, since $\ve\le 1/4$,
$$
r_{Q_i}\le t_i/2=2^{-i-1}\alpha\ve\le \alpha\ve/2\le
1/2,~~i=1,2,...\,\,,
$$
so that, by the doubling condition, $
\mu(Q_{i_0})\le\mu(2Q_i)\le C\mu(Q_i).$ Hence, by \rf{TI},
\be g(x,y)&:=&\sum\{t_i^{-p}r_{Q_i}^{n}\mu(Q_i)^{-2}:~
Q_i\in\pi_i, Q_i\ni x,y,~i=0,1,...\}\nn\\
&\le& C\,\mu(Q_{i_0})^{-2}\sum\{t_i^{-p}r_{Q_i}^{n}:~
Q_i\in\pi_i, Q_i\ni x,y,~i=0,1,...\}\nn\\
&\le& C\,\mu(Q_{i_0})^{-2}\sum\{t_i^{-p}t_i^{n}:~
Q_i\in\pi_i, Q_i\ni x,y,~i=0,1,...\}\nn\\
&\le& C\,\mu(Q_{i_0})^{-2}\sum\{t_i^{n-p}:~
i\le i_0\}\nn\\
&\le& C\,\mu(Q_{i_0})^{-2}t_0^{n-p}\le
C\,\mu(Q_{i_0})^{-2}r_{Q_{i_0}}^{n-p}\nn \ee
Observe that, by \rf{QDXY},
\bel{EDXY} \rxy\le \tfrac12\,\alpha^{-1}\diam
Q_{i_0}=r_{Q_{i_0}}/\alpha. \ee
Since $Q_{i_0}\ni x$, we obtain
$$
Q(x,\rxy)\subset (1+\alpha^{-1})Q_{i_0}\subset
2\alpha^{-1}Q_{i_0}.
$$
But
$$
\diam(2\alpha^{-1}Q_{i_0})=2\alpha^{-1}\diam Q_{i_0}
\le 2\alpha^{-1}t_{i_0} \le 2\ve\le 1,
$$
so that by the doubling condition
$$
\mu(Q(x,\rxy))\le\mu(2\alpha^{-1}Q_{i_0})\le
C\mu(Q_{i_0}).
$$
Combining this estimate with \rf{EDXY}, we obtain
$$
g(x,y)\le Cr_{Q_{i_0}}^{n-p}\mu(Q_{i_0})^{-2}\le
C\,\rxy^{n-p} \mu(Q(x,\rxy))^{-2}
$$
proving that
$$
I_1\le C\, \iint\limits_{\rxy<\ve}
\frac{|f(x)-f(y)|^p} {\rxy^{p-n}
\mu(Q(x,\rxy))^2}\,d\mu(x)\,d\mu(y).
$$
\par Prove the second inequality of the lemma. We have
\be I_2&:=&\iint\limits_{\rxy<\ve} \frac{|f(x)-f(y)|^p}
{\rxy^{p-n}\mu(Q(x,\rxy))^2}\,d\mu(x)\,d\mu(y)\nn\\&=&
\sum_{i=0}^\infty ~~\iint\limits_{t_{i+1}\le \rxy<t_i}
\frac{|f(x)-f(y)|^p} {\rxy^{p-n}\mu(Q(x,\rxy))^2}
\,d\mu(x)\,d\mu(y), \nn \ee
where $t_i:=2^{-i}\ve,$ $i=0,1,...$. Put
$$
A_i:=\iint\limits_{t_{i+1}\le \rxy<t_i}
\frac{|f(x)-f(y)|^p} {\rxy^{p-n}
\mu(Q(x,\rxy))^2}\,d\mu(x)\,d\mu(y)
$$
and prove that
\bel{EAPA} A_i\le C
t_i^{-p}\PRMA{\alpha/4}(f;2t_i:S)^p_{L_p(\mu)}. \ee
\par To this end we let $S_i$ denote the projection of the
set
$$
\{(x,y)\in S\times S:~t_{i+1}\le \rxy<t_i\}
$$
onto $S$. Thus
$$
S_i:=\{x\in S:~\exists\, x'\in S ~~{\rm
such~that}~~t_{i+1}\le
\rxx<t_i\},~~~~i=0,1,...\,.
$$
Then, by \rf{DXYA}, for each $x\in S_i$ there exists $x'\in
S$ and a cube $Q^{(x)}$ with
$$ t_{i+1}\le \diam Q^{(x)}<t_i $$
such that $\alpha Q^{(x)}\subset\RN\setminus S$ and
$$ x,x'\in Q^{(x)}. $$
\par Put $K^{(x)}:=Q(x,t_i)$. Since $x\in Q^{(x)}$ and
$\diam Q^{(x)}<t_i=\diam K^{(x)}$, we have $K^{(x)}\supset
Q^{(x)}$. In addition, since $\alpha Q^{(x)}\subset
Q^{(x)}\setminus S\subset K^{(x)}\setminus S$ and
$$
\diam (\alpha Q^{(x)})\ge \alpha t_{i+1}= \alpha
t_i/2=\tfrac14\,\alpha\diam K^{(x)},
$$
by Definition \reff{PR-Q}, $K^{(x)}$ is an
$\alpha/4$-porous cube (with respect to $S$).
\par By the Besicovitch covering theorem, see , e.g.
\cite{G}, there exists a countable subfamily
$$
V_i:=\{K_j=Q(x_j,t_i):~j\in J\}
$$
of the family of cubes $\{K^{(x)}:~x\in S_i\}$ which covers
$S_i$ and has the covering multiplicity $M_{V_i}\le C(n)$.
Consider a cube $K_j\in V_i$. Prove that if $x\in K_j$ and
$\rxy< t_i$, then
\bel{YKJ} y\in 2K_j. \ee
In fact, by \rf{DXYA}, there exists a cube $Q$ of $\diam
Q\le t_i$ such that $x,y\in Q$ and $\alpha
Q\subset\RN\setminus S$. Hence $\|x-y\|\le\diam Q\le t_i$
so that
$$
\|x_j-y\|\le \|x_j-x\|+\|x-y\|\le t_i+t_i=2t_i,
$$
proving \rf{YKJ}.
\par Also note that for each $x\in K_j$ we have
$$
2K_j\subset Q(x,2\diam K_j)=Q(x,4t_i) =Q(x,8t_{i+1}).
$$
But $8t_{i+1}\le 4\ve\le 1$ provided $\ve\le 1/4$ so that,
by the doubling condition, $\mu(2K_j)\le
C\mu(Q(x,t_{i+1})).$
\par Now we have
\be A_i&\le& \sum_{j\in J}~\intl_{K_j}~~\intl_{t_{i+1}\le
\rxy<t_i} \frac{|f(x)-f(y)|^p} {\rxy^{p-n}
\mu(Q(x,\rxy))^2}\,d\mu(y)\,d\mu(x)\nn\\
&\le& \sum_{j\in J}~\intl_{K_j}~~\intl_{t_{i+1}\le
\rxy<t_i} \frac{|f(x)-f(y)|^p}
{t_{i+1}^{p-n}\mu(Q(x,t_{i+1}))^2}\,d\mu(y)\,d\mu(x)\nn\\
&\le& C\,\sum_{j\in J} t_{i}^{n-p}\mu(2K_j)^{-2}
~\intl_{K_j}~~\intl_{t_{i+1}\le \rxy<t_i}
|f(x)-f(y)|^p\,d\mu(y)\,d\mu(x). \nn\ee
Hence, by \rf{YKJ}, we have
$$
A_i\le  C\,\sum_{j\in J} t_{i}^{n-p}\mu(2K_j)^{-2}
~\intl_{2K_j}\intl_{2K_j} |f(x)-f(y)|^p\,d\mu(x)\,d\mu(y).
$$
This and \rf{EMU} yield
\be A_i&\le&  C\,\sum_{j\in J}
t_{i}^{-p}|2K_j|\frac{1}{\mu(2K_j)^{-2}}
\iint\limits_{2K_j\times 2K_j}
|f(x)-f(y)|^p\,d\mu(x)\,d\mu(y)\nn\\
&\le&  C\,t_{i}^{-p}\sum_{j\in J}
|2K_j|\Ec(f;2K_j)^p_{L_p(\mu)}. \nn\ee
\par Recall that $V_i$ is a collection of equal
$\alpha/4$-porous cubes of diameter $t_i$ centered in $S$
with the covering multiplicity $M_{V_i}\le C(n).$ Then the
family $2V_i:=\{2K_j:~K_j\in V_i\}$ is a collection of
equal $\alpha/8$-porous cubes of diameter $2t_i$ with
centers in $S$.
\par It can be readily seen that the covering multiplicity
of the family $2V_i$ is also bounded by a constant
$C=C(n)$. Then, by Theorem \reff{TFM}, $2V_i$ can be
partitioned into at most $N(n)$ families of pairwise
disjoint cubes which allows us to assume that the cubes of
the collection $2V_i$ are pairwise disjoint. Hence,
$$
A_i\le  C\,t_{i}^{-p}\sup_\pi\sum_{K\in\pi}
|K|\Ec(f;K)^p_{L_p(\mu)}
$$
where $\pi$ runs over all {\it finite} subfamilies of
$2V_i$. Since every such $\pi$ is a packing, by Definition
\reff{AP-PRMU}, we have
$$
A_i\le C t_i^{-p} \PRMA{\alpha/8}(f;2t_i:S)^p_{L_p(\mu)}
$$
proving \rf{EAPA}.
\par Finally,
$$
I_2\le \sum_{i=0}^\infty A_i \le C\sum_{i=0}^\infty \left
(\frac{\PRMA{\alpha/8}(f;u_i:S)_{L_p(\mu)}}{u_i}\right)^p,
~~{\rm where}~~u_i:=2t_i, i=0,1,... .
$$
Since $\PRMA{\alpha/8}(f;\cdot:S)_{L_p(\mu)}$ is a
non-decreasing function, this implies the required
inequality
$$
I_2\le C\intl_0^{4\ve}\PRMA{\alpha/8}(f;u:S)_{L_p(\mu)}^p\,
\frac{du}{u^{p+1}},
$$
proving the second inequality of the lemma.\bx
\par Given a function $f\in C(S)$, $\alpha\in(0,1)$ and
$\ve>0$ we put
\bel{DJ1} J_1(f;t):=\left(~\iint\limits_{\|x-y\|<t}
\frac{|f(x)-f(y)|^p} {t^{p-n}\mu(Q(x,t))^2}
\,d\mu(x)\,d\mu(y)\right)^{\frac1p},~~~~t>0, \ee
and
\bel{DJ2} J_2(f;\alpha,\ve):=\left(~\iint\limits_{\rxy<\ve}
\frac{|f(x)-f(y)|^p} {\rxy^{p-n}\mu(Q(x,\rxy))^2}
\,d\mu(x)\,d\mu(y) \right)^{\frac1p}. \ee
\begin{lemma}\lbl{JEBIG} Let $p\in(n,\infty),$
$\alpha\in(0,1),$ and let $0<\ve'<\ve$. Then for every
function $f\in C(S)$ we have
\bel{J1BG}
J_1(f;t)\le
C\|f\|_{L_p(\mu)},~~~~t\in(\ve',\ve),
\ee
and
\bel{J2BG}
J_2(f;\alpha,\ve')\le J_2(f;\alpha,\ve)\le
C\{J_2(f;\alpha,\ve')+\|f\|_{L_p(\mu)}\}. \ee
Here $C$ is a positive constant depending only on
$n,p,\ve,\ve',$ and the constants $C_\mu,C'_\mu,C''_\mu$.
\end{lemma}
\par {\it Proof.} Put $\tve:=\min\{1,\ve'\}$. Then, by the
doubling condition,
$$
\mu(Q(x,1))\le C(\ve',n)
\mu(Q(x,\tve))\le C\mu(Q(x,\ve')),~~~ x\in
S.
$$
In turn, by \rf{B1}, $\mu(Q(x,1))\ge C'_\mu$ so that
\bel{EUPM} \mu(Q(x,\ve'))^{-1}\le C /C'_\mu. \ee
Hence, for every $t>\ve$ we have
\be J_1(f;t)^p&:=&t^{n-p}~\iint\limits_{\|x-y\|<t}
\frac{|f(x)-f(y)|^p} {\mu(Q(x,t))^2}
\,d\mu(x)\,d\mu(y)\nn\\&\le&
(\ve')^{n-p}~\iint\limits_{\|x-y\|<\ve}
\frac{|f(x)-f(y)|^p} {\mu(Q(x,\ve'))^2}
\,d\mu(x)\,d\mu(y)\nn\\&\le& C~\iint\limits_{\|x-y\|<\ve}
|f(x)-f(y)|^p \,d\mu(x)\,d\mu(y)\nn\\&\le&
C~\iint\limits_{\|x-y\|<\ve} (|f(x)|^p+|f(y)|^p)
\,d\mu(x)\,d\mu(y). \nn\ee
But
\be I_1&:=&\iint\limits_{\|x-y\|<\ve} (|f(x)|^p+|f(y)|^p)
\,d\mu(x)\,d\mu(y)=2\iint\limits_{\|x-y\|<\ve} |f(x)|^p
\,d\mu(y)\,d\mu(x)\nn\\&\le& 2\intl_S\intl_{y\in Q(x,\ve)}
|f(x)|^p \,d\mu(y)\,d\mu(x)=2\intl_S \mu(Q(x,\ve))|f(x)|^p
\,d\mu(x). \nn\ee
If $\ve \le 1$, then, by \rf{B1}, we have $\mu(Q(x,\ve))\le
\mu(Q(x,1))\le C''_\mu$. If $\ve>1$, then we can cover the
set $Q(x,\ve)\cap S$ by a finite family of cubes
$\pi=\{Q(x_i,1):~x_i\in Q(x,\ve)\cap S,
i=1,...,N(n,\ve)\}$, so that
$$
\mu(Q(x,\ve))\le\sum_{i=1}^{N(n,\ve)}
\mu(Q(x_i,1))\le N(n,\ve)C''_\mu.
$$
Hence,
\bel{JLP}
I_1\le 2\intl_S \mu(Q(x,\ve))|f(x)|^p \,d\mu(x)
\le C\|f\|_{L_p(\mu)}^p,
\ee
so that
$$
J_1(f;t)^p\le C\,I_1\le C\|f\|_{L_p(\mu)}^p,
~~~~t\in(\ve',\ve),
$$
proving inequality \rf{J1BG}.
\par Prove inequality \rf{J2BG}. The first inequality of
\rf{J2BG} is trivial. Let us prove the second one.
\par We have
$J_2(f;\alpha,\ve)^p=J_2(f;\alpha,\ve')^p+I_2,$ where
$$
I_2:=\iint\limits_{\ve'<\rxy<\ve}
\frac{|f(x)-f(y)|^p}
{\rxy^{p-n}\mu(Q(x,\rxy))^2}
\,d\mu(x)\,d\mu(y).
$$
By \rf{EUPM},
\be I_2&\le& \iint\limits_{\ve'<\rxy<\ve}
\frac{|f(x)-f(y)|^p} {(\ve')^{p-n}\mu(Q(x,\ve'))^2}
\,d\mu(x)\,d\mu(y)\nn\\&\le& C\iint\limits_{\ve'<\rxy<\ve}
|f(x)-f(y)|^p\,d\mu(x)\,d\mu(y)\le C\iint\limits_{\rxy<\ve}
|f(x)-f(y)|^p\,d\mu(x)\,d\mu(y). \nn\ee
By \rf{DXYA}, $\rxy\ge\|x-y\|$, so that
$$ I_2\le
C\iint\limits_{\|x-y\|<\ve}
|f(x)-f(y)|^p\,d\mu(x)\,d\mu(y)\le
C\iint\limits_{\|x-y\|<\ve}
(|f(x)|^p+|f(y)|^p)\,d\mu(x)\,d\mu(y)=C\,I_1.
$$
Hence, by \rf{JLP}, $I_2\le C\,I_1\le C\|f\|_{L_p(\mu)}^p$,
proving that
$$
J_2(f;\alpha,\ve)^p\le
J_2(f;\alpha,\ve')+C\|f\|_{L_p(\mu)}^p.
$$
\par The lemma is proved.\bx
\medskip
\par {\bf Proof of Theorem \reff{EXT-DIFF}.} The proof
easily follows from Theorem \reff{EXT-MU-NORM}, Proposition
\reff{FAPN} and Lemma \reff{WPN}. In fact, given a function
$f\in C(S)$, $\alpha\in(0,1)$ and $\ve>0$ we put
$$
J_3(\alpha,\ve):=\left\{\intl_0^\ve
\PRM(f;t:S)_{L_p(\mu)}^p\,
\frac{dt}{t^{p+1}}\right\}^{\frac{1}{p}}.
$$
\par Then, by Theorem \reff{EXT-MU-NORM},
\bel{EI0} \|f\|_{\WP|_S}\sim \|f\|_{L_p(\mu)}+ \sup_{0<t\le
\ve} \frac{\AP(f;t:S)_{L_p(\mu)}}{t} + J_3(\beta,\ve) \ee
provided $p\in(n,\infty)$, $\beta\in(0,1/7)$ and $\ve>0$.
In turn, by Proposition \reff{FAPN}, for every
$t\in(0,1/4]$ we have
$$
\frac1C\AP(f;t/4:S)_{L_p(\mu)}/t\le J_1(f;t)\le
C\AP(f;4t:S)_{L_p(\mu)}/t,
$$
where $J_1(f;t)$ is defined by \rf{DJ1}. Taking the
supremum in this inequality over all $t\in(0,\ve)$ where
$0<\ve\le1/4$, we obtain
$$
C^{-1}\sup_{0<t<\ve}\AP(f;t/4:S)_{L_p(\mu)}/t\le
\sup_{0<t<\ve}J_1(f;t)\le
C\sup_{0<t<\ve}\AP(f;4t:S)_{L_p(\mu)}/t,
$$
so that
$$
C^{-1}\sup_{0<t<\ve/4}\AP(f;t:S)_{L_p(\mu)}/t\le
\sup_{0<t<\ve}J_1(f;t)\le
C\sup_{0<t<4\ve}\AP(f;t:S)_{L_p(\mu)}/t.
$$
Applying Lemma \reff{APLP}, we obtain
\bel{E1}
\|f\|_{L_p(\mu)}+\sup_{0<t<\ve}\AP(f;t:S)_{L_p(\mu)}/t\sim
\|f\|_{L_p(\mu)}+\sup_{0<t<\ve}J_1(f;t),~~~~\ve\in(0,1/4].
\ee
\par In turn, by Lemma \reff{WPN}, we have
\bel{J12} C^{-1}J_3(2\alpha,\alpha\ve)\le
J_2(f;\alpha,\ve)\le C J_3(\alpha/8,4\ve). \ee
(Recall that $J_2(f;\alpha,\ve)$ is defined by \rf{DJ2}.)
On the other hand, by \rf{Y} and Lemma \reff{APLP},
$$ \PRM(f;t:S)_{L_p(\mu)}\le \AP(f;t:S)_{L_p(\mu)}\le
C(1+t^{\frac{n}{p}})\|f\|_{L_p(\mu)},$$
so that for every $\alpha\in(0,1)$ we have
$$
J_3(\alpha/8,4\ve)^p=J_3(\alpha/8,\ve)^p+
\intl_{\ve}^{4\ve}
\Ac_{p,\alpha/8}(f;t:S)_{L_p(\mu)}^p\,
\frac{dt}{t^{p+1}}\le
J_3(\alpha/8,\ve)^p+C\|f\|_{L_p(\mu)},
$$
and
$$
J_3(\alpha,\ve)^p=J_3(\alpha,\alpha\ve)^p+
\intl_{\alpha\ve}^{\ve}\PRM(f;t:S)_{L_p(\mu)}^p\,
\frac{dt}{t^{p+1}}\le
J_3(\alpha,\alpha\ve)^p+C\|f\|_{L_p(\mu)}.
$$
These inequalities and \rf{J12} imply
\bel{JFIN} C^{-1}(\|f\|_{L_p(\mu)}+J_3(2\alpha,\ve))\le
\|f\|_{L_p(\mu)}+J_2(f;\alpha,\ve)\le C(\|f\|_{L_p(\mu)}+
J_3(\alpha/8,\ve)). \ee
Observe that $\alpha/8, 2\alpha\in (0,1/7)$, so that, by
\rf{EI0}, we have
\be
\|f\|_{\WP|_S}&\sim& \|f\|_{L_p(\mu)}+ \sup_{0<t\le\ve}
\frac{\AP(f;t:S)_{L_p(\mu)}}{t} + J_3(2\alpha,\ve)\nn\\
&\sim& \|f\|_{L_p(\mu)}+ \sup_{0<t\le \ve}
\frac{\AP(f;t:S)_{L_p(\mu)}}{t} + J_3(\alpha/8,\ve).\nn \ee
These equivalences, \rf{E1} and inequalities \rf{JFIN}
yield
\bel{IE1}
\|f\|_{\WP|_S}\sim \|f\|_{L_p(\mu)}+
\sup_{0<t\le\ve} J_1(f;t) + J_2(f;\alpha,\ve)
\ee
proving Theorem \reff{EXT-DIFF} for $\ve\in(0,1/4]$.
\par Consider the case $\ve>\ve':=1/4$. By Lemma
\reff{JEBIG},
$$
\sup_{0<t\le\ve}J_1(f;t)\le C\{\sup_{0<t\le 1/4} J_1(f;t)
+\|f\|_{L_p(\mu)}\},
$$
and
$$
J_2(f;\alpha,1/4)\le J_2(f;\alpha,\ve)\le
C\{J_2(f;\alpha,1/4) +\|f\|_{L_p(\mu)}\},
$$
proving that equivalence \rf{IE1} holds for $\ve>1/4$ as
well.
\par Theorem \reff{EXT-DIFF} is completely proved.
\medskip
\par This theorem and Proposition \reff{WDS} imply the
following generalization of part (i) of Theorem
\reff{EXT-DIFF-S15}.
\begin{theorem}\lbl{EXT-DIFF-DSGM} Let $p\in(n,\infty)$,
$\ve>0$ and let $\alpha\in(0,1/14)$. Let $\sigma$ be a
measure supported on $\partial S$ which satisfies
conditions \rf{DNSG} and \rf{B1SIGMA}, and let
$\Omega:=int(S)$ be the interior of $S$. Then for every
$f\in C(S)$ the following equivalence
\be
\|f\|_{\WP|_S}&\sim& \|f|\,
_\Omega\|_{W^1_p(\Omega)}+\|f\|_{L_p(\sigma)}+ \sup_{0<t\le
\ve} \left(\,\,\,\iint\limits_{\|x-y\|<t}
\frac{|f(x)-f(y)|^p} {t^{p-n}\sigma(Q(x,t))^2}
\,d\sigma(x)d\sigma(y)\right)^{\frac1p}\nn\\
&+& \left(~\iint\limits_{\rho_{\alpha,\partial S}(x,y)<\ve}
\frac{|f(x)-f(y)|^p} {\rho_{\alpha,\partial
S}(x,y)^{p-n}\sigma(Q(x,\rho_{\alpha,\partial S}(x,y)))^2}
\,d\sigma(x)\,d\sigma(y) \right)^{\frac1p}
\nn
\ee
holds with constants depending only on $n,p,$ and
$C_\sigma,C'_\sigma,C''_\sigma$.
\end{theorem}
\par {\it Proof.} Let us apply Theorem \reff{EXT-DIFF} to
the set $\partial S$ and the measure $\sigma$ supported on
$\partial S$. By this theorem, the quantity $\|f|_{\partial
S}\|_{\WP|_{\partial S}}$ is equivalent to the sum of the
last three terms on the right-hand side of the equivalence
of Theorem \reff{EXT-DIFF-DSGM}. On the other hand, by
Proposition \reff{WDS},
$$
\|f\|_{\WP|_S}\sim \|f|_\Omega\|_{\WPO}+
\|f|_{\partial S}\|_{\WP|_{\partial S}},
$$
and the proof is finished.\bx
\begin{remark}\lbl{RH-S1} {\em The theorem remains true
whenever the function $\rho_{\alpha,\partial S}$ in its
formulation is replaced with the (bigger) function
$\rho_{\alpha,S}$. In fact, Remark \reff{EST-AFT} shows
that we can prove Theorem \reff{EXT-DIFF-DSGM} basing on
the remark's equivalence rather than on Proposition
\reff{WDS}. This allows us to accomplish the proof using
the functional $\PRM(f;\cdot:S)$ rather than
$\PRM(f;\cdot:\partial S)$ so that, by Lemma \reff{WPN},
the proof can be continued with the function
$\rho_{\alpha,S}$ rather than $\rho_{\alpha,\partial S}$.}
\end{remark}
\par The next result presents a generalization of Theorem
\reff{EXT-DIFF-15}.
\begin{theorem}\lbl{EXT-DF-LA} Let $p\in(n,\infty)$,
$\alpha\in(0,1/14)$, $\ve>0$, and let $\sigma$ be a measure
supported on $\partial S$ which satisfies conditions
\rf{DNSG} and \rf{B1SIGMA}. Then for every $f\in C(S)$ we
have
\be \|f\|_{\WP|_S}&\sim& \|f\|_{L_p(\mu)}+ \sup_{0<t\le
\ve} \left(~\iint\limits_{\|x-y\|<t} \frac{|f(x)-f(y)|^p}
{t^{p-n}\mu(Q(x,t))^2}
\,d\mu(x)\,d\mu(y)\right)^{\frac1p}+\|f\|_{L_p(\sigma)}\nn\\
&+& \left(~\iint\limits_{\rho_{\alpha,\partial S}(x,y)<\ve}
\frac{|f(x)-f(y)|^p} {\rho_{\alpha,\partial
S}(x,y)^{p-n}\sigma(Q(x,\rho_{\alpha,\partial S}(x,y)))^2}
\,d\sigma(x)\,d\sigma(y) \right)^{\frac1p}. \nn\ee
The constants of this equivalence depend only on
$n,p,\ve,\alpha,$, the constants $C_\mu,C'_\mu,C''_\mu$ and
$C_\sigma,C'_\sigma,C''_\sigma$ .
\end{theorem}
\par We will be needed the following generalization of the
equivalence \rf{SPACK}.
\begin{proposition}\lbl{NOM-AP} Let $\Omega$ be an open
subset of $\RN$ and let $1<p<\infty$. Then for every
function $f\in \LOPO$ we have
\bel{DM-LL} \|f\|_{L^1_p(\Omega)}\le C
\sup_{\pi}\,\{\sum_{Q\in\pi}r_Q^{n-p}\Ec_1(f;Q)_{L_p}^p\}
^{\frac1p}
\ee
where the supremum is taken over all packings\, $\pi$ of
\underline{equal} cubes lying in $\Omega$. Here the
constant $C$ depends only on $n$ and $p$.
\end{proposition}
\par {\it Proof}. Recall that the quantity
$\Ec_1(f;Q)_{L_p}$ is defined by \rf{ESP}.
\par It suffices to show that, for any packing $\pi'=\{K\}$
of cubes $K\subset\Omega$, there exists a packing
$\pi=\{Q\}$ of {\it equal} cubes $Q\subset\Omega$ such that
\bel{PACK} \sum_{K\in\pi'}\intl_{K}\|\nabla f(x)\|^p\,dx
\le C\sum_{Q\in\pi}r_Q^{n-p}\Ec_1(f;Q)_{L_p}^p. \ee
\par Note that, the proof of equivalence \rf{L1RN} given in
\cite{Z} can be easily extended to the case of an arbitrary
cube $K\subset\RN$. In other words, for every $f\in
L^1_p(int(K)), 1<p<\infty,$ we have
\bel{L1NQ}
\left(\intl_{K}\|\nabla
f(x)\|^p\,dx\right)^{\frac1p} \sim \sup_{0<t\le\diam
K}\omega(f;t)_{L_p(K)}/t
\ee
where $\omega(f;\cdot)_{L_p(K)}$ is the modulus of
smoothness in $L_p(K)$:
$$
 \omega (f;t)_{L_p(K)}:=\sup_{\|h\|<t} \left(
~\intl_{\{x\in K:~x+h\in K\}}|f(x+h)-f(x)|^{p}\,dx\right)
^{\frac{1}{p}}.
$$
\par It is also known that the function
$\omega(f;t)_{L_p(K)}/t$ is quasi-monotone; the reader can
easily prove that
\bel{QMON}
\omega(f;T)_{L_p(K)}/T\le
2\omega(f;t)_{L_p(K)}/t,~~~~0<t<T.
\ee
\par Finally, by Brudnyi's theorem \cite{Br1,Br3},
\bel{MSL} \omega(f;t)_{L_p(K)}^p\sim \sup_{\pi}\,
\sum_{Q\in\pi}|Q|\,\Ec_1(f;Q)_{L_p}^p, ~~~~0<t\le\diam K,
\ee
where $\pi=\{Q\}$ runs over all packings\, $\pi$ of {\it
equal} cubes $Q\subset K$ of diameter $\diam Q=t$.
\par The results cited above easily imply inequality
\rf{PACK}. In fact, in view of \rf{L1NQ}, for every cube
$K\in\pi'$ there exists $t_K\in(0,\diam K]$ such that
$$
\left(\intl_{K}\|\nabla f(x)\|^p\,dx\right)^{\frac1p}\le
C\omega(f;t_K)_{L_p(K)}/t_K
$$
with $C=C(n,p)$. Let $\ot:=\min\{t_K:~K\in\pi'\}$. Then, by
\rf{QMON},
\bel{NLKM} \left(\intl_{K}\|\nabla
f(x)\|^p\,dx\right)^{\frac1p} \le
C\omega(f;\ot)_{L_p(K)}/\ot. \ee
On the other hand, by \rf{MSL}, there exists a packing
$\pi_K=\{Q\}$ of {\it equal} cubes lying in $K$ of diameter
$\diam Q=2r_Q=\ot$ such that
$$
\omega(f;\ot)_{L_p(K)}^p\le C
\sum_{Q\in\pi_K}|Q|\,\Ec_1(f;Q)_{L_p}^p.
$$
This inequality and \rf{NLKM} yield
$$
\intl_{K}\|\nabla f(x)\|^p\,dx \le
(C/\ot^p)
\sum_{Q\in\pi_K}|Q|\Ec_1(f;Q)_{L_p}^p
\le
C\sum_{Q\in\pi_K}r_Q^{n-p}\Ec_1(f;Q)_{L_p}^p.
$$
(Recall that $\ot=2r_Q,~Q\in\pi_K$.)
\par It remains to put $\pi:=\cup\{\pi_K:~K\in\pi'\}$ and
the required inequality \rf{PACK} follows.\bx
\medskip
\par {\bf Proof Theorem \reff{EXT-DF-LA}.} We let
$\|f\|^*_S$ denote the right-hand side of theorem's
equivalence. Then, by Theorem \reff{EXT-DIFF-DSGM} and
Theorem \reff{EXT-DIFF}, $\|f\|^*_S\le C\|f\|_{\WP|_S}$.
\par Prove the converse inequality. Given a function $f\in
C(S)$ and a Borel measure $\nu$ on $\RN$ we put
$$
J(f;t:\nu):=\left(~\iint\limits_{\|x-y\|<t}
\frac{|f(x)-f(y)|^p} {t^{p-n}\nu(Q(x,t))^2}
\,d\nu(x)\,d\nu(y)\right)^{\frac1p},~~~~t>0.
$$
\par Then, by Proposition \reff{FAPN},
$$
J(f;t:\sigma)\le Ct^{-1}
\AP(f;4t:\partial S)_{L_p(\sigma)}
\le Ct^{-1} \AP(f;4t:S), ~~~t\in(0,1/4].
$$
Again, by Proposition \reff{FAPN},
\bel{APJ1} \AP(f;t:S)\le Ct J(f;4t:\mu),~~~
~~~t\in(0,1/16], \ee
so that
\bel{JSJM} J(f;t:\sigma)\le  Ct J(f;16t:\mu),~~~
~~~t\in(0,1/64]. \ee
\par Let $\Omega:=int(S)$. Let us estimate $\|f|\,
_\Omega\|_{W^1_p(\Omega)}$ via $\|f\|_{L_p(\mu)}$ and the
$A_p$-functional of $f$. By Lemma \reff{MUS},
\bel{LMS}
\|f\|_{L_p(\Omega)}\le\|f\|_{L_p(S)} \le
C\{\|f\|_{L_p(\mu)}+\AP(f;\ve:S)\}.
\ee
By Proposition \reff{NOM-AP}, we have
\bel{FEG} \|f\|_{\LOPO}\le
C\left\{\sup_{0<t\le\ve}\AP(f;t:S)/t+
\sup_{\pi}\,\left(\sum_{Q\in\pi}
r_Q^{n-p}\Ec_1(f;Q)_{L_p}^p\right) ^{\frac1p}\right\}
\ee
where the supremum is taken over all packings\, $\pi=\{Q\}$
of equal cubes of diameter $\diam Q=2r_Q\ge\ve$ lying in
$\Omega$. Since
$$
|Q|\Ec_1(f;Q)_{L_p}^p\le
|Q|\left(\frac{1}{Q}\intl_Q |f|^p\,dx\right)
=\intl_Q |f|^p\,dx,
$$
see \rf{ESP}, for every such a packing $\pi$ we have
$$
\sum_{Q\in\pi}r_Q^{n-p}\Ec_1(f;Q)_{L_p}^p\le 2^p\ve^{-p}
\sum_{Q\in\pi}|Q|\Ec_1(f;Q)_{L_p}^p\le C\sum_{Q\in\pi}
\intl_Q |f|^p\,dx\le C\|f\|_{L_p(\Omega)}^p.
$$
Combining this inequality with \rf{FEG} and \rf{LMS}, we
obtain
\be \|f\|_{\LOPO}&\le&
C\left\{\sup_{0<t\le\ve}\AP(f;t:S)/t+
\|f\|_{L_p(\Omega)}\right\}\nn\\ &\le&
C\left\{\sup_{0<t\le\ve}\AP(f;t:S)/t+
\|f\|_{L_p(\mu)}+\AP(f;\ve:S)\right\}\nn\\&\le&
C\left\{\sup_{0<t\le\ve}\AP(f;t:S)/t+
\|f\|_{L_p(\mu)}\right\}.\nn \ee
Hence
$$
\|f|\,_\Omega\|_{W^1_p(\Omega)}=\|f\|_{L_p(\Omega)}+
\|f\|_{\LOPO}\le
C\{\|f\|_{L_p(\mu)}+\sup_{0<t\le\ve}\AP(f;t:S)/t\},
$$
so that, by \rf{APJ1},
$$
\|f|\,_\Omega\|_{W^1_p(\Omega)}\le
C\{\|f\|_{L_p(\mu)}+\sup_{0<t\le\ve}J(f;t:\mu)\},
$$
provided $\ve\in(0,1/64]$. This and inequality \rf{JSJM}
yield
$$
\|f|\,_\Omega\|_{W^1_p(\Omega)}+
\sup_{0<t\le\ve}J(f;t:\sigma)\le
C\{\|f\|_{L_p(\mu)}+
\sup_{0<t\le\ve}J(f;t:\mu)\}\le C\|f\|^*_S.
$$
Finally, by Theorem \reff{EXT-DIFF-DSGM},
$$ \|f\|_{\WP|_S}\le C\{\|f|\,_\Omega\|_{W^1_p(\Omega)}+
\sup_{0<t\le\ve}J(f;t:\sigma)+\|f\|^*_S\}
\le C\|f\|^*_S,$$
proving the theorem.\bx
\medskip
\begin{remark}\lbl{RH-S}
{\em The function $\rho_{\partial
S}$ in the equivalence of Theorem \reff{EXT-DF-LA} can be
replaced with the (bigger) function $\rho_S(x,y)$. In fact,
such a replacement us possible in Theorem
\reff{EXT-DIFF-DSGM}, see Remark \reff{RH-S1}. Since the
proof of Theorem \reff{EXT-DF-LA} relies on the result of
Theorem \reff{EXT-DIFF-DSGM}, the same replacement can be
made in Theorem \reff{EXT-DF-LA} as well.}
\end{remark}
\medskip
\par {\bf Proof of Theorem \reff{EXT-CB-DM} and Theorem
\reff{EXT-THIN}.}
\par We let $\Bc(\RN)$ denote the class of all closed
subsets of $\RN$ satisfying the ball condition, see
Definition \reff {BALLCND}. It will be convenient for us to
use the following equivalent definition of the class
$\Bc(\RN)$:
\par We say that a {\it closed set $A\in\Bc(\RN)$ if there
exists a constant $\beta_A>0$ such that the following is
true: For every cube $Q$ of diameter at most $1$ centered
in $A$ there exists a cube $Q'\subset Q\setminus A$ such
that $\diam Q \le \beta_A\diam Q'.$}
\par The notion of a set satisfying ``the ball condition"
is well known in the literature and has been used in many
papers (in equivalent forms and with different names), see,
e.g. \cite{HM,JK,MV}. Observe that a set $A\in\Bc(\RN)$
satisfies the ball condition if and only if every cube $Q$
of $\diam Q\le 1$ centered in $A$ is $1/\beta_A$-porous
(with respect to $A$), see Definition \reff{PR-Q}. Also, it
can be readily seen that for every set $A\in\Bc(\RN)$ and
every $\alpha\in(0,\frac{1}{2\beta_A}]$ the following
inequality
\bel{RXY} \|x-y\|\le \rrxy{A}\le 4\|x-y\|, ~~~~x,y\in
A,~\|x-y\|\le 1/4, \ee
holds. (See \rf{DXYA} for the definition of
$\rho_{\alpha,A}$.)
\par The next theorem presents a slight generalization of
Theorem \reff{EXT-THIN}.
\begin{theorem}\lbl{EXT-THIN-EP} Let $p\in(n,\infty)$,
$\ve>0$, and let $S$ be a closed subset of $\RN$. Assume
that $int(S)=\emp$ and the ball condition is satisfied.
Then for every $f\in C(S)$
$$
\|f\|_{\WP|_S}\sim
\|f\|_{L_p(\mu)}+\left(\,\,\iint\limits_{\|x-y\|<\ve}
\frac{|f(x)-f(y)|^p} {\|x-y\|^{p-n}\mu(Q(x,\|x-y\|))^2}
\,d\mu(x)d\mu(y)\right)^{\frac1p}
$$
with constants of equivalence depending only on
$n,p,\ve,C_\mu,C'_\mu,C''_\mu,$ and $\beta_{S}$.
\end{theorem}
\par {\it Proof.} Since the interior of $S$ is empty,
$S=\partial S$ so that one can put $\mu=\sigma$. Since $S$
satisfies the ball condition, every cube $Q$ of $\diam Q\le
1$ centered in $S$ is $\alpha:=1/\beta_S$-porous so that,
by Definition \reff{DEF-AP} and Definition \reff{PR-AP},
$$
\AP(f;t:S)=\PRM(f;t:S),~~~~t\in(0,1].
$$
Thus the $A_p$-functional is majorized by the
$\Ac_{p,\alpha}-functional$ which allows us in the proof of
Theorem \reff{EXT-MU-NORM} to omit all estimates related to
the $A_p$-functional. This leads us to a variant of Theorem
\reff{EXT-MU-NORM} where the expression related to  $A_p$
is omitted as well; thus for every $f\in C(S)$ we have
$$ \|f\|_{\WP|_S}\sim \|f\|_{L_p(\mu)}+
+\left\{\intl_0^\ve\PRM(f;t:S)_{L_p(\mu)}^p\,
\frac{dt}{t^{p+1}}\right\}^{\frac{1}{p}} $$
provided $p\in(n,\infty),$ $0<\alpha<\min\{1/7,1/\beta_S)$,
and $\ve>0$.
\par We proof Theorem \reff{EXT-DIFF-DSGM} basing on this
version of Theorem \reff{EXT-MU-NORM}. This leads us to a
variant of Theorem \reff{EXT-DIFF-DSGM} (with $\sigma=\mu$)
where the third term in the theorem's equivalence can be
omitted. Recall that $int(S)=\emp$ so that the first term
can be omitted as well. We obtain that for every
$p\in(n,\infty)$, $\ve>0$,
$0<\alpha<\min\{1/14,1/\beta_S)$, and every $f\in C(S)$
$$ \|f\|_{\WP|_S}\sim \|f\|_{L_p(\mu)}+
\left(~\iint\limits_{\rho_{\alpha,S}(x,y)<\ve}
\frac{|f(x)-f(y)|^p} {\rho_{\alpha,S}(x,y)^{p-n}
\mu(Q(x,\rho_{\alpha,S}(x,y)))^2}
\,d\mu(x)\,d\mu(y) \right)^{\frac1p}. $$
\par Finally, by \rf{RXY}, $\rho_{\alpha,S}(x,y)\sim
\|x-y\|$ which allows us to replace in this equivalence
$\rho_{\alpha,S}(x,y)$ with $\|x-y\|$. The theorem is
proved.\bx
\medskip
\par This theorem and Proposition \reff{WDS} imply the
following generalization of Theorem \reff{EXT-CB-DM}.
\begin{theorem}\lbl{EXT-CB-CND} Let $p\in(n,\infty)$,
$\ve>0$, $\Omega:=int(S)$, and let $\sigma$ be a measure
with $\supp\sigma=\partial S$ satisfying conditions
\rf{DNSG} and \rf{B1SIGMA}. Suppose that $\partial S$
satisfies the ball condition. Then for every $f\in C(S)$
the following equivalence
\be \|f\|_{\WP|_S}&\sim&
\|f|\,_\Omega\|_{W^1_p(\Omega)}+\|f\|_{L_p(\sigma)}\nn\\
&+&\left(~\iint\limits_{\|x-y\|<\ve} \frac{|f(x)-f(y)|^p}
{\|x-y\|^{p-n}\sigma(Q(x,\|x-y\|))^2}
\,d\sigma(x)\,d\sigma(y) \right)^{\frac1p}\nn\ee
holds with constants depending only on
$n,p,\ve,\beta_{\partial S},$ and
$C_\sigma,C'_\sigma,C''_\sigma,$.
\end{theorem}
\par {\it Proof.} By Proposition \reff{WDS},
$$ \|f\|_{\WP|_S}\sim \|f|_\Omega\|_{\WPO}+ \|f|_{\partial
S}\|_{\WP|_{\partial S}} $$
provided $p>n$ and $f\in C(S)$. Since $\partial S$
satisfies the ball condition and $int(\partial S)$ is
empty, one can apply Theorem \reff{EXT-THIN-EP} to
$\partial S$ with $\mu=\sigma$. By this theorem, the
quantity $\|f|_{\partial S}\|_{\WP|_{\partial S}}$ is
equivalent to the sum of the last two terms in the
equivalence of Theorem \reff{EXT-CB-CND}. The proof is
finished.\bx
\medskip
\par Let us present two examples of description of the
traces of Sobolev functions.
\begin{example}\lbl{EXT-REG} {\em \par Let us apply the
results of Theorems \reff{EXT-DIFF} and \reff{EXT-CB-CND}
to $n$-sets in $\RN$, see \rf{DEF-D}. In this case the
measure $\mu$ in \rf{DEF-D} is the restriction to $S$ of
the Lebesgue measure. Thus $S$ is an $n$-set if there is a
constant $\theta=\theta_S\ge 1$ such that for every cube
$Q$ centered in $S$ of diameter at most $1$ we have $
|Q|\le \theta_S |Q\cap S|$. We call $n$-sets {\it regular}
subsets of $\RN$.
\par In \cite{S3,S4} we described the traces of Sobolev
functions to regular sets via the sharp maximal function
$$
f^{\sharp}_{1,S}(x):=\sup_{r>0}\frac{1}{r^{n+1}}
\intl_{Q(x,r)\cap S}|f-f_{Q(x,r)\cap S}|\,dx,~~~x\in S.
$$
We have proved that, for every regular set $S$ and every
function $f\in L_p(S)$, $p>1,$
$$
\|f\|_{\WP|_S}\sim
\|f\|_{L_p(S)}+\|f_{1,S}^{\sharp}\|_{L_p(S)}
$$
with constants of equivalence depending only on $n,p$ and
$\theta_S$.
\par Theorem \reff{EXT-DIFF} provides another description
of the trace space $\WP|_S$ which does not use the notion
of the sharp maximal function: Let $p\in(n,\infty)$,
$\alpha\in(0,1/14)$ and $f\in C(S)$. Then
\be \|f\|_{\WP|_S}&\sim& \|f\|_{L_p(S)}+ \sup_{0<t\le 1}
\left(~\iint\limits_{\|x-y\|<t}
\frac{|f(x)-f(y)|^p}{t^{p+n}}
\,dx\,dy\right)^{\frac1p}\nn\\
&+& \left(~\iint\limits_{\rxy<1} \frac{|f(x)-f(y)|^p}
{\rxy^{p+n}}\,dx\,dy \right)^{\frac1p}, \nn\ee
with constants of equivalence depending only on $n,p,$ and
$\theta_S$.
}\end{example}
\begin{example}\lbl{EXND} {\em \par Let us present an
example of a set in $S\subset\R^2$ satisfying the ball
condition such that the one dimensional Hausdorff measure
$\mu:=\Hc_1|S$ is doubling, but $S$ is not a $d$-set for
any $0<d<2$.
\par By $\ell_i,~i=0,1,...,$, we denote a line segment in
$\R^2$,
$$
\ell_i:=\{x=(x_1,x_2):~ x_1=2^{-i},~0\le x_2\le 4^{-i}\},
$$
and put
$ S:=\{0\}\cup(\cup_{i=0}^\infty\ell_i).$
\par Simple calculations show that for every square
$Q=Q(x,r)$, $x\in S, 0<r\le 1$, we have
$$
\mu(Q(x,r)):=\Hc_1(Q(x,r)\cap S)\sim \left \{
\begin{array}{lll}
r^2,& \|x\|\le r,\\\\
\|x\|^2,& \|x\|^2\le r<\|x\|,\\\\
r,& 0<r<\|x\|^2.
\end{array}
\right.
$$
Using this equivalence, one can easily see that $\mu$ is
doubling.
\par Clearly, $S$ satisfies the ball condition. It is also
obvious that $S$ is not a $d$-set for any $d$. Let us
calculate the trace norm $\|f\|_{W^1_p(\R^2)|_S}$.
\par It can be easily seen that for every $x=(x_1,x_2)$ and
$y=(y_1,y_2)\in S$ such that $\|x-y\|\le 1$ we have
$$
\mu(Q(x,\|x-y\|))\sim \left \{
\begin{array}{ll}
\|x-y\|,& x_1=y_1,\\\\
\|x-y\|^2,& x_1\ne y_1.
\end{array}
\right.
$$
Hence, by Theorem \reff{EXT-CB-CND},
\be \|f\|_{\WP|_S}&\sim& \|f\|_{L_p(\mu)}
+\left(~\iint\limits_{\{x,y\in S:\,x_1=y_1\}}
\frac{|f(x)-f(y)|^p} {\|x-y\|^{p}}\,d\mu(x)\,d\mu(y)
\right)^{\frac1p}\nn\\&+& \left(~\iint\limits_{\{x,y\in
S:\,x_1\ne y_1\}} \frac{|f(x)-f(y)|^p}
{\|x-y\|^{p+2}}\,d\mu(x)\,d\mu(y) \right)^{\frac1p}.
\nn \ee

\par Observe that, $d\mu(x)=dx_2$ (because the segments
$\ell_i$ are parallel to the axis $Ox_2$ and $\mu$ is the
one dimensional Hausdorff measure), $\|x-y\|=|x_2-y_2|$
whenever $x_1=y_1$, and $\|x-y\|=|x_1-y_1|$ whenever
$x_1\ne y_1$. This enables us to present the above
equivalence in the following form:
\be \|f\|_{\WP|_S}&\sim& \|f\|_{L_p(\mu)}
+\left(~\iint\limits_{\{x,y\in S:\,x_1=y_1\}}
\frac{|f(x_1,x_2)-f(x_1,y_2)|^p}
{|x_2-y_2|^{p}}\,dx_2\,dy_2 \right)^{\frac1p}\nn\\&+&
\left(~\iint\limits_{\{x,y\in S:\,x_1\ne y_1\}}
\frac{|f(x_1,x_2)-f(x_1,y_2)|^p}
{|x_1-y_1|^{p+2}}\,dx_2\,dy_2 \right)^{\frac1p}. \nn \ee
}\end{example}

\vspace*{10mm}
Department of Mathematics\\
Technion - Israel Institute of Technology\\
32000 Haifa\\
Israel\\\\
e-mail: pshv@tx.technion.ac.il\\\\
June 15, 2008

\begin{thebibliography}{ABCD}
\bibitem [A] {A} D. R. Adams, Sobolev spaces. New York,
    Academic Press, 1975.
\bibitem [BL] {BL} B. Bennewitz, J.Lewis, On weak reverse
    H\"{o}lder inequalities for nondoubling harmonic
    measures, Complex Var. Theory Appl. 49 (2004), no. 7-9,
    571--582.
\bibitem [BIN] {BIN} O. V. Besov, V. P. Il'in, S. M.
    Nikol'ski,  Integral Representations of Functions and
    Embedding Theorems, Nauka, Moskow, 1975; English
    edition: Winston and Sons, Washington D.C., vol. I.
    1978, vol. II, 1979.
\bibitem [BK] {BK} O. Besov, G. Kalyabin, Spaces of
    differentiable functions. Function spa\-ces,
    dif\-fe\-rential operators and nonlinear analysis
    (Teistungen, 2001), 3--21, Birkh\"auser, Basel, 2003.
\bibitem [BMP]{BMP} E. Bierstone, P. Milman and W.
    Pawlucki, Differentiable functions defined in closed
    sets. A problem of Whitney, Invent. Math.  151 (2003),
    no. 2, 329--352.
\bibitem [BM]{BM} E. Bierstone and P. Milman, $\Cc^m$
    norms on finite sets and $\Cc^m$ extension criteria,
    Duke Math. J. 137 (2007) 1--18.
\bibitem [Br1]{Br1} Yu. A. Brudnyi, Spaces that are
    definable by means of local approximations, Trudy
    Moscov. Math. Obshch., 24 (1971) 69--132; English
    Transl.: Trans. Moscow Math. Soc., 24 (1974) 73--139.
\bibitem [Br2]{Br2} Yu. A. Brudnyi, Investigation in the
    theory of local approximations, Doctoral Dissertation,
    Leningrad, 1977 (Russian).
\bibitem [Br3]{Br3} Yu. A. Brudnyi, Adaptive approximation
    of functions with singularities, Trudy Moscov. Mat.
    Obshch. {\bf 55,}  149--242 (1994); English transl.:
    Trans. Moscow Math. Soc. {\bf 55,} 123--186 (1995).
\bibitem [BrK]{BrK} Yu. A. Brudnyi, B.D. Kotljar, A certain
    problem of combinatorial geometry,  Sibirsk. Mat. Z. 11
    (1970) 1171–-1173; English transl. in Siberian Math. J.
    11 (1970), 870–871.
\bibitem [BS1]{BS1} Yu. Brudnyi and P. Shvartsman,
    Generalizations  of Whitney's  Extension Theorem,
    Intern. Math.  Research Notices (1994), no. 3,
    129--139.
\bibitem [BS2]{BS2} Yu. Brudnyi and P. Shvartsman, The
    Whitney Problem  of  Existence of a Linear Extension
    Operator,  J. Geom. Anal. 7, no. 4 (1997) 515--574.
\bibitem [BS3]{BS3} Yu. Brudnyi and P. Shvartsman,
    Whitney's Extension Problem for Multivariate
    $C^{1,\omega}$-functions, Trans. Amer. Math. Soc. 353
    (2001), no. 6, 2487--2512.
\bibitem [Bur]{Bur} V. I. Burenkov, Sobolev spaces on
    domains, Teubner-Texte zur Mathematik, 137. B. G.
    Teubner Verlagsgesellschaft, Stuttgart-Leipzig, 1998.
\bibitem [C1]{C1} A. P. Calder\'{o}n, Estimates for
    singular integral operators in terms of maximal
    functions, Studia Math. 44 (1972) 563--582.
\bibitem [CS]{CS}  A. P. Calder\'{o}n, R. Scott, Sobolev
    type inequalities for $p>0$,  Studia Math. 62 (1978)
    75--92.
\bibitem [DL]{DL} R. A. DeVore, G.G. Lorentz, Constructive
 approximation, Springer, New-York, 1993.
\bibitem [DS1]{DS1} R. A. DeVore,  R. C. Sharpley, Maximal
 functions measuring smoothness, Mem. Amer. Math. Soc.
 47 (1984), no. 293, viii+115 pp.
\bibitem [D1]{D1} E. M. Dynkin, Free interpolation by
    functions with derivatives in $H^1$, Zap. Nauchn. Sem.
    Leningrad Otdel. Mat. Inst. Steklov (LOMI) 126 (1983),
    77-87; English transl. in J. Soviet Math. 27 (1984),
    2475--2481.
\bibitem [D2]{D2} E. M. Dynkin, Homogeneous measures on
    subsets of $\RN$, Lect. Notes in Math., vol. 1043,
    Springer-Verlag, 1983, 698--700.
\bibitem [Dol]{Dol} V. L. Dolnikov, The partitioning of
    families of convex bodies, Sibirsk. Mat. \v{Z},
    12 (1971), 664–667 (Russian); English transl. in
    Siberian Math. J. 12 (1971), 473–475.
\bibitem [FJ]{FJ} W. Farkas,   N. Jacob, Sobolev spaces on
    Non-Smooth Domains and Dirichlet Forms Related to
    Subordinate Reflecting Diffusions, Math. Nachr.  224
    (2001) 75--104.
\bibitem [F1]{F1} C. Fefferman,  A sharp form of Whitney's
    extension theorem, Ann. of Math. (2), 161
    (2005), no. 1, 509--577.
\bibitem [F2]{F2} C. Fefferman,  A Generalized Sharp
    Whitney Theorem for Jets, Rev. Mat.
    Iberoamericana 21 (2005), no. 2, 577--688.
\bibitem [F3]{F3} C. Fefferman,  Whitney's Extension
    Problem for $C^m$, Ann. of Math. (2), 164 (2006),
    no. 1,  313-359.
\bibitem [F4]{F4} C. Fefferman,  $C^m$ Extension by Linear
    Operators, Ann. of Math. (2) 166 (2007), no. 3,
    779--835.
\bibitem [GV1]{GV1}  V. M. Gol'dstein , S. K.
    Vodop'janov,  Prolongement des fonctions de classe
    ${\cal L}\sp{1}\sb{p}$ et applications quasi
    conformes.,  C. R. Acad. Sci. Paris Ser. A-B 290
    (1980), no. 10, A453--A456.
\bibitem [GV2]{GV2}  V. M. Gol'dstein and S. K.
    Vodop'janov,  Prolongement des fonctions
    differentiables hors de domaines plans., C. R. Acad.
    Sci. Paris Ser. I Math. 293 (1981), no. 12, 581--584.
\bibitem [GK]{GK} L. Grafakos, J. Kinnunen, Sharp
    inequalities for maximal functions associated with
    general measures, Proc. Roy. Soc. Edinburgh Sect. A.,
    128A (1998) 717–-723.
\bibitem [Gud]{Gud} J. Gudayol,  Extension theorems of
 Whitney type by the use of integral operators,
arXiv:math/9803016v2 [math.CA] 23 Apr 1998.
\bibitem [G]{G} M. de Guzm\'{a}n, Differentiation of
    integrals in $\RN$, Lect. Notes in Math. 481,
    Springer-Verlag, 1975.
\bibitem [HM]{HM}  P. Haj\l asz,  O. Martio, Traces of
    Sobolev functions on fractal type sets and
    characterization of extension domains, J. Funct. Anal.
    143 (1997) 221--246.
\bibitem [JK]{JK} D. S. Jerison, C. E. Kenig, Boundary
    behavior of harmonic functions in non-tangentially
    accessible domains, Adv. in Math., 46 (1982), 80--147.
\bibitem [Jn]{Jn} P. W. Jones, Quasiconformal mappings and
    extendability of functions in Sobolev spaces, Acta
    Math., 147 (1981), 71--78.
\bibitem [J1]{J1}  A. Jonsson, The trace of potentials on
    general sets, Ark. Mat., 17 (1979), 1--18.
\bibitem [JW]{JW} A. Jonsson, H. Wallin, Function Spaces on
    Subsets of $\RN$, Harwood Acad. Publ., London, 1984,
    Mathematical Reports, Volume 2, Part 1.
\bibitem [J]{J}  A. Jonsson, Besov spaces on closed subsets
    of $\RN$,  Trans.  Amer. Math. Soc., 341 (1994)
    355--370.
\bibitem [LS]{LS}  J. Luukkainen and E. Saksman, Every
    complete doubling metric space carries a doubling
    measure, Proc.  Amer. Math. Sov., 126, no. 2, (1998),
    531--534.
\bibitem [K]{K} P. Koskela, Extensions and Imbeddings,
    J. Funct. Anal., 159 (1998) 369--383.
\bibitem [MV]{MV} O. Martio, J. V\"{a}is\"{a}l\"{a}, Global
    $L_p$-integrability of the derivative of quasiconformal
    mapping, Complex Variables Theory Appl. 9(1988), no. 4,
    309--319.
\bibitem [M]{M}  V.G. Maz'ja,  Sobolev spaces,
     Springer-Verlag, Berlin, 1985, xix+486 pp.
\bibitem [MP]{MP} V. Maz'ya,  S. Poborchi, Differentiable
    Functions on Bad Domains, Word Scientific, River Edge,
    NJ, 1997.
\bibitem [N]{N} S. M. Nikol'skii, Approximation of
 functions of several variables and imbedding theorems, Die
 Grundlehren  der Mathematischen Wissenschaften, Band 205.
Springer-Verlag, New York, Heidelberg, 1975. xviii+418 pp.
\bibitem [S3]{S3} P. Shvartsman, On extension of Sobolev
    functions defined on regular subsets of metric measure
    spaces, J. Approx. Theory, 144 (2007), 139–-161.
\bibitem [S4]{S4} P. Shvartsman, Local approximations and
    intrinsic characterizations of spaces of smooth
    functions on regular subsets of $\RN$, Math. Nachr. 279
    (2006), no.11, 1212–-1241.
\bibitem [S5]{S5} P. Shvartsman, A Whitney-type extension
 theorem for $W^k_p(\RN)$-functions, $p>n$. (to appear)
\bibitem [S6]{S6} P. Shvartsman, The Whitney extension
    problem and Lipschitz selections of set-valued mappings
    in jet-spaces, Trans. Amer. Math. Soc. (to appear)
\bibitem [St]{St} E. M. Stein, Singular integrals and
    differentiability properties of functions, Princeton
    Univ. Press, Princeton, New Jersey, 1970.
\bibitem [T1]{T1} H. Triebel, Local approximation spaces,
    Z. Analysis Anwendungen 8 (1989), 261--288.
\bibitem [T2]{T2} H. Triebel, Theory of function spaces.
    II. Monographs in Mathematics, 84. Birkh\"{a}user
    Verlag, Basel, 1992.
\bibitem [T3]{T3} H. Triebel, The structure of functions,
    Monographs in Mathematics, 97. Birkh\"{a}user Verlag,
    Basel, 2001.
\bibitem [VK]{VK} A. L. Vol'berg, S. V. Konyagin, On
    measures with doubling condition,  Math. USSR Izvestiya
    30 (1988) 629--638.
\bibitem [W1]{W1}  H. Whitney, Analytic extension of
 differentiable functions defined in closed sets,
 Trans.  Amer. Math. Soc. 36 (1934), 63--89.
\bibitem [W2]{W2} H. Whitney, Differentiable functions
 defined in closed sets. I., Trans. Amer. Math.
Soc. 36 (1934) 369--387.
\bibitem [Wu]{Wu} J-M. Wu, Hausdorff dimension and doubling
    measures on metric spaces, Proc. Amer. Math. Soc. 126
 (1998), 1453--1459.
\bibitem [Z]{Z} W. Ziemer, Weakly differentiable functions,
    Graduate Text in Mathematics, Springer-Verlag,
    New-York, 1989.
\bibitem [Zob1]{Zob1} N. Zobin, Whitney's problem on
    extendability of functions and an intrinsic metric,
    {\it Advances in Math.} {\bf 133} (1998) 96--132.
\end{thebibliography}
\end{document}